\documentclass[11pt]{amsart}

\usepackage{hyperref}
\hypersetup{nesting=true,debug=true,naturalnames=true}
\usepackage{graphicx,amssymb,upref}
\usepackage[usenames,dvipsnames]{pstricks}
\usepackage{epsfig}
\usepackage{color}
\usepackage[latin1]{inputenc}
\usepackage[T1]{fontenc}
\usepackage{amsmath, amsthm}
\usepackage[english]{babel}
\usepackage{amssymb}
\usepackage{latexsym}
\usepackage{graphicx,amssymb,upref}
\usepackage[all]{xy}
\usepackage{tikz}
\usetikzlibrary{matrix}

\date{\today}

\newtheorem{theorem}{Theorem}[section]
\newtheorem{lemma}[theorem]{Lemma}

\newtheorem{corollary}[theorem]{Corollary}
\theoremstyle{definition}
\newtheorem{definition}[theorem]{Definition}
\newtheorem{example}[theorem]{Example}

\newtheorem{remark}[theorem]{Remark}
\newtheorem{notation}[theorem]{Notation}

\newcommand{\ot}{\otimes}
\newcommand{\co}{\circ}

\newcommand{\ID}{I(q_D)}

\title[Projections of Hopf braces]{Projections of Hopf braces}
\title{Projections of Hopf braces} 

\begin{document}
	
\maketitle
	
	\begin{center}
	{\bf Jos\'e Manuel Fern\'andez Vilaboa$^{1}$, Ram\'on
	Gonz\'{a}lez Rodr\'{\i}guez$^{2}$, Brais Ramos P\'erez$^{3}$ and Ana Bel\'en Rodr\'{\i}guez Raposo$^{4}$}.
	\end{center}
	
	\vspace{0.4cm}
	\begin{center}
	{\small $^{1}$ [https://orcid.org/0000-0002-5995-7961].}
	\end{center}
	\begin{center}
	{\small  CITMAga, 15782 Santiago de Compostela, Spain.}
	\end{center}
	\begin{center}
	{\small  Universidade de Santiago de Compostela. Departamento de Matem\'aticas,  Facultade de Matem\'aticas, E-15771 Santiago de Compostela, Spain. 
	\\ email: josemanuel.fernandez@usc.es.}
	\end{center}
	\vspace{0.2cm}
	
	\begin{center}
	{\small $^{2}$ [https://orcid.org/0000-0003-3061-6685].}
	\end{center}
	\begin{center}
	{\small  CITMAga, 15782 Santiago de Compostela, Spain.}
	\end{center}
	\begin{center}
	{\small  Universidade de Vigo, Departamento de Matem\'{a}tica Aplicada II,  E. E. Telecomunicaci\'on,
	E-36310  Vigo, Spain.
	\\email: rgon@dma.uvigo.es}
	\end{center}
    \vspace{0.2cm}
    
    \begin{center}
   	{\small $^{3}$ [https://orcid.org/0009-0006-3912-4483].}
   \end{center}
	\begin{center}
	{\small  CITMAga, 15782 Santiago de Compostela, Spain. \\}
	\end{center}
    \begin{center}
	{\small  Universidade de Santiago de Compostela. Departamento de Matem\'aticas,  Facultade de Matem\'aticas, E-15771 Santiago de Compostela, Spain. 
	\\email: braisramos.perez@usc.es}
	\end{center}
	\vspace{0.2cm}
	
	\begin{center}
	{\small $^{4}$ [https://orcid.org/0000-0002-8719-5159]}
	\end{center}
	\begin{center}
	{\small Universidade de Santiago de Compostela, Departamento de Did\'acticas Aplicadas, Facultade C. C. Educaci\'on, E-15782 Santiago de Compostela, Spain.
	\\email: anabelen.rodriguez.raposo@usc.es}
	\end{center}

%	\begin{center}
%	{\small $^{d}$ Corresponding author}
%	\end{center}

\begin{abstract}  
This paper is devoted to the study of Hopf braces projections  in a monoidal setting. Given a cocommutative Hopf brace ${\mathbb H}$ in a strict symmetric monoidal category {\sf C}, we define the braided monoidal category of left Yetter-Drinfeld modules over ${\mathbb H}$. For a Hopf brace ${\mathbb A}$ in  this category, we introduce the concept of bosonizable Hopf brace and we prove that its bosonization ${\mathbb A}\blacktriangleright\hspace{-0.1cm}\blacktriangleleft {\mathbb H}$ is a new Hopf brace in {\sf C} that gives rise to a projection of Hopf braces satisfying certain properties. Finally, taking these properties into account, we introduce the notions of v$_{i}$-strong projection over ${\mathbb H}$, $i=1,2,3,4$,  and we prove that there is a categorical equivalence between the categories of bosonizable Hopf braces in the category of left Yetter-Drinfeld modules over ${\mathbb H}$ and the category of v$_{4}$-strong projections over ${\mathbb H}$.
\end{abstract}

\vspace{0.2cm}



\section*{Introduction}

The study of non-degenerate set-theoretical  solutions of the Yang-Baxter equation with the involutive property is the origin of the notion of brace introduced by W. Rump in \cite{Rump}. A brace is  an abelian group $(G, .)$ with another group structure $(G, \star)$  satisfying for all $g,h,t\in G$ the following condition:
$$g\star (h . t) \star g = (g\star h) . (g\star t)$$

In \cite{CEDO} we can find an equivalent notion of brace and, taking inspiration from it,  recently, L. Guarnieri and L. Vendramin introduced  in \cite{GV} a generalization of braces, called skew braces, as a tool to find non-degenerate bijective solutions of the Yang-Baxter equation not necessarily involutive. Following the definition of L. Guarneri and L. Vendramin, a  skew brace is  a group $(G, .)$ with an additional group structure $(G, \star)$  satisfying 
$$g\star (h . t) = (g\star h) .g^{-1}. (g\star t),$$
for all $g,h,t\in G,$
and it is easy to see that Rump's  braces are examples of skew braces.

In this way, the latest extension of the notion of brace was proposed by I. Angiono, C. Galindo and L. Vendramin in \cite{AGV} with the name of Hopf braces. Hopf braces are the quantum version of skew braces, provide solutions of the Yang-Baxter equation and, as was pointed by the authors, give the right setting for considering left symmetric algebras as Lie-theoretical analogs of braces. If $(H,\epsilon, \delta)$ is a coalgebra,  a Hopf brace structure over $H$ consist on the following data: A Hopf algebra structure $$H_{1}=(H,1,\cdot,\epsilon, \delta,\lambda), $$ 
and  a Hopf algebra structure $$H_{2}=(H,1_\circ,\circ,\epsilon, \delta,s)$$
satisfying the following compatibility:
$$
g\circ(h.k)=(g_1\circ h).\lambda(g_2).(g_3\circ k),\;\; g,h,k\in H.
$$

In any Hopf brace, $1_\circ=1$ and, in this introduction, we will denote a Hopf brace by ${\mathbb  H}=(H_{1}, H_{2})$ or also, in a more reduced form, as ${\mathbb  H}$.

Bearing in mind that the notion of Hopf brace is closely linked to that of Hopf algebra, recently,  A. Agore proposed in \cite{AGORE} a  method  to construct new examples of  Hopf braces working with matched pairs of Hopf algebras $(A, H,  \triangleright, \triangleleft)$ where $H$ is cocommutative. Finally, as has been proved in \cite{AGV} (see also \cite{GGV} and \cite{GONROD}) there exists a strong connection between Hopf braces and invertible $1$-cocycles that induces a categorical equivalence between the categories of Hopf braces and bijective 1-cocycles. 

On the other hand, in the category of vector spaces over a field ${\mathbb K}$, a well known result by D. E. Radford gives the  conditions for the tensor product of two Hopf algebras $Z\otimes
X$ (equipped with smash product algebra and smash coproduct coalgebra) to be a Hopf algebra, and characterizes such objects
via bialgebra projections (see \cite{RAD}). S. Majid in \cite{MAJ2}
interpreted this result in the modern context of braided
categories and stated that there is a categorical equivalence between the category of  Hopf algebras in the category of
left-left Yetter-Drinfeld modules  over $X$ and the category of  Hopf algebra
projections for $X$.  The concrete details of this equivalence are the following: Let  $X$ be a Hopf algebra  and let $(X,Y,f,h)$ be a Hopf algebra
projection over $X$, i.e., $Y$ is a Hopf algebra, $f:X\rightarrow Y$ and $h:Y\rightarrow X$ are morphisms of  Hopf algebras and the following identity holds $h\circ f=id_{X}$. 
Let $I(q_{Y})$ be the image of the idempotent morphism $q_{Y}:Y\rightarrow Y$ defined by the convolution product of the identity of $Y$ and the composition $f\co \lambda_{X}\co h$ where $\lambda_{X}$ is the antipode of $X$. Then, the object $I(q_{Y})$ (the algebra of coinvariants) is a Hopf algebra in the category of left-left Yetter-Drinfeld modules  over $X$ denoted by  $\;^{X}_{X}{\sf Y}
{\sf D}$. Conversely,
if $A$ is a Hopf algebra in $\;^{X}_{X}{\sf Y}{\sf D}$, let $Y=A\blacktriangleright\hspace{-0.1cm}\blacktriangleleft X$ be the
smash (co)product (co)algebra, i.e., $Y$ is the bosonization of $A$ (see
Proposition 4.15 in \cite{MAJ1}). Then $Y=A\blacktriangleright\hspace{-0.1cm}\blacktriangleleft X$ is  a Hopf algebra  and, if $1_{A}$ is the unit of $A$ and $\varepsilon_{A}$ its counit,  $f:X\rightarrow Y$, $f(x)=1_{A}\otimes x$,  and $h:Y\rightarrow X$, $h(a\otimes x)=\varepsilon_{A}(a)x$, are morphisms of Hopf algebras such
that $h\circ f=id_{X}$. These constructions are mutually inverse in the following way: For any  Hopf algebra
projection   $(X,Y,f,h)$,  there exists an isomorphism of Hopf algebras between $Y$ and $I(q_{Y})\blacktriangleright\hspace{-0.1cm}\blacktriangleleft X$ and, for any Hopf algebra $A$ in $\;^{X}_{X}{\sf Y}{\sf D}$, $A=I(q_{A\blacktriangleright\hspace{-0.1cm}\blacktriangleleft X})$. Later, Bespalov proved the  same
results for braided categories with split idempotents in \cite{BES1} and, in collaboration with Drabant, he continued the development of Radford and Majid  theory in this setting (see \cite{BESPA1}, \cite{BESPA2} and \cite{BESPA3}). 

In  \cite{BN}, D. Bulacu and E. Nauwelaerts explained in detail how the
above ideas can be generalized to quasi-Hopf algebras, and in \cite{NikaRamon4}, J. N. Alonso \'Alvarez, J. M. Fern\'{a}ndez Vilaboa and  R. Gonz\'{a}lez Rodr\'{\i}guez obtain a similar categorical equivalence  for weak Hopf algebras in a braided monoidal  setting. Continuing in this line of generalization, the study of projections of Hopf braces begins in the work of H. Zhu in \cite{Zhu} where a method  to build Hopf braces is given based on the new notion of left-compatible ${\mathbb H}$-module. Following the work of H. Zhu, if ${\mathbb  H}$ is a Hopf brace, a left $H_{1}$-module $(M,\triangleright)$ is called a left module over the  Hopf brace ${\mathbb  H}$ if $(M,\blacktriangleright)$  is a left $H_{2}$-module and the following identities 
\begin{itemize}
\item[(i)] $g\blacktriangleright (h\triangleright m)=[(g_{1}\circ h)\cdot \lambda(g_{2})]\triangleright (g_{3}\blacktriangleright m),$
\item[(ii)] $g_1\blacktriangleright m\otimes g_2= (g_{1}\cdot \lambda(g_{3}))]\triangleright 
(g_4\blacktriangleright m) \otimes g_2,$
\end{itemize}
hold for all $g,h\in H$ and $m\in M$. 

It is a relevant fact that the  condition (ii) is used by H. Zhu  to prove that the category of left modules over the Hopf brace ${\mathbb  H}$  is monoidal and, if  $H$ is cocommutative, (ii) always hold. However, this condition presents one  problem: In a general context  the trivial object $(H, \triangleright =\cdot, \blacktriangleright=\circ)$  it is not an example of left module over the  Hopf brace ${\mathbb  H}$. 

Taking into account the above, in \cite{Zhu} the author also introduce the definitions of subbialgebra and left compatible module over a Hopf brace.   These notions are the foundations that support the definition of left Yetter-Drinfeld module for a Hopf brace introduced in \cite[Definition 4.7]{Zhu} and also an analogue of Radford's result about Hopf algebras (see  Remark \ref{timo}). Subsequently in \cite{Zhu2}, H. Zhu and Z. Ying expanded the study of the projection problem for Hopf braces introducing the notion of compatible Hopf brace: Roughly speaking, if $H$ is a Hopf algebra with bijective antipode, a Hopf brace ${\mathbb R}$ in the category of left-left Yetter-Drinfeld modules over $H$ is called compatible if  ${\mathbb R}\otimes H$ equipped with smash product algebra and smash coproduct coalgebra is a Hopf brace. Then, the main result proved in \cite{Zhu2} asserts the following: Let $H$ be a Hopf algebra with bijective antipode and let ${\mathbb A}=(A_{1}, A_{2})$ be a Hopf brace with a projection $(H,{\mathbb A},f,h)$ such that $f(g)\cdot a=f(g)\circ a$ for all $g\in H$ and $a\in A$. Then, there exists a compatible Hopf brace ${\mathbb R}$ such that ${\mathbb A}$ is isomorphic to the smash product algebra and smash coproduct coalgebra  of ${\mathbb R}$ with $H$ as Hopf braces. However, as it is proved in Example \ref{ex-YD-Zhu} of the present paper,  Yetter-Drinfeld modules in the sense of H. Zhu has a trivial coaction in the cocommutative case.

Taking into account the final lines of the last paragraph, the main motivation of this paper is to give a different approach to the study of projections of Hopf braces based on the notion of left module for a Hopf brace introduced by R. Gonz\'{a}lez Rodr\'{\i}guez  in \cite[Definition 2.10]{RGON}.  Note that \cite[Definition 2.10]{RGON} is weaker than the one introduced  in \cite{Zhu} and permits to include  $(H, \triangleright =\cdot, \blacktriangleright=\circ)$ as an example of left module over the Hopf brace ${\mathbb  H}$. Moreover, in the cocommutative setting, \cite[Definition 3.1]{Zhu} and \cite[Definition 2.10]{RGON}  are equivalent.  Using the quoted notion of left module, we introduce a suitable category of Yetter-Drinfeld modules for a  Hopf brace that allows the study of Hopf braces projections and the bosonization process for these algebraic "Hopf" objects in a generic and global way at least in the cocommutative case.

The paper is organized as follows. In the first section we recall the basic notions that we will need in the rest of the article and we will review the bosonization process in a strict braided monoidal setting. Section 2 is devoted to studying   Hopf braces and their categories of modules (following   \cite{RGON}) in a braided setting.  In Section 3, we define the category of modules Yetter-Drinfeld  associated to a Hopf brace ${\mathbb H}$, denoted by $\:^{\mathbb H}_{\mathbb H}{\sf YD}$, and we prove that, if the base category is symmetric and the Hopf brace is cocommutative, this category is braided monoidal. In the last section we introduce the categories of projections, strong projections and v$_{i}$-strong projection, $i=1,2,3,4$,  over ${\mathbb H}$ and we prove that for any projection there exists two idempotents with the same image that will play the role of algebra of coinvariants. Also,  in Theorem  \ref{th-proj-mod1}, we prove  that  strong projections provide examples of left modules and we show that some constructions  of \cite{AGORE} give examples of v$_{1}$-strong projections for Hopf braces.  Moreover,  in Theorem \ref{mainequiv} we determine the conditions under which a Hopf brace ${\mathbb A}$ in $\:^{\mathbb H}_{\mathbb H}{\sf YD}$ is bosonizable in the following sense: ${\mathbb A}$ is bosonizable if when  we apply the bosonization process to  ${\mathbb A}$, the new object ${\mathbb A}\blacktriangleright\hspace{-0.1cm}\blacktriangleleft {\mathbb H}$ is a Hopf brace in the base category. Taking all this into account, in Theorem \ref{ex111} we show that $({\mathbb H}, {\mathbb A}\blacktriangleright\hspace{-0.1cm}\blacktriangleleft {\mathbb H}, x=\eta_{A}\otimes H, y=\varepsilon_{A}\otimes H)$ is a ${\rm v}_{1}$-strong projection of Hopf braces and moreover, if a projection of Hopf braces is v$_{1}$-strong, its algebra of coinvariants is an  object in $\:^{\mathbb H}_{\mathbb H}{\sf YD}$ (see Theorem \ref{YDbrace3}). On the other hand,  we prove  that, if a projection is v$_{2}$-strong, its algebra of coinvariants determines a Hopf brace in $\:^{\mathbb H}_{\mathbb H}{\sf YD}$ (see Theorem \ref{311YD}), if a projection is v$_{3}$-strong, its algebra of coinvariants determines a bosonizable Hopf brace in $\:^{\mathbb H}_{\mathbb H}{\sf YD}$ (see Theorem \ref{312YD}) and, finally, if a projection of Hopf braces  $({\mathbb H}, {\mathbb D}, x,y)$ is v$_{4}$-strong  and ${\mathbb I(q_{D})}$ is the Hopf brace associated to its algebra of coinvariants, the Hopf brace  ${\mathbb I(q_{D})}\blacktriangleright\hspace{-0.1cm}\blacktriangleleft {\mathbb H}$  is isomorphic to ${\mathbb D}$ (see Theorem \ref{313YD}). Therefore, as a consequence of these theorems, in Corollary \ref{fin} we prove that the categories of v$_{4}$-strong projections of Hopf braces with ${\mathbb H}$ fixed and the category of bosonizable Hopf braces in  $^{\mathbb H}_{\mathbb H}{\sf YD}$ are equivalent. Finally, note that $({\mathbb H}, {\mathbb A}\blacktriangleright\hspace{-0.1cm}\blacktriangleleft {\mathbb H}, x=\eta_{A}\otimes H, y=\varepsilon_{A}\otimes H)$ is an example of v$_{i}$-strong projection of Hopf braces, $i=2,3,4$, when ${\mathbb A}$ is a bosonizable Hopf brace (see Theorems \ref{ex113}, \ref{ex114} and \ref{ex115}).

\section{Preliminaries}

In this paper we will work in a monoidal setting. Following \cite{Mac}, recall that a monoidal category is a category ${\sf  C}$ together with a functor $\ot :{\sf  C}\times {\sf  C}\rightarrow {\sf  C}$, called tensor product, an object $K$ of ${\sf C}$, called the unit object, and  families of natural isomorphisms 
$$a_{M,N,P}:(M\ot N)\ot P\rightarrow M\ot (N\ot P),\;\;\;r_{M}:M\ot K\rightarrow M, \;\;\; l_{M}:K\ot M\rightarrow M,$$
in ${\sf  C}$, called  associativity, right unit and left unit constraints, respectively, satisfying the Pentagon Axiom and the Triangle Axiom, i.e.,
$$a_{M,N, P\ot Q}\co a_{M\ot N,P,Q}= (id_{M}\ot a_{N,P,Q})\co a_{M,N\ot P,Q}\co (a_{M,N,P}\ot id_{Q}),$$
$$(id_{M}\ot l_{N})\co a_{M,K,N}=r_{M}\ot id_{N},$$
where for each object $X$ in ${\sf  C}$, $id_{X}$ denotes the identity morphism of $X$. A monoidal category is called strict if the constraints of the previous paragraph are identities. It is a well-known  fact (see for example \cite{Christian}) that every non-strict monoidal category is monoidal equivalent to a strict one. Then we can assume without loss of generality that the category is strict. This lets us to treat monoidal categories as if they were strict and, as a consequence, the results proved in an strict setting hold for every non-strict symmetric monoidal category, for example the category ${\mathbb F}$-{\sf Vect} of vector spaces over a field ${\mathbb F}$, or the category $R$-{\sf Mod} of left modules over a commutative ring $R$.
For
simplicity of notation, given objects $M$, $N$, $P$ in ${\sf  C}$ and a morphism $f:M\rightarrow N$, in most cases we will write $P\ot f$ for $id_{P}\ot f$ and $f \ot P$ for $f\ot id_{P}$.

A braiding for a strict monoidal category ${\sf  C}$ is a natural family of isomorphisms 
$$c_{M,N}:M\ot N\rightarrow N\ot M$$ subject to the conditions 
$$
c_{M,N\ot P}= (N\ot c_{M,P})\co (c_{M,N}\ot P),\;\;
c_{M\ot N, P}= (c_{M,P}\ot N)\co (M\ot c_{N,P}).
$$

A strict braided monoidal category ${\sf  C}$ is a strict monoidal category with a braiding. These categories were introduced by Joyal and Street in \cite{JS1} (see also  \cite{JS2}) motivated by the theory of braids and links in topology. Note that, as a consequence of the definition, the equalities $c_{M,K}=c_{K,M}=id_{M}$ hold, for all object  $M$ of ${\sf  C}$.  If the braiding satisfies that  $c_{N,M}\co c_{M,N}=id_{M\ot N},$ for all $M$, $N$ in ${\sf  C}$, we will say that ${\sf C}$  is symmetric. In this case, we call the braiding $c$ a symmetry for the category ${\sf  C}$.

Throughout this paper ${\sf  C}$ denotes a strict braided monoidal category with tensor product $\ot$, unit object $K$ and braiding $c$. Following \cite{Bespa}, we also assume that  every idempotent morphism in ${\sf  C}$ splits, i.e., for any morphism $q:X\rightarrow X$ such that $q\co q=q$, there exist an object $I(q)$, called the image of $q$, and morphisms $i:I(q)\rightarrow X$, $p:X\rightarrow I(q)$ such that $q=i\co p$ and $p\co i=id_{I(q)}$. The morphisms $p$ and $i$ will be called a factorization of $q$. Note that $I(q)$, $p$ and $i$ are unique up to isomorphism. The categories satisfying this property constitute a broad class that includes, among others, the categories with epi-monic decomposition for morphisms and categories with equalizers or coequalizers.

\begin{definition}
{\rm 
An algebra in ${\sf  C}$ is a triple $A=(A, \eta_{A},
\mu_{A})$ where $A$ is an object in ${\sf  C}$ and
 $\eta_{A}:K\rightarrow A$ (unit), $\mu_{A}:A\otimes A
\rightarrow A$ (product) are morphisms in ${\sf  C}$ such that
$\mu_{A}\circ (A\otimes \eta_{A})=id_{A}=\mu_{A}\circ
(\eta_{A}\otimes A)$, $\mu_{A}\circ (A\otimes \mu_{A})=\mu_{A}\circ
(\mu_{A}\otimes A)$. Given two algebras $A= (A, \eta_{A}, \mu_{A})$
and $B=(B, \eta_{B}, \mu_{B})$, a morphism  $f:A\rightarrow B$ in {\sf  C} is an algebra morphism if $\mu_{B}\circ (f\otimes f)=f\circ \mu_{A}$, $ f\circ
\eta_{A}= \eta_{B}$. 

If  $A$, $B$ are algebras in ${\sf  C}$, the tensor product
$A\otimes B$ is also an algebra in
${\sf  C}$ where
$\eta_{A\otimes B}=\eta_{A}\otimes \eta_{B}$ and $\mu_{A\otimes
	B}=(\mu_{A}\otimes \mu_{B})\circ (A\otimes c_{B,A}\otimes B).$
}
\end{definition}

\begin{definition}
{\rm 
A coalgebra  in ${\sf  C}$ is a triple ${D} = (D,
\varepsilon_{D}, \delta_{D})$ where $D$ is an object in ${\sf
C}$ and $\varepsilon_{D}: D\rightarrow K$ (counit),
$\delta_{D}:D\rightarrow D\otimes D$ (coproduct) are morphisms in
${\sf  C}$ such that $(\varepsilon_{D}\otimes D)\circ
\delta_{D}= id_{D}=(D\otimes \varepsilon_{D})\circ \delta_{D}$,
$(\delta_{D}\otimes D)\circ \delta_{D}=
 (D\otimes \delta_{D})\circ \delta_{D}.$ If ${D} = (D, \varepsilon_{D},
 \delta_{D})$ and
${ E} = (E, \varepsilon_{E}, \delta_{E})$ are coalgebras, a morphism 
$f:D\rightarrow E$ in  {\sf  C} is a coalgebra morphism if $(f\otimes f)\circ
\delta_{D} =\delta_{E}\circ f$, $\varepsilon_{E}\circ f
=\varepsilon_{D}.$ 

Given  $D$, $E$ coalgebras in ${\sf  C}$, the tensor product $D\otimes E$ is a
coalgebra in ${\sf  C}$ where $\varepsilon_{D\otimes
E}=\varepsilon_{D}\otimes \varepsilon_{E}$ and $\delta_{D\otimes
E}=(D\otimes c_{D,E}\otimes E)\circ( \delta_{D}\otimes \delta_{E}).$
}
\end{definition}

\begin{definition}
	{\rm 
 Let ${D} = (D, \varepsilon_{D},
\delta_{D})$ be a coalgebra and let $A=(A, \eta_{A}, \mu_{A})$ be an
algebra. By ${\mathcal  H}(D,A)$ we denote the set of morphisms
$f:D\rightarrow A$ in ${\sf  C}$. With the convolution operation
$f\ast g= \mu_{A}\circ (f\otimes g)\circ \delta_{D}$, ${\mathcal  H}(D,A)$ is an algebra where the unit element is $\eta_{A}\circ \varepsilon_{D}=\varepsilon_{D}\otimes \eta_{A}$.
}
\end{definition}

\begin{definition}
{\rm 
 Let  $A$ be an algebra. The pair
$(M,\varphi_{M})$ is a left $A$-module if $M$ is an object in
${\sf  C}$ and $\varphi_{M}:A\otimes M\rightarrow M$ is a morphism
in ${\sf  C}$ satisfying $\varphi_{M}\circ(
\eta_{A}\ot M)=id_{M}$, $\varphi_{M}\circ (A\ot \varphi_{M})=\varphi_{M}\circ
(\mu_{A}\ot M)$. Given two left ${A}$-modules $(M,\varphi_{M})$
and $(N,\varphi_{N})$, $f:M\rightarrow N$ is a morphism of left
${A}$-modules if $\varphi_{N}\circ (A\ot f)=f\circ \varphi_{M}$.  

The  composition of morphisms of left $A$-modules is a morphism of left $A$-modules. Then left $A$-modules form a category that we will denote by $\;_{\sf A}${\sf Mod}.

Let  $D$ be a coalgebra. The pair
$(M,\rho_{M})$ is a left $D$-comodule if $M$ is an object in ${\sf  C}$ and $\rho_{M}:M\rightarrow D\ot M$ is a morphism
in ${\sf  C}$ satisfying $(\varepsilon_{D}\ot M)\co \rho_{M}=id_{M}$, $(D\ot \rho_{M})\co \rho_{M}=(\delta_{D}\ot M)\co \rho_{M}$. Given two left ${D}$-comodules $(M,\rho_{M})$
and $(N,\rho_{N})$, $f:M\rightarrow N$ is a morphism of left  ${D}$-comodules if $(D\ot f)\co \rho_{M}=\rho_{N}\co f$.  

The  composition of morphisms of left $D$-comodules is a morphism of left $C$-comodules. Then left $D$-comodules form a category that we will denote by $\;_{\sf D}${\sf Comod}.}

\end{definition}

\begin{definition}
{\rm 
We say that $X$ is a
bialgebra  in ${\sf  C}$ if $(X, \eta_{X}, \mu_{X})$ is an
algebra, $(X, \varepsilon_{X}, \delta_{X})$ is a coalgebra, and
$\varepsilon_{X}$ and $\delta_{X}$ are algebra morphisms
(equivalently, $\eta_{X}$ and $\mu_{X}$ are coalgebra morphisms). Moreover, if there exists a morphism $\lambda_{X}:X\rightarrow X$ in ${\sf  C}$,
called the antipode of $X$, satisfying that $\lambda_{X}$ is the inverse of $id_{X}$ in ${\mathcal  H}(X,X)$, i.e., 
\begin{equation}
\label{antipode}
id_{X}\ast \lambda_{X}= \eta_{X}\circ \varepsilon_{X}= \lambda_{X}\ast id_{X},
\end{equation}
we say that $X$ is a Hopf algebra. A morphism of Hopf algebras is an algebra-coalgebra morphism. Note that, if $f:X\rightarrow Y$ is a Hopf algebra morphism the following equality holds:
$$
%\label{morant}
\lambda_{Y}\co f=f\co \lambda_{X}.
$$

With the composition of morphisms in {\sf C} we can define a category whose objects are  Hopf algebras  and whose morphisms are morphisms of Hopf algebras. We denote this category by ${\sf  Hopf}$.

A Hopf algebra is commutative if $\mu_{X}\co c_{X,X}=\mu_{X}$ and cocommutative if $c_{X,X}\co \delta_{X}=\delta_{X}.$ It is easy to see that in both cases $\lambda_{X}\circ \lambda_{X} =id_{X}$ and then $\lambda_{X}$ is an isomorphism with inverse $\lambda_{X}^{-1}=\lambda_{X}$.
}
\end{definition}

If $X$ is a Hopf algebra,  the antipode is antimultiplicative and anticomultiplicative 
\begin{equation}
\label{a-antip}
\lambda_{X}\co \mu_{X}=  \mu_{X}\co (\lambda_{X}\ot \lambda_{X})\co c_{X,X},\;\;\;\; \delta_{X}\co \lambda_{X}=c_{X,X}\co (\lambda_{X}\ot \lambda_{X})\co \delta_{X}, 
\end{equation}
and leaves the unit and counit invariant, i.e., 
\begin{equation}
\label{u-antip}
\lambda_{X}\co \eta_{X}=  \eta_{X},\;\; \varepsilon_{X}\co \lambda_{X}=\varepsilon_{X}.
\end{equation}

Also $X$ becomes a left $X$-module by the adjoint action which is defined by 
$$\varphi_{X}^{ad}=\mu_{X}\co (\mu_{X}\ot \lambda_{X})\co (X\ot c_{X,X})\co (\delta_{X}\ot X), $$
and a left $X$-comodule by the adjoint coaction $$\rho_{X}^{ad}=(\mu_{X}\ot X)\co (X\ot c_{X,X})\co (\delta_{X}\ot \lambda_{X})\co \delta_{X}.$$

In the following definition we recall the notions of left (co)module (co)algebra. The notions of  right (co)module (co)algebra are similar.

\begin{definition}
{\rm 
Let $X$ be a Hopf algebra. An algebra $A$  is said to be a left $X$-module algebra if $(A, \varphi_{A})$ is a left $X$-module and $\eta_{A}$, $\mu_{A}$ are morphisms of left $X$-modules, i.e.,
\begin{equation}
\label{mod-alg}
\varphi_{A}\circ (X\otimes \eta_{A})=\varepsilon_{X}\otimes \eta_{A},\;\;\varphi_{A}\circ (X\otimes \mu_{A})=\mu_{A}\circ \varphi_{A\otimes A},
\end{equation}
where  $\varphi_{A\otimes A}=(\varphi_{A}\otimes \varphi_{A})\circ (X\otimes c_{X,A}\otimes A)\circ (\delta_{X}\otimes A\otimes A)$ is the left action on $A\otimes A$. For example, $X$ with the adjoint action  $\varphi_{X}^{ad}$ is a left  $X$-module algebra.

On the other hand,  $A$ is said to be a left $X$-comodule algebra if $(A,\rho_{A})$ is a left $X$-comodule and $\eta_{A}$ and $\mu_{A}$ are morphisms of left $X$-comodules, i.e.,
\begin{equation}
\label{comod-alg} \rho_{A}\circ \eta_{A}=\eta_{X}\otimes \eta_{A},\;\; \rho_{A}\circ \mu_{A}=(X\otimes \mu_{A})\circ \rho_{A\otimes A}
\end{equation}
where $\rho_{A\otimes A}=(\mu_{X}\ot A\ot A)\co (X\ot c_{A,X}\ot A)\co (\rho_{A}\ot \rho_{A})$ is the coaction on $A\ot A$.  

In a similar way we can define the notion of left $X$-module coalgebra  and left $X$-comodule coalgebra. Then, a coalgebra $B$  is said to be a left $X$-module coalgebra if $(B,\varphi_{B})$ is a left $X$-module and  $\varepsilon_{B}$ and $\delta_{B}$ are morphisms of left $X$-modules, i.e.,
$$
%\label{mod-coalg}
\varepsilon_{B}\circ \varphi_{B}=\varepsilon_{X}\otimes \varepsilon_{B},\;\;\delta_{B}\co \varphi_{B}=\varphi_{B\ot B}\co (X\ot \delta_{B}).
$$

Finally, a coalgebra $B$  is said to be a left $X$-comodule coalgebra if $(B,\rho_{B})$ is a left $X$-comodule and  $\varepsilon_{B}$ and $\delta_{B}$ are morphisms of left $X$-comodules, i.e.,
\begin{equation}
\label{comod-coalg}
(X\otimes \varepsilon_{B})\circ \rho_{B}=\eta_{X}\otimes \varepsilon_{B},\;\;(X\otimes \delta_{B})\circ \rho_{B}=\rho_{B\otimes B}\circ \delta_{B}.
\end{equation}
For example, $X$ with the adjoint coaction  $\rho_{X}^{ad}$ is a left  $X$-comodule coalgebra.
}
\end{definition}

If $(A,\varphi_{A})$ is a left $X$-module algebra, 
$$A\sharp X=(A\otimes X, \eta_{A\sharp X}=\eta_{A}\otimes \eta_{X}, \mu_{A\sharp X}=(\mu_{A}\otimes \mu_{X})\circ (A\otimes \Psi_{A}^{X}\otimes X))$$ where 
$\Psi_{A}^{X}=(\varphi_{A}\otimes X)\co  (X\otimes c_{X,A})\circ (\delta_{X}\otimes A),$ is an algebra called the smash product of $A$ and $X$. Similarly, if $(B, \rho_{B})$ is a left $X$-comodule coalgebra, we can define the coalgebra smash coproduct of $B$ and $X$ as $$B\propto X=(B\otimes X,\varepsilon_{B\propto X}=\varepsilon_{B}\otimes \varepsilon_{X}, \delta_{B\propto X}=(B\otimes \Omega_{X}^{B}\ot X)\circ (\delta_{B}\otimes \delta_{X})),$$ where $\Omega_{X}^{B}= (\mu_{X}\otimes B)\circ ( X\otimes c_{B,X})\circ
(\rho_{B}\otimes X).$

\begin{definition}
{\rm 
Let $X$ be a Hopf algebra in ${\sf  C}$. We shall denote by $^{\sf X}_{\sf X}{\sf  Y}{\sf  D}$ the category of left Yetter-Dinfeld modules over $X$. More concretely, a triple $M=(M, \varphi_{M}, \rho_{M})$ is an object in $^{\sf X}_{\sf X}{\sf  Y}{\sf  D}$ if $(M,\varphi_{M})$ is a left $X$-module, $(M,\rho_{M})$ is a left $X$-comodule and the following identity
$$
%\label{YD1}
(\mu_{X}\otimes M)\circ (X\otimes c_{M,X})\circ ((\rho_{M}\circ \varphi_{M})\otimes X)\circ (X\otimes c_{X,M})\circ (\delta_{X}\otimes M)
$$
$$=(\mu_{X}\otimes \varphi_{M})\circ (X\otimes c_{X,X}\otimes M)\circ (\delta_{X}\otimes \rho_{M}).$$
holds. The morphisms in $\;^{\sf X}_{\sf X}{\sf  Y}{\sf  D}$ are morphisms of left modules and comodules.

For example, for any Hopf algebra $X$, $(X, \varphi_{X}^{ad}, \rho_{X}=\delta_{X}),\;\; (X, \varphi_{X}=\mu_{X}, \rho_{X}^{ad})$ are left Yetter-Drinfeld modules over $X$. Also, any left $X$-module $(M,\varphi_{M})$ over a cocommutative Hopf algebra $X$ is a Yetter-Drinfeld module with the trivial left coaction $\rho_{M}=\eta_{X}\ot M$. Finally, the triple $(M, \varphi_{M}=\varepsilon_{X}\ot M, \rho_{M}=\eta_{X}\ot M)$ is a left Yetter-Drinfeld module for all Hopf algebra $X$. 
}
\end{definition}

The category $^{\sf X}_{\sf X}{\sf  Y}{\sf  D}$ is  strict monoidal with the usual tensor product in ${\sf  C}$. For
$M$, $N$ in $\;^{X}_{X}{\sf  Y}{\sf  D}$, $M\otimes N$ has the tensor module and comodule structures given by $$\varphi_{M\otimes N}=(\varphi_{M}\otimes \varphi_{N})\circ (X\otimes c_{X,M}\otimes N)\circ (\delta_{X}\otimes M\otimes N)$$
and $$\rho_{M\otimes N}=(\mu_{X}\otimes M\otimes N)\circ (X\otimes c_{M,X}\otimes N)\circ (\rho_{M}\otimes \rho_{N}).$$

If the antipode of $X$ is an isomorphism,  $\;^{\sf X}_{\sf X}{\sf  Y}{\sf  D}$ is a braided monoidal category where the braiding $t_{M,N}:M\otimes N\rightarrow N\otimes N,$ is given by $t_{M,N}=(\varphi_{N}\otimes M)\circ (X\otimes c_{M,N})\circ (\rho_{M}\otimes N).$ It is immediate to see that $t_{M,N}$ is natural and it is an isomorphism with inverse
$$t_{M,N}^{-1}=c_{M,N}^{-1}\circ (\varphi_{N}\otimes M)\circ (\lambda_{X}^{-1}\otimes N\otimes M)\circ (c_{X,N}^{-1}\otimes M)\circ (N\otimes \rho_{M}).$$

Then, if X is a Hopf algebra with $\lambda_{X}$ isomorphism,  a Hopf algebra  in  $^{\sf X}_{\sf X}{\sf  Y}{\sf  D}$ is an object $(A, \varphi_{A}, \rho_{A})$ in $^{\sf X}_{\sf X}{\sf  Y}{\sf  D}$ such that it is an algebra-coalgebra in ${\sf  C}$ with an endomorphism $\lambda_{A}:A\rightarrow A$ satisfying the following: $(A,\varphi_{A})$ is a left $X$-module (co)algebra,  $(A,\rho_{A})$ is a left $X$-comodule (co)algebra, $\lambda_{A}$ is a morphism of left $X$-modules and left $X$-comodules,  for  $\varepsilon_{A}$,  $\delta_{A}$ the following identities 
$$\varepsilon_{A}\co\eta_{A}=id_{K}, \; \varepsilon_{A}\co \mu_{A}=\varepsilon_{A}\ot \varepsilon_{A}, \; 
\eta_{A}\ot \eta_{A}=\delta_{A}\co \eta_{A},$$
$$\delta_{A}\co \mu_{A}=(\mu_{A}\otimes \mu_{A})\circ (A\otimes t_{A,A}\otimes A)\circ (\delta_{A}\ot \delta_{A}), $$
hold  and, finally, $\lambda_{A}$ is the inverse of $id_{A}$ in ${\mathcal  H}(A,A)$. Then, the Hopf algebra $X$ with $\varphi_{X}=\varepsilon_{X}\ot X$ and $\rho_{X}=\eta_{X}\ot X$ is a Hopf algebra in $^{\sf X}_{\sf X}{\sf  Y}{\sf  D}$. Note that in this case $ t_{X,X}=c_{X,X}$.

In the following paragraphs of this section we briefly summarize some results from \cite{MN}, \cite{MAJ2} and \cite{RAD} about projections of Hopf algebras and the bosonization process in a monoidal setting. 

\begin{definition}
{\rm A projection of Hopf algebras in {\sf C} is a 4-tupla $(X,Y, f,h)$ where $X$, $Y$ are Hopf algebras, and $f:X\rightarrow Y$, $h:Y\rightarrow X$ are Hopf algebra morphisms such that $h\co f=id_{X}$. 
		
A morphism between  projections of Hopf algebras $(X,Y, f,h)$ and $(X^{\prime},Y^{\prime}, f^{\prime},h^{\prime})$ is a pair $(r,s)$, where $r:X\rightarrow X^{\prime}$, $s:Y\rightarrow Y^{\prime}$ are Hopf algebra morphisms such that
$$
%\label{proj-morph}
s\circ f=f^{\prime}\circ r, \;\; r\circ h=h^{\prime}\circ s. 
$$

With the obvious composition of morphisms we can define a category whose objects are  Hopf algebra projections and whose morphisms are morphisms of Hopf algebra projections. We denote this category by ${\sf  P}({\sf  Hopf})$.

It is obvious that there exists a functor ${\sf P}:{\sf Hopf}\rightarrow {\sf P(Hopf)}$ defined on objects by $${\sf P}(X)=(X,X, id_{X}, id_{X})$$ and on morphisms by ${\sf P}(f)=(f,f)$.
}
\end{definition}

Let $(X, Y, f, h)$ be an object in ${\sf  P}({\sf  Hopf})$.   The morphism 
$q_{Y}=id_{Y}\ast (f\co \lambda_{X}\co h)$ is  idempotent and, as a consequence, there exist an epimorphism $p_{Y}$, a monomorphism $i_{Y}$, and an object $I(q_{Y})$ such that $q_{Y}=i_{Y}\circ p_{Y}$ and $p_{Y}\circ i_{Y}=id_{I(q_{Y})}$. As a consequence,
$$
\setlength{\unitlength}{1mm}
\begin{picture}(101.00,12.00)
\put(20,7){\vector(1,0){20}}
\put(46,8){\vector(1,0){27}}
\put(46,6){\vector(1,0){27}}
\put(13,7){\makebox(0,0){$I(q_{Y})$}}
\put(43,7){\makebox(0,0){$Y$}}
\put(80,7){\makebox(0,0){$Y\otimes X$}}
\put(27.5,10){\makebox(0,0){$i_{Y}$}}
\put(60,11){\makebox(0,0){$(Y\ot h)\co \delta_{Y}$}}
\put(60,3){\makebox(0,0){$Y\otimes \eta_{X}$}}
\end{picture}
$$
is an equalizer diagram and $I(q_{Y})$ is a left $X$-module algebra where  the algebra structure is defined by
\begin{equation}
\label{alg-id} 
\eta_{I(q_{Y})}=p_{Y}\co \eta_{Y},\;\;\; \mu_{I(q_{Y})}=p_{Y}\co \mu_{Y}\co (i_{Y}\ot i_{Y}),
\end{equation}
i.e., $\eta_{I(q_{Y})}$ is the unique morphism such that $i_{Y}\co \eta_{I(q_{Y})}=\eta_{Y}$ and $\mu_{I(q_{Y})}$ is the unique morphism such that 
\begin{equation}
\label{m-i}
i_{Y}\co \mu_{I(q_{Y})}=\mu_{Y}\co (i_{Y}\ot i_{Y}).
\end{equation} 

The action $\varphi_{I(q_{Y})}: X\ot I(q_{Y})\rightarrow I(q_{Y})$ is 
$
\varphi_{I(q_{Y})}=p_{Y}\co\mu_{Y}\co (f\ot i_{Y}),
$
and then $\varphi_{I(q_{Y})}$ is the unique morphism such that 
$$
%\label{phi2}
i_{Y}\co \varphi_{I(q_{Y})}=\varphi_{Y}^{ad}\co (f\ot i_{Y}).
$$

On the other hand,
$$
\setlength{\unitlength}{1mm}
\begin{picture}(101.00,10.00)
\put(20.00,7.00){\vector(1,0){25.00}}
\put(20.00,5.00){\vector(1,0){25.00}}
\put(55.00,6.00){\vector(1,0){21.00}}
\put(32.00,11.00){\makebox(0,0)[cc]{$\mu_{Y}\co (Y\ot f)$ }}
\put(32.00,1.00){\makebox(0,0)[cc]{$Y\otimes \varepsilon_{X} $ }}
\put(65.00,9.00){\makebox(0,0)[cc]{$p_{Y} $ }}
\put(13.00,6.00){\makebox(0,0)[cc]{$ Y\otimes X$ }}
\put(50.00,6.00){\makebox(0,0)[cc]{$ Y$ }}
\put(83.00,6.00){\makebox(0,0)[cc]{$I(q_{Y}) $ }}
\end{picture}
$$
is a coequalizer diagram and, as a consequence, $I(q_{Y})$ is a left $X$-comodule coalgebra with
\begin{equation}
\label{coalg-id} 
\varepsilon_{I(q_{Y})}=\varepsilon_{Y}\co i_{Y},\;\; \delta_{I(q_{Y})}=(p_{Y}\otimes p_{Y})\circ \delta_{Y}\co i_{Y}
\end{equation}
and coaction $\rho_{I(q_{Y})}: I(q_{Y})\rightarrow X\ot I(q_{Y})$ defined by 
$
\rho_{I(q_{Y})}=(h\otimes p_{Y})\circ \delta_{Y}\circ i_{Y}.
$

In this case $\varepsilon_{I(q_{Y})}$ is the unique morphism such that $\varepsilon_{I(q_{Y})}\co p_{Y}=\varepsilon_{Y}$, $\delta_{I(q_{Y})}$ is the unique morphism such that 
$$
%\label{co-i}
\delta_{I(q_{Y})}\co p_{Y}=(p_{Y}\ot p_{Y})\co \delta_{Y},
$$
and  the coaction $\rho_{I(q_{Y})}$ is the unique morphism satisfying 
$$
%\label{rho2}
\rho_{I(q_{Y})}\co p_{Y}=(h\ot p_{Y})\co \rho^{ad}_{Y}.
$$

The algebra-coalgebra $I(q_{Y})$, with the action $\varphi_{I(q_{Y})}$ and the coaction $\rho_{I(q_{Y})}$, is a Hopf algebra in   $^{\sf X}_{\sf X}{\sf  Y}{\sf  D}$ with antipode $\lambda_{I(q_{Y})}=\varphi_{I(q_{Y})}\circ (X\otimes (p_{Y}\circ \lambda_{Y}\circ i_{Y}))\circ \rho_{I(q_{Y})}.$ 

 Also, using that $i_{Y}$ is an equalizer morphism and $p_{Y}$ is a coequalizer, we obtain the following identities:
\begin{equation}
\label{id1}
p_{Y}\co \mu_{Y}\co (Y\ot q_{Y})=p_{Y}\co \mu_{Y}, \;\;\;\;\;
(Y\ot q_{Y})\co \delta_{Y}\co i_{Y}= \delta_{Y}\co i_{Y}. 
\end{equation}

For the  Hopf algebra $I(q_{Y})$  in  $^{\sf X}_{\sf X}{\sf  Y}{\sf  D}$ we can apply the monoidal version of the construction introduced by Radford in \cite{RAD}, and extended to the quantum  setting  by Majid \cite{MAJ2}, producing a Hopf algebra $I(q_{Y})\blacktriangleright\hspace{-0.1cm}\blacktriangleleft X$ in ${\sf  C}$, called by Majid the bosonization of $I(q_{Y})$, with the following structure: The Hopf algebra $I(q_{Y})\blacktriangleright\hspace{-0.1cm}\blacktriangleleft X$ is the smash product $I(q_{Y})\sharp X$ as algebra, the smash coproduct $I(q_{Y})\propto X$ as coalgebra, and the antipode is defined by $$\lambda_{I(q_{Y})\blacktriangleright\hspace{-0.1cm}\blacktriangleleft X}=\Psi^{X}_{I(q_{Y})}\circ (\lambda_{X}\otimes
\lambda_{I(q_{Y})})\circ \Omega_{X}^{I(q_{Y})}.$$

Moreover, 
\begin{equation}
\label{nnu}
\nu_{Y}=(p_{Y}\ot h)\co \delta_{Y}:Y\rightarrow I(q_{Y})\blacktriangleright\hspace{-0.1cm}\blacktriangleleft X
\end{equation} 
is a Hopf algebra isomorphism with inverse with inverse $\nu_{Y}^{-1}=\mu_{Y}\co (i_{Y}\ot f):I(q_{Y})\blacktriangleright\hspace{-0.1cm}\blacktriangleleft X\rightarrow Y.$

The existence of the previous isomorphism  is the main tool  to obtain a categorical equivalence between the category of Hopf algebras in $^{\sf X}_{\sf X}{\sf  Y}{\sf  D}$ and the category of Hopf algebra projections associated to a fixed $X$ with invertible antipode. This categorical equivalence is a corollary of the more general result proved in \cite{NikaRamon4} for weak Hopf algebras.

Finally, note that $i_{Y}$ is a coalgebra morphism iff 
\begin{equation}
\label{iz-coal}
(q_{Y}\ot Y)\co \delta_{Y}\co i_{Y}=\delta_{Y}\co i_{Y}.
\end{equation}

Equivalently, $i_{Y}$ is a coalgebra morphism iff  $\rho_{I(q_{Y})}=\eta_{X}\ot I(q_{Y})$ (see \cite{MN}). Therefore in this case, $\varepsilon_{I(q_{Y})\blacktriangleright\hspace{-0.1cm}\blacktriangleleft X}=\varepsilon_{I(q_{Y})}\ot \varepsilon_{X},\;\; \delta_{I(q_{Y})\blacktriangleright\hspace{-0.1cm}\blacktriangleleft X}=\delta_{I(q_{Y})\ot X}$ and $\lambda_{I(q_{Y})\blacktriangleright\hspace{-0.1cm}\blacktriangleleft X}=\Psi^{X}_{I(q_{Y})}\circ (\lambda_{X}\otimes
\lambda_{I(q_{Y})})\circ c_{I(q_{Y}),X}.$

Note that, if $Y$ is cocommutative, condition (\ref{iz-coal}) always holds. This fact was proved by Sweedler in \cite{SW} for projections of Hopf algebras in a category of vector spaces. On the other hand, there exist examples where  $i_{Y}$ it is not a coalgebra morphism (see \cite{BCM} for the complete details). In any case, if $i_Y$ is a coalgebra morphism, we have that $I(q_{Y})$ is a Hopf algebra in ${\sf C}$ because $\rho_{I(q_{Y})}$ is trivial.

Similarly, $p_{Y}$ is an algebra morphism iff 
\begin{equation}
%\label{pz-al}
p_{Y}\co \mu_{Y}\co (q_{Y}\ot Y)=p_{Y}\co \mu_{Y}.
\end{equation}

Equivalently, $p_{Y}$ is an algebra morphism iff  $\varphi_{I(q_{Y})}=\varepsilon_{X}\ot I(q_{Y})$ (see \cite{MN}). Therefore in this case, $\eta_{I(q_{Y})\blacktriangleright\hspace{-0.1cm}\blacktriangleleft X}=\eta_{I(q_{Y})}\ot \eta_{X},\;\; \mu_{I(q_{Y})\blacktriangleright\hspace{-0.1cm}\blacktriangleleft X}=\mu_{I(q_{Y})\ot X}$ and 
$\lambda_{I(q_{Y})\blacktriangleright\hspace{-0.1cm}\blacktriangleleft X}=c_{X,I(q_{Y})}\circ(\lambda_{X}\otimes\lambda_{I(q_{Y})})\circ\Omega_{X}^{I(q_{Y})}.$ Also, if $p_Y$ is an algebra morphism, we have that $I(q_{Y})$ is a Hopf algebra in ${\sf C}$ because $\varphi_{I(q_{Y})}$ is trivial.

Finally, we have the following result.

\begin{lemma}
\label{cole}
Let $(X, Y, f, h)$ be an object in ${\sf  P}({\sf  Hopf})$. If $Y$ is cocommutative, the morphism $q_{Y}$ is a coalgebra morphism. Also, under these conditions, the following equality 
\begin{equation}
\label{lxy}
i_{Y}\circ \lambda_{I(q_{Y})}=\lambda_{Y}\circ i_{Y}.
\end{equation}
holds.
\end{lemma}

\begin{proof}
Trivially $q_{Y}$ preserves the counit. On the other hand, 
\begin{itemize}
	\item[ ]$\hspace{0.38cm}\delta_{Y}\circ q_{Y} $ 
	\item [ ]$= \mu_{Y\otimes Y}\circ (\delta_{Y}\otimes (\delta_{Y}\circ f\circ \lambda_{X}\circ h))\circ \delta_{Y}$ {\scriptsize (by the condition of algebra morphism for $\delta_{Y}$)}
	\item [ ]$=\mu_{Y\otimes Y}\circ (\delta_{Y}\otimes ((( f\circ \lambda_{X}\circ h)\otimes (f\circ \lambda_{X}\circ h))\circ\delta_{Y}))\circ  \delta_{Y} $ {\scriptsize (by (\ref{a-antip}), the condition of Hopf algebra morphisms }
	\item[ ]$\hspace{0.38cm}$ {\scriptsize for $f$ and $h$ and the cocommutativity of $\delta_{Y}$)}
	\item [ ]$= ((\mu_{Y}\circ (Y\otimes (f\circ \lambda_{X}\circ h)))\otimes (\mu_{Y}\circ (Y\otimes (f\circ \lambda_{X}\circ h))))\circ (Y\otimes (c_{Y,Y}\circ \delta_{Y})\otimes Y) \circ(Y\otimes\delta_{Y})\circ\delta_{Y} $ 
	\item[ ]$\hspace{0.38cm}$ {\scriptsize (by the coassociativity of $\delta_{Y}$ and the naturality of $c$)}
	\item [ ]$=( q_{Y}\otimes  q_{Y})\circ \delta_{Y}$ {\scriptsize (by the coassociativity  and cocommutativity of $\delta_{Y}$)}
\end{itemize}
holds, and as a consequence $q_{Y}$ is a coalgebra morphism.

If $Y$ is cocommutative we have that $\rho_{I(q_{Y})}=\eta_{X}\ot I(q_{Y})$ and this implies that $\lambda_{I(q_{Y})}= p_{Y}\circ \lambda_{Y}\circ i_{Y}.$ Then,
\begin{itemize}
	\item[ ]$\hspace{0.38cm} i_{Y}\circ \lambda_{I(q_{Y})}$ 
	\item [ ]$=q_Y\circ \lambda_Y\circ i_Y $ {\scriptsize (by $\lambda_{I(q_{Y})}= p_{Y}\circ \lambda_{Y}\circ i_{Y}$)}
	\item [ ]$=\mu_Y\circ (\lambda_Y\otimes (f\circ h\circ \lambda_{Y}\circ \lambda_{Y}))\circ \delta_{Y}\circ i_Y $ {\scriptsize (by (\ref{a-antip}), the condition of Hopf algebra morphisms  for $f$ and $h$}
	\item[ ]$\hspace{0.38cm}$ {\scriptsize and the cocommutativity of $\delta_{Y}$)}
	\item [ ]$=\mu_Y\circ (\lambda_Y\otimes (f\circ h))\circ \delta_{Y}\circ i_Y  $ {\scriptsize (by $\lambda_{Y}\circ \lambda_{Y}=id_{Y}$)}
	\item [ ]$=\mu_Y\circ (\lambda_Y\otimes (f\circ \eta_{X}))\circ i_Y  $ {\scriptsize (by the equalizer condition for $i_{Y}$)}
	\item [ ]$=\lambda_{Y}\circ i_{Y}$ {\scriptsize (by the unit properties).}
\end{itemize}
\end{proof}

\section{Modules  for Hopf braces}

The main objective of this section is to present the main properties of the modules, in the sense of \cite{RGON},  associated with a Hopf brace. We begin the section with the definition of Hopf brace  in a braided monoidal category ${\sf C}$.

\begin{definition}
\label{H-brace}
{\rm Let $H=(H, \varepsilon_{H}, \delta_{H})$ be a coalgebra in {\sf C}. Let's assume that there are two algebra structures $(H, \eta_{H}^1, \mu_{H}^1)$, $(H, \eta_{H}^2, \mu_{H}^2)$ defined on $H$ and suppose that there exist two endomorphism of $H$ denoted by $\lambda_{H}^{1}$ and $\lambda_{H}^{2}$. We will say that 
$$(H, \eta_{H}^{1}, \mu_{H}^{1}, \eta_{H}^{2}, \mu_{H}^{2}, \varepsilon_{H}, \delta_{H}, \lambda_{H}^{1}, \lambda_{H}^{2})$$
is a Hopf brace in {\sf C} if:
\begin{itemize}
\item[(i)] $H_{1}=(H, \eta_{H}^{1}, \mu_{H}^{1},  \varepsilon_{H}, \delta_{H}, \lambda_{H}^{1})$ is a Hopf algebra in {\sf C}.
\item[(ii)] $H_{2}=(H, \eta_{H}^{2}, \mu_{H}^{2},  \varepsilon_{H}, \delta_{H}, \lambda_{H}^{2})$ is Hopf algebra in {\sf C}.
\item[(iii)] The  following equality holds:
$$\mu_{H}^{2}\co (H\ot \mu_{H}^{1})=\mu_{H}^{1}\co (\mu_{H}^{2}\ot \Gamma_{H_{1}} )\co (H\ot c_{H,H}\ot H)\co (\delta_{H}\ot H\ot H),$$
\end{itemize}
where  $$\Gamma_{H_{1}}=\mu_{H}^{1}\co (\lambda_{H}^{1}\ot \mu_{H}^{2})\co (\delta_{H}\ot H).$$

Following \cite{RGON}, a Hopf brace will be denoted by ${\mathbb H}=(H_{1}, H_{2})$ or  in a simpler way by ${\mathbb H}$.

}
\end{definition}

The previous definition is the general notion of Hopf brace in a braided monoidal setting. If we restrict it to a category of Yetter-Drinfeld modules over a Hopf algebra which antipode is an isomorphism we obtain the definition of braided Hopf brace introduced by H. Zhu and Z. Ying in \cite[Definition 2.1]{Zhu2}.

\begin{definition}
%\label{H-coco}
 If  ${\mathbb H}$ is a Hopf brace in {\sf C}, we will say that ${\mathbb H}$ is cocommutative if $\delta_{H}=c_{H,H}\circ \delta_{H}$, i.e., $H_{1}$ and $H_{2}$ are cocommutative Hopf algebras in {\sf C}.
 
Note that by \cite[Corollary 5]{Sch}, if $H$ is a  cocommutative Hopf algebra  in the  braided monoidal category {\sf C}, the identity 
\begin{equation}
\label{ccb}
c_{H,H}\circ c_{H,H}=id_{H\otimes H} 
\end{equation}
holds.
\end{definition}

\begin{definition}
%\label{mor}
{\rm  Given two Hopf braces ${\mathbb H}$  and  ${\mathbb B}$ in {\sf C}, a morphism $x$ in {\sf C} between the two underlying objects is called a morphism of Hopf braces if both $x:H_{1}\rightarrow B_{1}$ and $x:H_{2}\rightarrow B_{2}$ are algebra-coalgebra morphisms.
		
Hopf braces together with morphisms of Hopf braces form a category which we denote by {\sf HBr}. 

\begin{theorem}
\label{2-th1}
	There exists a functor between the categories {\sf Hopf} and {\sf HBr}.
\end{theorem}

\begin{proof}
	If $H$ is a Hopf algebra, ${\mathbb H}_{triv}=(H,H, \eta_{H}, \mu_{H}, \eta_{H},\mu_{H}, \varepsilon_{H}, \delta_{H}, \lambda_{H}, \lambda_{H})$ is an object in {\sf HBr}. On the other hand, if $x:H\rightarrow B$ is a morphism of Hopf algebras, we have that the pair $(x,x)$ is a morphism in {\sf HBr} between ${\mathbb H}_{triv}$  and ${\mathbb B}_{triv}$. Therefore, there exists a functor $${\sf H}^{\prime}:{\sf Hopf}\rightarrow {\sf HBr}$$ defined on objects by ${\sf H}^{\prime}(H)={\mathbb H}_{triv}$ and on morphisms by ${\sf H}^{\prime}(x)=(x,x)$.
\end{proof}
		
}
\end{definition}

Let ${\mathbb H}$  be a Hopf brace in {\sf C}. Then 
$$
%\label{eb1}
\eta_{H}^{1}=\eta_{H}^2, 
$$
holds and, by \cite[Lemma 1.7]{AGV},  in this braided setting  the equality
\begin{equation}
\label{agv1}
\Gamma_{H_{1}}\circ (H\otimes \lambda_{H}^1)=\mu_{H}^{1}\circ ((\lambda_{H}^1\circ \mu_{H}^{2})\otimes H)\circ (H\otimes c_{H,H}) \circ (\delta_{H}\otimes H)
\end{equation}
also holds. Moreover, in our braided context \cite[Lemma 1.8]{AGV} and \cite[Remark 1.9]{AGV} hold and then we have that $(H,\eta_{H}^{1}, \mu_{H}^{1})$ is a left $H_{2}$-module algebra with action $\Gamma_{H_{1}}$ and $\mu_{H}^2$ admits the following expression:
\begin{equation}
\label{eb2}
\mu_{H}^2=\mu_{H}^{1}\circ (H\otimes \Gamma_{H_{1}})\circ (\delta_{H}\otimes H). 
\end{equation}

Finally, taking into account that every Hopf brace is an example of Hopf truss, by \cite[Theorem 6.4]{BRZ1}, we have that $(H,\eta_{H}^{1}, \mu_{H}^{1})$ is also a left $H_{2}$-module algebra with action 
$$\Gamma_{H_{1}}^{\prime}=\mu_{H}^{1}\circ (\mu_{H}^{2}\otimes \lambda_{H}^{1})\circ (H\otimes c_{H,H})\circ (\delta_{H}\otimes H)$$
because the symmetry  is not needed in the proof as in the case of $\Gamma_{H_{1}}$.

Finally, by the naturality of $c$ and the coassociativity of $\delta_{H}$, we obtain that 
$$
%\label{alt}
\mu_{H}^{1}\co (\mu_{H}^{2}\ot \Gamma_{H_{1}} )\co (H\ot c_{H,H}\ot H)\co (\delta_{H}\ot H\ot H)
$$
$$=\mu_{H}^{1}\co (\Gamma_{H_{1}}^{\prime}\otimes  \mu_{H}^{2})\co (H\ot c_{H,H}\ot H)\co (\delta_{H}\ot H\ot H)$$
and then (iii) of Definition \ref{H-brace} is equivalent to 
\begin{equation}
\label{iii-eq}
\mu_{H}^{2}\co (H\otimes \mu_{H}^{1})
=\mu_{H}^{1}\co (\Gamma_{H_{1}}^{\prime}\otimes  \mu_{H}^{2})\co (H\ot c_{H,H}\ot H)\co (\delta_{H}\ot H\ot H).
\end{equation}
Therefore, the equality 
$$
%\label{iii-eq1}
\mu_{H}^{2}
=\mu_{H}^{1}\co (\Gamma_{H_{1}}^{\prime}\otimes  H)\co (H\ot c_{H,H})\co (\delta_{H}\ot H)
$$
holds.

\begin{lemma}
%\label{dphi}
Let ${\mathbb H}$  be a Hopf brace in {\sf C}. If ${\mathbb H}$ is cocommutative, $\Gamma_{H_{1}}$ is a coalgebra morphism.
\end{lemma}

\begin{proof}
Trivially $\varepsilon_{H}\circ \Gamma_{H_{1}}=\varepsilon_{H}\otimes \varepsilon_{H}$. Moreover, 
\begin{itemize}
\item[ ]$\hspace{0.38cm}\delta_{H}\circ \Gamma_{H_{1}} $ 
\item [ ]$= \mu_{H_{1}\otimes H_{1}}\circ (((\lambda_{H}^{1}\otimes \lambda_{H}^{1})\circ c_{H,H}\circ \delta_{H})\otimes (\mu_{H_{2}\otimes H_{2}}\circ (\delta_{H}\otimes \delta_{H})))\circ (\delta_{H}\otimes H)$ {\scriptsize (by the condition of}
\item[ ]$\hspace{0.38cm}${\scriptsize  coalgebra morphisms for $\mu_{H}^{1}$ and $\mu_{H}^{2}$ and (\ref{a-antip}))}
\item [ ]$=(\Gamma_{H_{1}}\otimes \Gamma_{H_{1}})\circ \delta_{H\otimes H} ${\scriptsize (by the naturality of $c$ and the cocommutativity and coassociativity conditions)}
\end{itemize}
\end{proof}

\begin{lemma}
\label{dphi-p}
Let ${\mathbb H}$  be a Hopf brace in {\sf C}. If ${\mathbb H}$ is cocommutative,  $\Gamma_{H_{1}}^{\prime}$ is a coalgebra morphism.
\end{lemma}

\begin{proof}
As in the case of $\Gamma_{H_{1}}$, trivially $\varepsilon_{H}\circ \Gamma_{H_{1}}^{\prime}=\varepsilon_{H}\otimes \varepsilon_{H}$. Moreover,
\begin{itemize}
\item[ ]$\hspace{0.38cm}\delta_{H}\circ \Gamma_{H_{1}}^{\prime} $ 
\item [ ]$= \mu_{H_{1}\otimes H_{1}}\circ ( (\mu_{H_{2}\otimes H_{2}}\circ (\delta_{H}\otimes \delta_{H}))\otimes ((\lambda_{H}^{1}\otimes \lambda_{H}^{1})\circ  \delta_{H}))\circ (H\otimes c_{H,H})\circ (\delta_{H}\otimes H)$ 
\item[ ]$\hspace{0.38cm}${\scriptsize (by the condition of coalgebra morphisms for $\mu_{H}^{1}$ and $\mu_{H}^{2}$, (\ref{a-antip}) and cocommutativity of $\delta_{H}$)}
\item [ ]$= ((\mu_{H}^{1}\circ (\mu_{H}^{2}\otimes \lambda_{H}^{1}))\otimes (\mu_{H}^{1}\circ (\mu_{H}^{2}\otimes \lambda_{H}^{1})))\circ (H \otimes H\otimes c_{H,H}\otimes H\otimes H)\circ (H\otimes c_{H,H}\otimes (c_{H,H}\circ c_{H,H})\otimes H)$
\item[ ]$\hspace{0.38cm}\circ (\delta_{H}\otimes c_{H,H}\otimes c_{H,H})\circ (H\otimes \delta_{H\otimes H})\circ (\delta_{H}\otimes H) ${\scriptsize (by the naturality of $c$)}
\item [ ]$= ((\mu_{H}^{1}\circ (\mu_{H}^{2}\otimes \lambda_{H}^{1}))\otimes (\mu_{H}^{1}\circ (\mu_{H}^{2}\otimes \lambda_{H}^{1})))\circ (H \otimes H\otimes c_{H,H}\otimes H\otimes H)\circ (H\otimes c_{H,H}\otimes H\otimes H\otimes H)$
\item[ ]$\hspace{0.38cm}\circ (\delta_{H}\otimes c_{H,H}\otimes c_{H,H})\circ (H\otimes \delta_{H\otimes H})\circ (\delta_{H}\otimes H) ${\scriptsize (by (\ref{ccb}))}
\item [ ]$= ((\mu_{H}^{1}\circ (\mu_{H}^{2}\otimes \lambda_{H}^{1}))\otimes (\mu_{H}^{1}\circ (\mu_{H}^{2}\otimes \lambda_{H}^{1})))\circ (H\otimes ((c_{H,H}\otimes H)\circ (H\otimes c_{H,H})\circ ((c_{H,H}\circ \delta_{H})\otimes H))$
\item[ ]$\hspace{0.38cm}\otimes c_{H,H}) \circ (\delta_{H}\otimes c_{H,H}\otimes H)\circ (\delta_{H}\otimes \delta_{H})$ {\scriptsize (by the naturality of $c$ and the coassociativity condition)}
\item [ ]$=(\Gamma_{H_{1}}^{\prime}\otimes \Gamma_{H_{1}}^{\prime})\circ \delta_{H\otimes H} $ {\scriptsize (by the naturality of $c$
 and the cocommutativity and coassociativity conditions)}

\end{itemize}
\end{proof}

Following \cite{RGON} we recall the notion of left module for a Hopf brace.

\begin{definition}
\label{l-mod}
{\rm Let ${\mathbb H}$ be a Hopf brace. A left ${\mathbb H}$-module is a triple $(M,\psi_{M}^{1}, \psi_{M}^{2})$, where $(M,\psi_{M}^{1})$ is a left $H_{1}$-module, $(M, \psi_{M}^{2})$  is a left $H_{2}$-module and the following identity 
\begin{equation}
\label{mod-l1}
\psi_{M}^{2}\co (H\ot \psi_{M}^{1})=\psi_{M}^1\co (\mu_{H}^{2}\ot \Gamma_{M})\co (H\ot c_{H,H}\ot M)\co (\delta_{H}\ot H\ot M)
\end{equation}
holds, where 
$$\Gamma_{M}=\psi_{M}^{1}\co (\lambda_{H}^1\ot \psi_{M}^{2})\co (\delta_{H}\ot M).$$
		
Given two left ${\mathbb H}$-modules  $(M,\psi_{M}^1, \psi_{M}^{2})$  and  $(N,\psi_{N}^1, \psi_{N}^{2})$, a morphism $f:M\rightarrow N$  is called a morphism of left ${\mathbb H}$-modules if  $f$ is a morphism of left $H_{1}$-modules and left $H_{2}$-modules. Left ${\mathbb H}$-modules  with morphisms of left ${\mathbb H}$-modules  form a category which we denote by $\;_{\mathbb H}${\sf Mod}. 
}
\end{definition}

\begin{example}
\label{exatriv}
Let ${\mathbb H}$ be a Hopf brace. The triple $(H, \mu_{H}^{1}, \mu_{H}^{2})$ is an example of  left ${\mathbb H}$-module. Also, if $K$ is the unit object of ${\sf C}$, $(K, \psi_{K}^1=\varepsilon_{H}, \psi_{K}^2=\varepsilon_{H})$ is a left ${\mathbb H}$-module called the trivial module.  Moreover, $(H, \psi_{H}^1=\varepsilon_{H}\otimes H, \psi_{H}^2=\mu_{H}^2)$ is an object in $\;_{\mathbb H}${\sf Mod} and we have a functor  $${\sf T}:\;_{\sf H_2}{\sf Mod}\;\rightarrow \; _{\mathbb H}{\hspace{-0.01cm}}{\sf Mod}$$ defined on objects by ${\sf T}((M, \psi_{M}))=(M, \psi_{M}^1=\varepsilon_{H}\otimes M, \psi_{M}^2=\psi_{M})$ and by the identity on morphisms. In this setting, there exists a forgetful functor 
$${\sf W}:{\hspace{0.01cm}}_{\mathbb H}{\hspace{-0.01cm}}{\sf Mod}\rightarrow{\hspace{0.01cm}}_{\sf H_2}{\sf Mod}$$ defined on objects by ${\sf W}((M, \psi_{M}^1, \psi_{M}^{2}))=(M, \psi_{M}^2)$ and by the identity on morphisms. Obviously, ${\sf W}\circ {\sf T}={\sf id}_{_{\sf H_2}{\sf Mod}}.$
		
Let $H=(H, \eta_{H}, \mu_{H},  \varepsilon_{H}, \delta_{H}, \lambda_{H})$ be a Hopf algebra. Then $(H,\mu_{H}, \mu_{H})$ is an example of left ${\mathbb H}$-module for the Hopf brace ${\mathbb H}_{triv}$. Also, if $(M, \psi_{M})$ is a left $H$-module, the triple $(M, \psi_{M}, \psi_{M})$ is a left ${\mathbb H}_{triv}$-module. Then, we have a functor $${\sf J}:\;_{\sf H}{\sf Mod}\;\rightarrow \; _{\mathbb{H}_{triv}}{\hspace{-0.01cm}}{\sf Mod}$$ defined on objects by ${\sf J}((M, \psi_{M}))=(M, \psi_{M}, \psi_{M})$ and by the identity on morphisms. Also, there exists a forgetful functor $${\sf U}:\;_{\mathbb{H}_{triv}}{\sf Mod}\;\rightarrow {\hspace{-0.01cm}}_{\sf H}{\sf Mod}$$ defined on objects by ${\sf U}((M, \psi_{M}^1, \psi_{M}^{2}))=(M, \psi_{M}^1)$ and by the identity on morphisms. Then, ${\sf U}\circ {\sf J}={\sf id}_{_{\sf H}{\sf Mod}}$ holds trivially.
\end{example}

\begin{remark}
\label{chino}
As was pointed in \cite{RGON}, Definition \ref{l-mod} is weaker than the one introduced by H. Zhu in \cite{Zhu}. For this author, if ${\mathbb H}$ is a Hopf brace, a left ${\mathbb H}$-module is a triple $(M,\psi_{M}^1, \psi_{M}^2)$, where $(M,\psi_{M}^1)$ is a left $H_{1}$-module, $(M, \psi_{M}^2)$  is a left $H_{2}$-module, and the equalities (\ref{mod-l1}) and 
\begin{equation}
\label{ch}
(\psi_{M}^{2}\ot H)\co (H\ot c_{H,M})\co (\delta_{H}\ot M)=(\psi_{M}^1\ot H)\co (H\ot c_{H,M})\co (\delta_{H}\ot \Gamma_{M})\co (\delta_{H}\ot M)
\end{equation}
hold (see \cite[Definition 3.1, Lemma 3.2]{Zhu}). Thus, for an arbitrary Hopf brace ${\mathbb H}$, a left ${\mathbb H}$-module in the sense of Zhu is a left ${\mathbb H}$-module in our sense. Moreover, if ${\mathbb H}$ is cocommutative, (\ref{ch}) hold for any left ${\mathbb H}$-module as in Definition \ref{l-mod}. As a consequence, in the cocommutative setting, \cite[Definition 3.1]{Zhu} and Definition \ref{l-mod}  are equivalent. 

For every Hopf algebra $H$, the first example of a left module over $H$ is the algebra $H$ taking as action the product $\mu_{H}$. In the case that we intend to introduce a coherent definition of left module for a Hopf brace ${\mathbb H}$, the same should still be true for  ${\mathbb H}$. If we work with Definition \ref{l-mod}, trivially, $(H,\mu_{H}^1, \mu_{H}^2)$ is a left ${\mathbb H}$-module but if we work with the definition introduced by Zhu the triple $(H,\mu_{H}^1, \mu_{H}^2)$ is a left ${\mathbb H}$-module iff
\begin{equation}
\label{ch1}
(\mu_{H}^{2}\ot H)\co (H\ot c_{H,H})\co (\delta_{H}\ot H)=(\mu_{H}^1\ot H)\co (H\ot c_{H,H})\co (\delta_{H}\ot \Gamma_{H_1})\co (\delta_{H}\ot H)
\end{equation}
holds. If equality (\ref{ch1}) is satisfied, the following identity 
\begin{equation}
\label{ch11}
(\mu_{H}^{1}\otimes H)\circ (H\otimes ((\Gamma_{H_1}\otimes H)\circ (H\otimes c_{H,H})\circ (\delta_{H}\otimes H)))\circ (\delta_{H}\otimes H)
\end{equation}
$$=(\mu_{H}^{1}\otimes H)\circ (H\otimes ((\Gamma_{H_1}\otimes H)\circ (H\otimes c_{H,H})\circ ((c_{H,H}\circ\delta_{H})\otimes H)))\circ (\delta_{H}\otimes H)$$
holds because of (\ref{eb2}), the naturality of $c$ and the coassociativity of $\delta_{H}$. Then, composing in (\ref{ch11}) with $((\lambda_{H}^1\otimes H)\circ \delta_{H})\otimes H$ on the right and with $\mu_{H}^{1}\otimes H$ on the left we obtain the identity 
$$
(\mu_{H}^{1}\otimes H)\circ ((\lambda_{H}^1\ast id_{H})\otimes ((\Gamma_{H_1}\otimes H)\circ (H\otimes c_{H,H})\circ (\delta_{H}\otimes H)))\circ (\delta_{H}\otimes H)$$
$$=(\mu_{H}^{1}\otimes H)\circ ((\lambda_{H}^1\ast id_{H})\otimes ((\Gamma_{H_1}\otimes H)\circ (H\otimes c_{H,H})\circ ((c_{H,H}\circ\delta_{H})\otimes H)))\circ (\delta_{H}\otimes H).$$

This implies that 
$$
%\label{clcc}
(\Gamma_{H_1}\otimes H)\circ (H\otimes c_{H,H})\circ (\delta_{H}\otimes H)=(\Gamma_{H_1}\otimes H)\circ (H\otimes c_{H,H})\circ ((c_{H,H}\circ\delta_{H})\otimes H)
$$
holds. Therefore, if $(H,\mu_{H}^1, \mu_{H}^2)$ is a left ${\mathbb H}$-module in the sense of Zhu and the category {\sf C} is symmetric (for example, the category of vector spaces over a field ${\mathbb K}$), we have that $(H, \Gamma_{H_1})$ is in the cocommutativity class of $H$ (see \cite{CCH} for the definition) and, obviously, this does not always have to happen. In other words, under certain circumstances, for example, the lack of cocommutativity, the category of left modules over a Hopf brace introduced by Zhu could have as its only objects the base object of the category {\sf C} and its tensor products with the trivial action. 

\end{remark}

\begin{remark}
%\label{yo} 
Using the naturality of $c$ and the coassociativity of $\delta_{H}$, it is easy to  show that (\ref{mod-l1}) is equivalent to 
\begin{equation}
\label{mod-l1p}
\psi_{M}^{2}\co (H\ot \psi_{M}^{1})=\psi_{M}^1\co (\Gamma_{H_1}^{\prime}\ot \psi_{M}^{2})\co (H\ot c_{H,H}\ot M)\co (\delta_{H}\ot H\ot M).
\end{equation}
\end{remark}

\begin{lemma}
%\label{lemG}
Let ${\mathbb H}$ be a Hopf brace and let $(M,\psi_{M}^1, \psi_{M}^2)$ be a left ${\mathbb H}$-module. Then, the following equality holds:
\begin{equation}
\label{GMH1}
\Gamma_{M}\circ (H\otimes \psi_{M}^1)=\psi_{M}^1\circ (\Gamma_{H_1}\otimes \Gamma_{M})\circ (H\otimes c_{H,H}\otimes M)\circ (\delta_{H}\otimes H\otimes M).
\end{equation}
Also, $(M, \Gamma_{M})$ is a left $H_{2}$-module.
\end{lemma}

\begin{proof} Let $(M,\psi_{M}^1, \psi_{M}^2)$ be a left ${\mathbb H}$-module.Then the equality (\ref{GMH1}) follows from:
\begin{itemize}
\item[ ]$\hspace{0.38cm}\Gamma_{M}\circ (H\otimes \psi_{M}^1) $
\item [ ]$=\psi_{M}^1\circ (\lambda_{H}^1\otimes (\psi_{M}^1\co (\mu_{H}^{2}\ot \Gamma_{M})\co (H\ot c_{H,H}\ot M)\co (\delta_{H}\ot H\ot M)))\circ (\delta_{H}\otimes H\otimes H)  $ {\scriptsize (by (\ref{mod-l1}))}
\item [ ]$= \psi_{M}^1\circ (\Gamma_{H_1}\otimes \Gamma_{M})\circ (H\otimes c_{H,H}\otimes M)\circ (\delta_{H}\otimes H\otimes M)${\scriptsize (by the coassociativity of $\delta_{H}$ and  the condition }
\item[ ]$\hspace{0.38cm}$ {\scriptsize   of left $H_{1}$-module for $M$).}
\end{itemize}

On the other hand, trivially $\Gamma_{M}\circ (\eta_{H}\otimes M)=id_{M}$ and 
\begin{itemize}
\item[ ]$\hspace{0.38cm} \Gamma_{M}\circ (H\otimes \Gamma_{M}) $
\item [ ]$=\Gamma_{M}\circ (H\otimes (\psi_{M}^1\circ (\lambda_{H}^1\otimes \psi_{M}^2)\circ (\delta_{H}\otimes M))) $ {\scriptsize (by definition of $\Gamma_{M}$)}
\item [ ]$=\psi_{M}^1\circ (\Gamma_{H_1}\otimes \Gamma_{M})\circ (H\otimes c_{H,H}\otimes M)\circ (\delta_{H}\otimes ((\lambda_{H}^1\otimes \psi_{M}^2)\circ (\delta_{H}\otimes M))) $ {\scriptsize (by (\ref{GMH1}))}
\item [ ]$=\psi_{M}^1\circ ((\Gamma_{H_1}\circ (H\otimes \lambda_{H}^{1}))\otimes \Gamma_{M})\circ (H\otimes c_{H,H}\otimes M)\circ (\delta_{H}\otimes ((H\otimes \psi_{M}^2)\circ (\delta_{H}\otimes M)))$ {\scriptsize (by }
\item[ ]$\hspace{0.38cm}$ {\scriptsize naturality of $c$)}
\item [ ]$=\psi_{M}^1\circ ((\mu_{H}^{1}\circ ((\lambda_{H}^1\circ \mu_{H}^{2})\otimes H)\circ (H\otimes c_{H,H}) \circ (\delta_{H}\otimes H))\otimes (\psi_{M}^1\circ (\lambda_{H}^1\otimes \psi_{M}^2)\circ (\delta_{H}\otimes M)))$
\item[ ]$\hspace{0.38cm}\circ (H\otimes c_{H,H}\otimes M)\circ (\delta_{H}\otimes ((H\otimes \psi_{M}^2)\circ (\delta_{H}\otimes M))) $ {\scriptsize (by (\ref{agv1}) and the definition of $\Gamma_{M}$)}
\item [ ]$= \psi_{M}^1\circ ((\mu_{H}^1\circ ((\lambda_{H}^1\circ \mu_{H}^2)\otimes (\mu_{H}^1\circ (H\otimes \lambda_{H}^1)))\otimes (\psi_{M}^2\circ (\mu_{H}^2\otimes M)))\circ (H\otimes c_{H,H}\otimes \delta_{H}\otimes H\otimes M)$
\item[ ]$\hspace{0.38cm}\circ (\delta_{H}\otimes c_{H,H}\otimes H\otimes M)\circ (\delta_{H}\otimes \delta_{H}\otimes M)$ {\scriptsize (by the condition of left $H_{1}$ and $H_{2}$-module for $M$ and the}
\item[ ]$\hspace{0.38cm}${\scriptsize associativity of $\mu_{H}^1$)}
\item [ ]$=\psi_{M}^1\circ ((\mu_{H}^1\circ ((\lambda_{H}^1\circ \mu_{H}^2)\otimes (\mu_{H}^1\circ (H\otimes \lambda_{H}^1)))\otimes (\psi_{M}^2\circ (\mu_{H}^2\otimes M)))\circ (((H\otimes c_{H,H}\otimes H\otimes H\otimes H)$
\item[ ]$\hspace{0.38cm}\circ (\delta_{H}\otimes c_{H,H}\otimes H\otimes H)\circ (H\otimes \delta_{H\otimes H})\circ (\delta_{H}\otimes H))\otimes M) $ {\scriptsize (by naturality of $c$)}
\item [ ]$=\psi_{M}^1\circ ((\mu_{H}^1\circ ((\lambda_{H}^{1}\circ \mu_{H}^2)\otimes (id_{H}\ast \lambda_{H}^1))\circ (H\otimes c_{H,H})\circ (\delta_{H}\otimes H))\otimes (\psi_{M}^2\circ (\mu_{H}^2\otimes M)))\circ (\delta_{H\otimes H}\otimes M) $ 
\item[ ]$\hspace{0.38cm}$ {\scriptsize (by naturality of $c$ and coassociativity of $\delta_{H}$)}
\item [ ]$= \psi_{M}^1\circ (\lambda_{H}^1\otimes \psi_{M}^2) \circ (((\mu_{H}^2\otimes \mu_{H}^2)\circ \delta_{H\otimes H})\otimes M)$ {\scriptsize (by (\ref{antipode}) and unit and counit properties)}
\item [ ]$= \Gamma_{M}\circ (\mu_{H}^2\otimes M)$ {\scriptsize (by the condition of coalgebra morphism for $\mu_{H}^2$)}
\end{itemize}

Therefore, $(M,\Gamma_{M})$ is a left $H_{2}$-module.
\end{proof}

\begin{theorem}\label{th-mon-cat2}
 Let's assume that {\sf C} is symmetric with natural isomorphism of symmetry $c$. Let ${\mathbb H}$ be  a cocommutative Hopf brace in {\sf C}. Then the category of left modules over ${\mathbb H}$ is symmetric monoidal with unit object the trivial left module over ${\mathbb H}$.
\end{theorem}

\begin{proof} Let $(M, \psi_{M}^1, \psi_{M}^2)$, $(N, \psi_{N}^1, \psi_{N}^2)$ be objects in $_{\mathbb H}{\sf Mod}$. The tensor product is defined by $(M\otimes N, \psi_{M\otimes N}^1, \psi_{M\otimes N}^2)$ where $\psi_{M\otimes N}^1$ and $\psi_{M\otimes N}^2$ are the corresponding module tensor structures. In fact, $(M\otimes N, \psi_{M\otimes N}^1)$ is a left $H_{1}$-module, $(M\otimes N, \psi_{M\otimes N}^2)$ is a left $H_{2}$-module due to the monoidal character of the category of modules over a Hopf algebra. On the other hand, the identity  
\begin{equation}
\label{gg}
\Gamma_{M\otimes N}=(\Gamma_{M}\otimes \Gamma_{N})\circ (H\otimes c_{H,M}\otimes N)\circ (\delta_{H}\otimes M\otimes N)
\end{equation}
holds because
\begin{itemize}
\item[ ]$\hspace{0.38cm} \Gamma_{M\otimes N}$
\item [ ]$= ((\psi_{M}^1\circ (H\otimes \psi_{M}^2))\otimes (\psi_{N}^1\circ (H\otimes \psi_{N}^2)))\circ (H\otimes ((H\otimes c_{H,M}\otimes H)\circ (c_{H,H}\otimes c_{H,M}))\otimes N)$
\item[ ]$\hspace{0.38cm}\circ (((((\lambda_{H}^1\otimes \lambda_{H}^1)\circ \delta_{H})\otimes \delta_{H})\circ \delta_{H})\otimes M\otimes N)$ {\scriptsize (by (\ref{a-antip}), the cocommutativity of $\delta_{H}$ and the naturality of}
\item[ ]$\hspace{0.38cm}${\scriptsize  $c$)}
\item [ ]$= ((\psi_{M}^1\circ (\lambda_{H}^1\otimes \psi_{M}^2))\otimes (\psi_{N}^1\circ (\lambda_{H}^1\otimes \psi_{N}^2)))\circ (H\otimes ((H\otimes c_{H,M}\otimes H)\circ (H\otimes H\otimes c_{H,M}))\otimes N)$
\item[ ]$\hspace{0.38cm}\circ (((H\otimes (c_{H,H}\circ\delta_{H})\otimes H)\circ (\delta_{H}\otimes H)\circ \delta_{H})\otimes M\otimes N)${\scriptsize (by coassociativity of $\delta_{H}$ and the naturality }
\item[ ]$\hspace{0.38cm}${\scriptsize  of $c$)}
\item [ ]$= (\Gamma_{M}\otimes \Gamma_{N})\circ (H\otimes c_{H,M}\otimes N)\circ (\delta_{H}\otimes M\otimes N)$ {\scriptsize (by the coassociativity and cocommutativity of $\delta_{H}$ and}
\item[ ]$\hspace{0.38cm}${\scriptsize  the naturality of $c$)}
\end{itemize}

Then
\begin{itemize}
\item[ ]$\hspace{0.38cm}\psi_{M\otimes N}^1\co (\mu_{H}^{2}\ot \Gamma_{M\otimes N})\co (H\ot c_{H,H}\ot M\otimes N)\co (\delta_{H}\ot H\ot M\otimes N) $
\item [ ]$=(\psi_{M}^1\otimes \psi_{N}^1)\circ (H\otimes c_{H,M}\otimes N)\circ (((\mu_{H}^2\otimes \mu_{H}^2)\circ \delta_{H\otimes H})\otimes ((\Gamma_{M}\otimes \Gamma_{N})\circ (H\otimes c_{H,M}\otimes N)$
\item[ ]$\hspace{0.38cm}\circ (\delta_{H}\otimes M\otimes N)))\circ (((H\otimes c_{H,H})\circ (\delta_{H}\otimes H))\otimes M\otimes N)$ {\scriptsize (by (\ref{gg}) and the condition of coalgebra}
\item[ ]$\hspace{0.38cm}${\scriptsize  morphism of $\mu_{H}^2$)}
\item [ ]$= ((\psi_{M}^1\circ (\mu_{H}^2\otimes \Gamma_{M})\circ (H\otimes c_{H,H}\otimes M))\otimes (\psi_{N}^1\circ (\mu_{H}^2\otimes \Gamma_{N})\circ (H\otimes c_{H,H}\otimes N)))$
\item[ ]$\hspace{0.38cm}\circ (H\otimes ((H\otimes H\otimes c_{H,M}\otimes H\otimes H)\circ (H\otimes c_{H,H}\otimes c_{H,M}\otimes H)\circ ((c_{H,H}\circ \delta_{H})\otimes c_{H,H}\otimes c_{H,M}))\otimes N)$
\item[ ]$\hspace{0.38cm}\circ (((\delta_{H}\otimes H)\circ \delta_{H})\otimes \delta_{H}\otimes M\otimes N) $ {\scriptsize (by the coassociativity of $\delta_{H}$, the naturality of $c$ and $c_{H,H}\circ c_{H,H}=id_{H\otimes H}$)}
\item [ ]$= ((\psi_{M}^2\circ (H\otimes \psi_{M}^1) )\otimes (\psi_{N}^1\circ (\mu_{H}^2\otimes \Gamma_{N})\circ (H\otimes c_{H,H}\otimes N)))$
\item[ ]$\hspace{0.38cm}\circ (H\otimes ((H\otimes c_{H,M}\otimes H\otimes H)\circ (c_{H,H}\otimes c_{H,M}\otimes H)\circ (H\otimes c_{H,H}\otimes c_{H,M}))\otimes N)$
\item[ ]$\hspace{0.38cm}\circ (((\delta_{H}\otimes H)\circ \delta_{H})\otimes \delta_{H}\otimes M\otimes N) $ {\scriptsize (by (\ref{mod-l1}) for $M$, the coassociativity and cocommutativity of $\delta_{H}$)}
\item [ ]$=((\psi_{M}^2\circ (H\otimes \psi_{M}^1))\otimes (\psi_{N}^2\circ (H\otimes \psi_{N}^1)))\circ  (H\otimes ((H\otimes c_{H,M}\otimes H)\circ (c_{H,H}\otimes c_{H,M}))\otimes N)$
\item[ ]$\hspace{0.38cm}\circ (\delta_{H}\otimes \delta_{H}\otimes M\otimes N) $ {\scriptsize (by the coassociativity of $\delta_{H}$, the naturality of $c$ and (\ref{mod-l1}) for $N$)}
\item [ ]$=\psi_{M\otimes N}^{2}\co (H\ot \psi_{M\otimes N}^{1}) $ {\scriptsize (by the naturality of $c$)}.
\end{itemize}

The unit object is $(K, \psi^{1}_{K}=\varepsilon_{H}, \psi^{2}_{K}=\varepsilon_{H})$ and the natural isomorphism of symmetry is $c$ because, if $H$ is cocommutative and ${\sf C}$ is symmetric, $c$ is a morphism of left modules over ${\mathbb H}$. 	
\end{proof}

\section{Yetter-Drinfeld modules for Hopf braces}

The first goal of this section is to introduce a suitable notion of Yetter-Drinfeld module for a Hopf brace that can be useful in the study of Hopf brace projections.  To do it, previously it is necessary to define some intermediate objects called weak left Yetter-Drinfeld modules.
		
\begin{definition}
\label{YDbrace}
Let ${\mathbb H}$ be a Hopf brace in ${\sf C}$. A weak left Yetter-Drinfeld module over ${\mathbb H}$ is a quadruple $(M,\psi_{M}^{1}, \psi_{M}^{2}, \rho_{M})$ such that 
\begin{itemize}
\item[(i)] $(M,\psi_{M}^{1}, \psi_{M}^{2})\in \;_{\mathbb H}{\sf Mod}$.
\item[(ii)] $(M,\psi_{M}^{1},  \rho_{M})\in \:^{\sf H_1}_{\sf H_1}{\sf YD}$.
\item[(iii)] $(M,\psi_{M}^{2},  \rho_{M})\in \:^{\sf H_2}_{\sf H_2}{\sf YD}$.
\item[(iv)] The following equality 
$$(\mu_{H}^{1}\otimes M)\circ (H\otimes c_{M,H})\circ (\rho_{M}\otimes H)=
(\mu_{H}^{2}\otimes M)\circ (H\otimes c_{M,H})\circ (\rho_{M}\otimes H)$$
holds.
\end{itemize}

With the obvious morphisms, i.e., morphisms of left ${\mathbb H}$-modules and left $H$-comodules, weak left Yetter-Drinfeld modules over ${\mathbb H}$ form a category that we will denote by $\:^{\mathbb H}_{\mathbb H}{\sf WYD}$.

\end{definition}

\begin{theorem}
\label{WYD1}
Let's assume that {\sf C} is symmetric with natural isomorphism of symmetry $c$. Let ${\mathbb H}$ be  a cocommutative Hopf brace in {\sf C}. Then the category  $\:^{\mathbb H}_{\mathbb H}{\sf WYD}$ is monoidal.
\end{theorem}

\begin{proof}
Let $(M,\psi_{M}^{1}, \psi_{M}^{2}, \rho_{M})$, $(N,\psi_{N}^{1}, \psi_{N}^{2}, \rho_{N})$ be objects in $\:^{\mathbb H}_{\mathbb H}{\sf WYD}$. Then the tensor product is defined by 
$$(M\otimes N, \psi_{M\otimes N}^{1},  \psi_{M\otimes N}^{2},  \rho_{M\otimes N}).$$
where, by (iv) of Definition \ref{YDbrace}, 
$$ \rho_{M\otimes N}=(\mu_{H}^1\otimes M\otimes N)\circ (H\otimes c_{M,H}\otimes N)\circ (\rho_{M}\otimes \rho_{N})=(\mu_{H}^2\otimes M\otimes N)\circ (H\otimes c_{M,H}\otimes N)\circ (\rho_{M}\otimes \rho_{N}).$$

By Theorem \ref{th-mon-cat2} we have that $(M\otimes N, \psi_{M\otimes N}^{1},  \psi_{M\otimes N}^{2})$ is an object in $\;_{\mathbb H}{\sf Mod}$. Moreover by the monoidal structure of the categories of Yetter-Drinfeld modules associated to a Hopf algebra we have that   
$(M\otimes N, \psi_{M\otimes N}^{1}, \rho_{M\otimes N})$ belongs to $ \:^{\sf H_1}_{\sf H_1}{\sf YD}$ and $(M\otimes N, \psi_{M\otimes N}^{2}, \rho_{M\otimes N})\in \:^{\sf H_2}_{\sf H_2}{\sf YD}$. Finally, (iv) of Definition (\ref{YDbrace}) also holds because
\begin{itemize}
	\item[ ]$\hspace{0.38cm} (\mu_{H}^{1}\otimes M)\circ (H\otimes c_{M\otimes N,H})\circ (\rho_{M\otimes N}\otimes H) $
	\item [ ]$=(\mu_{H}^{1}\otimes M\otimes N)\circ (H\otimes c_{M,H}\otimes N)\circ (\rho_{M}\otimes ((\mu_{H}^{1}\otimes N)\circ (H\otimes c_{N,H})\circ (\rho_{N}\otimes H))) $ {\scriptsize (by associativity }
	\item[ ]$\hspace{0.38cm}${\scriptsize of $\mu_{H}^{1}$ and naturality of $c$)}
	\item [ ]$=(\mu_{H}^{2}\otimes M\otimes N)\circ (H\otimes c_{M,H}\otimes N)\circ (\rho_{M}\otimes ((\mu_{H}^{2}\otimes N)\circ (H\otimes c_{N,H})\circ (\rho_{N}\otimes H))) $ {\scriptsize (by (iv) of  }
	\item[ ]$\hspace{0.38cm}${\scriptsize Definition \ref{YDbrace} for $M$ and $N$)}
	\item [ ]$=(\mu_{H}^{2}\otimes M)\circ (H\otimes c_{M\otimes N,H})\circ (\rho_{M\otimes N}\otimes H)  $ {\scriptsize (by associativity of $\mu_{H}^{2}$ and naturality of $c$).}
\end{itemize}

Finally, it is easy to show that  the unit object is $(K,\psi^{1}_{K}=\varepsilon_{H},\psi^{2}_{K}=\varepsilon_{H}, \rho_{K}=\eta_{H})$. 

\end{proof}

\begin{definition}
\label{YDbrace-new}
Let ${\mathbb H}$ be a Hopf brace in ${\sf C}$. We define the category of left Yetter-Drinfeld modules over ${\mathbb H}$, denoted by $\:^{\mathbb H}_{\mathbb H}{\sf YD}$, as the full subcategory of $\:^{\mathbb H}_{\mathbb H}{\sf WYD}$ whose objects $(M,\psi_{M}^{1}, \psi_{M}^{2}, \rho_{M})$ satisfy that 
$$
t_{M,N}^2=(\psi_{N}^{2}\otimes M)\circ (H\otimes c_{M,N})\circ (\rho_{M}\otimes N)
$$
is a morphism of left $H_{1}$-modules for all $(N,\psi_{N}^{1}, \psi_{N}^{2}, \rho_{N}) \in \:^{\mathbb H}_{\mathbb H}{\sf WYD}$
\end{definition}

\begin{remark}
1) Note that, under the conditions of the previous definition,  $t_{M,N}^2$ is a morphism of left $H_{2}$-modules because  $(M,\psi_{M}^{2}, \rho_{M})$ and $(N,\psi_{N}^{2}, \rho_{N})$ are left Yetter-Drinfeld modules over $H_2$.   Moreover, if the antipode of $H_{2}$ is an isomorphism,  $t_{M,N}^2$ is the braiding of the category $\:^{H_2}_{H_2}{\sf YD}$.

2) Let ${\mathbb H}$ be a cocommutative Hopf brace in ${\sf C}$. Then, by Example \ref{exatriv}, (\ref{ccb}), the unit-counit properties and the naturality of $c$, we have that 
$$(H, \psi_{H}^{1}=\varepsilon_{H}\otimes H, \psi_{H}^{2}=\mu_{H}^{2}, \rho_{H}=\eta_{H}\otimes H)$$
is an object in $\:^{\mathbb H}_{\mathbb H}{\sf YD}$.

Then, if  {\sf C} is symmetric, using similar arguments to the previous paragraph we can assure that there exists a functor 
$${\sf S}:\;_{\sf H_2}{\sf Mod}\;\rightarrow  \;^{\mathbb H}_{\mathbb H}{\sf YD}$$ defined on objects by ${\sf S}((M, \psi_{M}))=(M, \psi_{M}^1=\varepsilon_{H}\otimes M, \psi_{M}^2=\psi_{M}, \rho_{M}=\eta_{H}\otimes M)$ and by the identity on morphisms. Also, as in Example \ref{exatriv}, there exists a forgetful functor $${\sf V}:\;^{\mathbb H}_{\mathbb H}{\sf YD}\rightarrow{\hspace{0.01cm}}_{\sf H_2}{\sf Mod}$$ defined on objects by ${\sf V}((M, \psi_{M}^1, \psi_{M}^{2}, \rho_{M}))=(M, \psi_{M}^2)$ and by the identity on morphisms. Obviously, ${\sf V}\circ {\sf S}={\sf id}_{_{\sf H_2}{\sf Mod}}.$

3) Assume that ${\mathbb H}$ is a cocommutative Hopf brace in ${\sf C}$.  From the previous point we know that  $$(H, \psi_{H}^{1}=\varepsilon_{H}\otimes H, \psi_{H}^{2}=\mu_{H}^{2}, \rho_{H}=\eta_{H}\otimes H)$$ is an object in the category $\:^{\mathbb H}_{\mathbb H}{\sf YD}$. Let $(M,\psi_{M}^{1}, \psi_{M}^{2}, \rho_{M})$ be an object in the same category. Thus, by definition, 
$$t^{2}_{M,H}=(\mu_{H}^{2}\otimes M)\circ (H\otimes c_{M,H})\circ (\rho_{M}\otimes H)$$ 
is a morphism of left $H_{1}$-modules. This fact is equivalent to the following equality 
$$\psi_{H\otimes M}^1\circ (H\otimes t^{2}_{M,H})=t^{2}_{M,H}\circ \psi_{M\otimes H}^1$$
and, using the naturality of $c$ and the properties of the counit, we can prove that the previous identity is equivalent to
\begin{equation}
\label{diss1}
(\mu_{H}^{2}\otimes H)\circ (H\otimes c_{M,H})\circ (((H\otimes \psi_{M}^1)\circ (c_{H,H}\otimes M)\circ (H\otimes \rho_{M}))\otimes H)
\end{equation}
$$
=(\mu_{H}^{2}\otimes H)\circ (H\otimes c_{M,H})\circ ((\rho_{M}\circ \psi_{M}^1)\otimes H).
$$

Therefore, composing on the right of (\ref{diss1}) with $H\otimes M\otimes \eta_{H}$, we obtain that 
\begin{equation}
\label{diss11}
(H\otimes \psi_{M}^1)\circ (c_{H,H}\otimes M)\circ (H\otimes \rho_{M})=\rho_{M}\circ \psi_{M}^1,
\end{equation}
or, in other words, $(M, \psi_{M}^1, \rho_{M})$ is a Long dimodule. This category was introduced by Long in \cite{Long} to study the Brauer group of $H$-dimodule algebras for a commutative and cocommutative Hopf algebra $H$. Later, the notion was extended by considering two arbitrary Hopf algebras $H$ and $B$, introducing the category of left-left $H$-$B$-Long dimodules, denoted by ${\sf \;_H^B  Long}$. In this category the objects are triples $(M, \varphi_M, \rho_M)$ such that $(M, \varphi_M)$ is  a left $H$-module and $(M,  \rho_M)$ is a left $B$-comodule satisfying the axiom
\begin{equation}
\label{Long-dm-1}
\rho_M\circ \varphi_M=(B\otimes   \varphi_M)\circ (c_{H,B}\otimes   M)\circ (H\otimes   \rho_M), 
\end{equation}

The morphisms in ${\sf \;_H^B  Long}$ are morphisms  of left  $H$-modules and left $B$-comodules.

Then, in our setting, taking into account that (\ref{diss11}) is exactly (\ref{Long-dm-1}) for $H=B=H_{1}$, we have a functor $${\sf L}:\;^{\mathbb H}_{\mathbb H}{\sf YD}\;\rightarrow  {\sf \;_{H_{1}}^{H_{1}}  Long}$$ defined on objects by ${\sf L}((M,\psi_{M}^{1}, \psi_{M}^{2}, \rho_{M}))=(M, \psi_{M}^1, \rho_{M})$ and by the identity on morphisms.

\end{remark}

\begin{example}
\label{ex-YD-Zhu}
In \cite{Zhu}, Haixin Zhu gives a definition of a Yetter-Drinfeld category making use of the notion of the  Hopf brace subbialgebra. Recall that in \cite{Zhu} the base category is a category of vector spaces over a fixed field ${\mathbb K}$. In this context, a subbialgebra $H'$ of a Hopf brace  $\mathbb{H}$ (\cite[Definition 4.1]{Zhu}) is a subbialgebra $H'$ of $H_{1}$ (i.e., if $j_{H^{\prime}}:H'\rightarrow H$ is the inclusion morphism, $j_{H^{\prime}}$ is an algebra-coalgebra morphism) such that 
\begin{equation}
\label{echin}
\mu_{H}^1\circ (j_{H^{\prime}}\otimes H) =  \mu_{H}^2\circ (j_{H^{\prime}}\otimes H). 
\end{equation}

Note, by (\ref{echin}), $H'$ also is a subbialgebra of $H_{2}$ because 
$$\mu_{H}^{2}\circ (j_{H^{\prime}}\otimes j_{H^{\prime}})=\mu_{H}^{1}\circ (j_{H^{\prime}}\otimes j_{H^{\prime}})=j_{H^{\prime}}\circ \mu_{H'}.$$

Once the subbialgebra is defined, the author considers the so-called compatible modules (see \cite[Definition 4.5]{Zhu}), that are modules $(M,\psi_{M}^1, \psi_{M}^2)$ over $\mathbb H$,  in the sense of Remark \ref{chino}, satisfying
\[\psi_{M}^1\circ (j_{H^{\prime}}\otimes M) = \psi_{M}^2\circ (j_{H^{\prime}}\otimes M). \]  

Now, the objects of the category of Yetter-Drinfeld modules $_\mathbb{H}^{\sf H'}{\sf YD}$ are compatible left modules over $\mathbb H$ that have a comodule structure $\rho_M:M\to H\ot M$ such that (see \cite[Definition 4.7]{Zhu}):
\begin{enumerate}
\item[(i)] $\rho_M(m) \in H'\ot M, \;\forall \; m\in M$,
\item[(ii)] $(M, \psi_{M}^1, \rho_M)\in \:^{\sf H_1}_{\sf H_1}{\sf YD}$,
\item[(iii)] $(M, \psi_{M}^2, \rho_M)\in \:^{\sf H_2}_{\sf H_2}{\sf YD}$.
\end{enumerate} 

Observe that the first condition, together with (\ref{echin}), imply (iv) of Definition \ref{YDbrace}. And as was observed in Remark \ref{chino}, any module in the sense of \cite{Zhu} is a module in the sense of Definition \ref{l-mod}. Finally, the morphisms of the category $_\mathbb{H}^{\sf H'}{\sf YD}$ are morphisms of left $\mathbb H$-modules and of $H$-comodules, as stated in Definition \ref{YDbrace}. Thus, $_\mathbb{H}^{\sf H'}{\sf YD}$ is a full subcategory of $_\mathbb{H}^{\mathbb H}{\sf WYD}$. 

Note that, if $(M,\psi_{M}^1, \psi_{M}^2, \rho_{M})$ and $(N,\psi_{N}^1, \psi_{N}^2, \rho_{N})$ are  objects in $_\mathbb{H}^{\sf H^{\prime}}{\sf YD}$ and the antipodes of $H_{1}$ and $H_{2}$ are isomorphisms, by the condition of compatible module and (i), we have that $t_{M,N}^1=t_{M,N}^2$ where $t_{M,N}^i$ is the braiding of $\:^{\sf H_i}_{\sf H_i}{\sf YD}$ for $i=1,2$. Therefore, $t_{M,N}^2$ is a morphism of left $H_{1}$-modules, and this implies that $_\mathbb{H}^{\sf H'}{\sf YD}$ is a full subcategory of $_\mathbb{H}^{\mathbb H}{\sf YD}$.

On the other hand, if $H^{\prime}$ is a subbialgebra of $\mathbb{H}$, (i) holds trivially for $\rho_{H}=\eta_{H}\otimes H$ and, as a consequence, if we assume the cocommutativity condition for $\mathbb{H}$,   
$$(H, \psi_{H}^{1}=\varepsilon_{H}\otimes H, \psi_{H}^{2}=\mu_{H}^{2}, \rho_{H}=\eta_{H}\otimes H)$$
is an object in $_\mathbb{H}^{\sf H'}{\sf YD}$. Then, taking into account that we know that   $t_{M,H}^1=t_{M,H}^2$ holds for any $(M,\psi_{M}^1, \psi_{M}^2, \rho_{M})$ in $_\mathbb{H}^{\sf H'}{\sf YD}$, the identities 
\begin{equation}
%\label{trich}
c_{M,H}=(\varepsilon_{H}\otimes c_{M,H})\circ (\rho_{M}\otimes H)=(\mu_{H}^{2}\otimes M)\circ (H\otimes c_{M,H})\circ (\rho_{M}\otimes H)
\end{equation}
also hold and, as a consequence, composing on the right with $M\otimes \eta_{H}$, we obtain that 
$$\rho_{M}=\eta_{H}\otimes M.$$

Thus, under cocommutative conditions the coaction is trivial and then $_\mathbb{H}^{\sf H'}{\sf YD}$ and $\;_{\mathbb H}{\sf Mod}$ are isomorphic categories.

The previous categorical isomorphism implies that the definition introduced in \cite{Zhu} does not provide a consistent non-trivial theory of  Yetter-Drinfeld modules for Hopf braces. Note that the condition of cocommutativity is used systematically in \cite{Zhu} and \cite{Zhu2}.

\end{example}

\begin{theorem}
%\label{WYDTope}
Let's assume that {\sf C} is symmetric with natural isomorphism of symmetry $c$. Let ${\mathbb H}$ be  a cocommutative Hopf brace in {\sf C}. Then the category  $\:^{\mathbb H}_{\mathbb H}{\sf YD}$ is braided monoidal.
\end{theorem}
\begin{proof} To prove the theorem we only need to show that if $(M,\psi_{M}^1, \psi_{M}^2, \rho_{M})$ and $(N,\psi_{N}^1, \psi_{N}^2, \rho_{N})$ are  objects in $_\mathbb{H}^{\mathbb{H}}{\sf YD}$ and $(P,\psi_{P}^1, \psi_{P}^2, \rho_{P})$ is an object in $_\mathbb{H}^{\mathbb{H}}{\sf WYD}$ we have that $t_{M\otimes N, P}^{2}$ is a morphism of left $H_{1}$-modules. Indeed:
\begin{itemize}
	\item[ ]$\hspace{0.38cm}\psi_{P\otimes M\otimes N}^1\circ (H\otimes t_{M\otimes N, P}^{2}) $
	\item [ ]$= ((\psi_{P\otimes M}^1\circ (H\otimes t_{M, P}^{2}))\otimes \psi_{N}^{1})\circ (H\otimes M\otimes c_{H,P}\otimes N)\circ (((H\otimes c_{H,M})\circ (\delta_{H}\otimes M))\otimes t_{N, P}^{2}) $ 
	\item[ ]$\hspace{0.38cm}${\scriptsize (by the naturality of $c$ and the coassociativity of $\delta_{H}$)}
	\item [ ]$= ((t_{M, P}^{2}\circ \psi_{M\otimes P}^1)\otimes \psi_{N}^{1})\circ (H\otimes M\otimes c_{H,P}\otimes N)\circ (((H\otimes c_{H,M})\circ (\delta_{H}\otimes M))\otimes t_{N, P}^{2})$ {\scriptsize (by the }
	\item[ ]$\hspace{0.38cm}${\scriptsize condition of morphism of left $H_{1}$-modules for $t_{M, P}^{2}$)}
	\item [ ]$=(t_{M, P}^{2}\otimes N)\circ (\psi_{M}^{1}\otimes (\psi_{P\otimes N}^{1}\circ (H\otimes t_{N, P}^{2})))\circ (((H\otimes c_{H,M})\circ (\delta_{H}\otimes M))\otimes N\otimes P)  $ {\scriptsize (by the}
	\item[ ]$\hspace{0.38cm}${\scriptsize naturality of $c$ and the coassociativity of $\delta_{H}$)}
	\item [ ]$=t_{M\otimes N, P}^{2} \circ \psi_{M\otimes N\otimes P}^1$ {\scriptsize (by the condition of morphism of left $H_{1}$-modules for $t_{N, P}^{2}$).}
\end{itemize}
\end{proof}

\begin{notation}
Under the conditions of the previous theorem we know that  the braiding of $_\mathbb{H}^{\mathbb{H}}{\sf YD}$ is $t^{2}$, i.e., the braiding of $_{\sf H_{2}}^{\sf H_{2}}{\sf YD}$. Taking into account this fact, from this moment and to simplify the notation the braiding of the category $_\mathbb{H}^{\mathbb{H}}{\sf YD}$ will be denoted by $t$. 
\end{notation}

To finish this section we will prove that if $\mathbb{H}$ is a cocommutative Hopf brace in a symmetric monoidal category, $\:^{\mathbb H}_{\mathbb H}{\sf YD}$ can be seen as a type of categorical center.  To fix the notation and make the reading more self-contained we will first remember the notion of center of a monoidal category.

\begin{definition}
	%\label{center}
	Let ${\sf D}$ be a strict monoidal category with tensor product $\boxtimes$ and unit object  $I$. The center (or left center) ${\sf Z}({\sf D})$ is the category with the following objects and morphisms: An object is a pair $(M, \tau_{M,-})$, with $M\in {\sf D}$ and $\tau_{M,-}:M\boxtimes -\rightarrow -\boxtimes M$ a natural isomorphism satisfying the following condition for all $N,P\in {\sf D}$:
	\begin{equation}
	\label{cen1}
	\tau_{M,N\boxtimes P}=(N\boxtimes \tau_{M,P})\circ (\tau_{M,N}\boxtimes P).
	\end{equation}
	A morphism between $(M, \tau_{M,-})$ and $(M^{\prime}, \tau_{M^{\prime},-})$ consists on a morphism $f:M\rightarrow M^{\prime}$ in ${\sf D}$ such that 
	\begin{equation}
	\label{cen2}
	(N\boxtimes f)\circ \tau_{M,N}=\tau_{M^{\prime},N}\circ (f\boxtimes N).
	\end{equation}
	for all $N\in {\sf D}$.
	
	Note that, as a consequence of the strict character of ${\sf D}$, we have that $\tau_{M,I}=id_{M}$ for all $M\in {\sf D}$.
	
	The center ${\sf Z}({\sf D})$ is a strict braided monoidal category.
	The tensor product is 
	$$(M, \tau_{M,-})\boxtimes (M^{\prime}, \tau_{M^{\prime},-})=(M\boxtimes M^{\prime}, \tau_{M\boxtimes M^{\prime},-})$$
	with 
	$$\tau_{M\boxtimes M^{\prime},N}=
	(\tau_{M,N}\boxtimes M^{\prime})\circ (M\boxtimes \tau_{M^{\prime},N})$$
	and the unit object is $(I, \tau_{I,-}={\sf id}_{\sf D}).$
	
	The braiding is given by 
	$$\tau_{M,M^{\prime}}: (M, \tau_{M,-})\boxtimes (M^{\prime}, \tau_{M^{\prime},-})\rightarrow (M^{\prime}, \tau_{M^{\prime},-})\boxtimes (M, \tau_{M,-}).$$ 
	
\end{definition}

\begin{example}
	%\label{excen}
	Let $X$ be a Hopf algebra in ${\sf C}$. The category of left $X$-modules is a monoidal category where the tensor product of two objects $(M, \psi_{M})$, $(N, \psi_{N})$ is defined by 
	$(M\otimes N, \psi_{M\otimes N})$ with $\psi_{M\otimes N}$ the tensor module structure. The unit object is $(K, \psi_{K}=\varepsilon_{X})$. Then, ${\sf Z}(_{\sf X}{\sf Mod})$ is a strict braided monoidal category where the objects can be identified with triples $(M,\psi_{M}, \tau_{M,-})$
	where $(M,\psi_{M})$ is an object in $\;_{\sf X}{\sf Mod}$  and $\tau_{M,N}:M\otimes N\rightarrow N\otimes M$ is a family of natural isomorphisms in $\;_{\sf X}{\sf Mod}$ satisfying (\ref{cen1}). Also in this case the morphisms in ${\sf Z}(_{\sf X}{\sf Mod})$ are morphisms in $\;_{\sf X}{\sf Mod}$ satisfying (\ref{cen2}).
\end{example}

\begin{definition}
	%\label{center1} 
	Let $X$ be a Hopf algebra in ${\sf C}$. We define the small center of $\;_{\sf X}{\sf Mod}$ as the full subcategory ${\sf SZ}(_{\sf X}{\sf Mod})$ of ${\sf Z}(_{\sf X}{\sf Mod})$ with objects $(M,\psi_{M}, \tau_{M,-})$ satisfying 
	\begin{equation}
	\label{cen3}
	\tau_{M,N}=(\psi_{N}\otimes M) \circ (X\otimes c_{M,N})\circ ((\tau_{M,X}\circ (M\otimes \eta_{X}))\otimes N),
	\end{equation}
	for all $(N,\psi_{N})$ in $\;_{\sf X}{\sf Mod}$.
	
	Note that, if $(M,\psi_{M}, \tau_{M,-})$, $(M^{\prime},\psi_{M^{\prime}}, \tau_{M^{\prime},-})$ are objects in the category ${\sf SZ}(_{\sf X}{\sf Mod})$, the tensor product $(M\otimes M^{\prime},\psi_{M\otimes M^{\prime}}, \tau_{M\otimes M^{\prime},-})$ also is. Indeed, let $(N,\psi_{N})$ in $\;_{\sf X}{\sf Mod}$, then:
	\begin{itemize}
		\item[ ]$\hspace{0.38cm} (\psi_{N}\otimes M\otimes M^{\prime})\circ (X\otimes c_{M\otimes M^{\prime}, N})\circ ((\tau_{M\otimes M^{\prime},X}\circ (M\otimes M^{\prime}\otimes \eta_{X}))\otimes N)) $
		\item[ ]$=(\psi_{N}\otimes M\otimes M^{\prime})\circ(X\otimes c_{M,N}\otimes M^{\prime})\circ (\tau_{M,X}\otimes c_{M^{\prime}, N})\circ (M\otimes (\tau_{ M^{\prime},X}\circ ( M^{\prime}\otimes \eta_{X}))\otimes N)$ 
		\item[ ]$\hspace{0.38cm}$ {\scriptsize (by definition of $\tau_{M\otimes M^{\prime},-}$)}
		\item[ ]$=(\psi_{N}\otimes M\otimes M^{\prime})\circ(\mu_X\otimes c_{M,N}\otimes M^{\prime})\circ (X\otimes c_{M,X}\otimes  c_{M^{\prime},N})\circ ((\tau_{ M,X}\circ ( M\otimes \eta_{X}))\otimes X\otimes M^{\prime}\otimes N)$
		\item[ ]$\hspace{0.38cm}\circ (M\otimes (\tau_{ M^{\prime},X}\circ ( M^{\prime}\otimes \eta_{X}))\otimes N)${\scriptsize (by (\ref{cen3}) for $(X, \mu_{X})$)}
		\item[ ]$=((\psi_{N}\circ (X\otimes \psi_{N}))\otimes M\otimes M^{\prime})\circ(X\otimes X\otimes c_{M,N}\otimes M^{\prime})\circ (X\otimes c_{M,X}\otimes  c_{M^{\prime},N})$
		\item[ ]$\hspace{0.38cm}\circ  ((\tau_{ M,X}\circ ( M\otimes \eta_{X}))\otimes (\tau_{ M^{\prime},X}\circ ( M^{\prime}\otimes \eta_{X}))\otimes N)${\scriptsize (by the condition of left $X$-module for $N$)}
		\item[ ]$= (((\psi_{N}\otimes M) \circ (X\otimes c_{M,N})\circ ((\tau_{M,X}\circ (M\otimes \eta_{X}))\otimes N))\otimes M^{\prime})\circ (M\otimes ((\psi_{N}\otimes M^{\prime}) \circ (X\otimes c_{M^{\prime},N})$
		\item[ ]$\hspace{0.38cm}\circ ((\tau_{M^{\prime},X}\circ (M^{\prime}\otimes \eta_{X}))\otimes N))) ${\scriptsize (by naturality of $c$)}
		\item[ ]$=\tau_{M\otimes M^{\prime},N} ${\scriptsize (by (\ref{cen3}) for $(N, \psi_{N})$)}.
	\end{itemize}

	Therefore ${\sf SZ}(_{\sf X}{\sf Mod})$ is a braided monoidal subcategory of ${\sf Z}(_{\sf X}{\sf Mod})$. As a consequence, the inclusion functor is a braided strong monoidal functor.

\end{definition}

\begin{example}
%\label{excen2}	
	In the previous definition, if ${\sf C}$ is the category of $R$-modules over a commutative ring $R$, the equality (\ref{cen3}) always holds as was proved in \cite[Theorem XIII.5.2]{Christian}. Then, in this setting, ${\sf SZ}(_{\sf X}{\sf Mod})={\sf Z}(_{\sf X}{\sf Mod})$. 
\end{example}

\begin{theorem}
%\label{center11}	
	Let $X$ be a Hopf algebra in ${\sf C}$ such that $\lambda_{X}$ is an isomorphism. Then the category ${\sf SZ}(_{\sf X}{\sf Mod})$  is isomorphic 
	to the category of left Yetter-Drinfeld modules over $X$ as braided monoidal categories.
\end{theorem}

\begin{proof}
By (\ref{cen3}) the proof follows the one proposed in \cite[Theorem XIII.5.2]{Christian}. Then, we will restrict the proof of this Theorem to a brief description of the connecting functors. Take $(M,\psi_{M}, \tau_{M,-})$ in ${\sf SZ}(_{\sf X}{\sf Mod})$. The morphism $\rho_{M}=\tau_{M,X}\circ (M\otimes \eta_{X})$ makes $(M,\rho_{M})$ in a left $X$-comodule and, by (\ref{cen3}), we obtain that $(M,\psi_{M},\rho_{M})$ is a left Yetter-Drinfeld module over $X$. Conversely, if $(N,\psi_{N},\rho_{N})$ is a left Yetter-Drinfeld module over $X$, the natural isomorphism is defined by 
$\tau_{N,P}=t_{N,P}$ where $t_{N,P}$ is the braiding of $\:^{\sf X}_{\sf X}{\sf YD}$ and $(P,\psi_{P},\rho_{P})$ is an arbitrary object in $\:^{\sf X}_{\sf X}{\sf YD}$.
\end{proof}

Let ${\mathbb H}$ be a cocommutative Hopf brace in symmetric monoidal category ${\sf C}$. By Theorem \ref{WYD1}, we know that the category of left $\:^{\mathbb H}_{\mathbb H}{\sf WYD}$ is a  monoidal category. Then, ${\sf Z}(^{\mathbb H}_{\mathbb H}{\sf WYD})$ is a strict braided monoidal category where the objects can be identified with quadruples $(M,\psi_{M}^1,\psi_{M}^2, \rho_{M}, \tau_{M,-})$
where $(M,\psi_{M}^1,\psi_{M}^2, \rho_{M})$ is an object in $\;^{\mathbb H}_{\mathbb H}{\sf WYD}$  and $$\tau_{M,N}:M\otimes N\rightarrow N\otimes M$$ is a family of natural isomorphisms in $^{\mathbb H}_{\mathbb H}{\sf WYD}$ satisfying (\ref{cen1}). Also in this case, the morphisms in ${\sf Z}(^{\mathbb H}_{\mathbb H}{\sf WYD})$ are morphisms in $\;^{\mathbb H}_{\mathbb H}{\sf WYD}$ satisfying (\ref{cen2}). Note that, in this setting, $\tau_{M,N}$ is a morphism of left $H_{1}$-modules, left $H_{2}$-modules and left $H$-comodules.

\begin{definition}
%\label{center2} 
Let's assume that ${\sf C}$ is symmetric. Let ${\mathbb H}$ be a cocommutative Hopf brace in ${\sf C}$.  We define the small center of $\;^{\mathbb H}_{\mathbb H}{\sf WYD}$ as the full subcategory ${\sf SZ}(^{\mathbb H}_{\mathbb H}{\sf WYD})$ of ${\sf Z}(^{\mathbb H}_{\mathbb H}{\sf WYD})$ with objects $(M,\psi_{M}^1, \psi_{M}^2, \rho_{M}, \tau_{M,-})$ satisfying 
\begin{equation}
%\label{cenx}
\rho_{M}=\tau_{M,H}\circ (M\otimes \eta_{H}),
\end{equation} 
\begin{equation}
%\label{cen4}
\tau_{M,N}=(\psi_{N}^2\otimes M) \circ (H\otimes c_{M,N})\circ ((\tau_{M,H}\circ (M\otimes \eta_{H}))\otimes N)
\end{equation}
for all $(N,\psi_{N}^1, \psi_{N}^2, \rho_{N})$ in $\;^{\mathbb H}_{\mathbb H}{\sf WYD}$.
	
Note that, as in the Hopf algebra case,  if $(M,\psi_{M}^1, \psi_{M}^2,\rho_{M}, \tau_{M,-})$, $(M^{\prime},\psi_{M^{\prime}}^1, \psi_{M^{\prime}}^2, \rho_{M^{\prime}}, \tau_{M^{\prime},-})$ are objects in the category ${\sf SZ}(^{\mathbb H}_{\mathbb H}{\sf WYD})$, the tensor product $(M\otimes M^{\prime},\psi_{M\otimes M^{\prime}}^1, \psi_{M\otimes M^{\prime}}^2, \rho_{M\otimes M^{\prime}}, \tau_{M\otimes M^{\prime},-})$ also is. Then,  ${\sf SZ}(^{\mathbb H}_{\mathbb H}{\sf WYD})$ is a braided monoidal subcategory of ${\sf Z}(^{\mathbb H}_{\mathbb H}{\sf WYD})$. As a consequence, the inclusion functor is a braided strong monoidal functor.
\end{definition}

\begin{theorem}
%\label{YDbrace2} 
Let's assume that ${\sf C}$ is symmetric. Let ${\mathbb H}$ be a cocommutative Hopf brace in ${\sf C}$. Then the categories $\:^{\mathbb H}_{\mathbb H}{\sf YD}$ and ${\sf SZ}(^{\mathbb H}_{\mathbb H}{\sf WYD})$ are isomorphic as braided monoidal categories. 
\end{theorem}

\begin{proof}
The proof follows directly from the definitions of  the categories $\:^{\mathbb H}_{\mathbb H}{\sf YD}$ and ${\sf SZ}(^{\mathbb H}_{\mathbb H}{\sf WYD})$ because if   $(M,\psi_{M}^1, \psi_{M}^2,\rho_{M}, \tau_{M,-})$ is an object in ${\sf SZ}(^{\mathbb H}_{\mathbb H}{\sf WYD})$, we have that $(M,\psi_{M}^1, \psi_{M}^2,\rho_{M}) $ is an object in $\:^{\mathbb H}_{\mathbb H}{\sf YD}$. Conversely, if $(V,\psi_{V}^1, \psi_{V}^2,\rho_{V}) $ is an object in $\:^{\mathbb H}_{\mathbb H}{\sf YD}$, we obtain that 
$(V,\psi_{V}^1, \psi_{V}^2,\rho_{V}, \tau_{V,-}=t^2_{V,-}) $ is an object in ${\sf SZ}(^{\mathbb H}_{\mathbb H}{\sf WYD})$.  Then it is easy to prove that this correspondence defines a pair of inverse  functors and the braided monoidal isomorphism.
\end{proof}

\section{Projections of Hopf braces}

As emphasized in the introduction of this paper, the notions of  Yetter-Drinfeld module and  projection of Hopf algebras are strongly linked. In the following pages we will try to study this connection in the context of Hopf braces. Once the notion of Yetter-Drinfeld module for a Hopf brace  was introduced in Definition \ref{YDbrace-new}, in the following definition we present the notion of projection for Hopf braces.

\begin{definition}
	%\label{ppHBr}
	A projection of Hopf braces in {\sf C} is a 4-tuple  $({\mathbb H}, {\mathbb D}, x,y)$, where ${\mathbb H}$, ${\mathbb D} $  are Hopf braces in {\sf C}, $x:{\mathbb H}\rightarrow {\mathbb D}$, $y:{\mathbb D}\rightarrow {\mathbb H}$ are morphisms of Hopf braces in ${\sf C}$ and the following equality $y\circ x=id_{\mathbb H}$ holds.

	A morphism between two projections of Hopf braces  $({\mathbb H}, {\mathbb D}, x,y)$ and  $({\mathbb H}^{\prime}, {\mathbb D}^{\prime}, x^{\prime},y^{\prime})$ is a pair $(z ,t)$ where $z:{\mathbb H}\rightarrow {\mathbb H}^{\prime}$,  $t:{\mathbb D}\rightarrow {\mathbb D}^{\prime}$ are morphisms in {\sf HBr} and the following equalities hold:
	\begin{equation}
		\label{eqpHBr}
		x^{\prime}\circ z=t\circ x, \;\;\;\; y^{\prime}\circ t=z\circ y.
	\end{equation}
	
	With this morphisms and the previous objects we can define the category of projections of Hopf braces. We will denote this category by {\sf P(HBr)}.
	
Note that (\ref{eqpHBr}) implies that 
\begin{equation}
%\label{eqpHBrx}
z= y^{\prime}\circ t\circ x.
\end{equation}

\end{definition}

\begin{remark}
	\label{qd1qd2}
	If $({\mathbb H}, {\mathbb D}, x,y)$ is a projection of Hopf braces in {\sf C}, we have two projections of Hopf algebras $(H_{1}, D_{1}, x,y)$ and $(H_{2}, D_{2}, x,y)$. Then, with $q_{D}^1$ and $q_{D}^2$ we will denote the associated idempotent morphisms.  Note that, if $q_{D}^1=i_{D}^1\circ p_{D}^1 $ and $q_{D}^2=i_{D}^2\circ p_{D}^2$, with  $ p_{D}^1 \circ  i_{D}^1 =id_{I(q_{D}^1)}$ and $p_{D}^2 \circ  i_{D}^2 =id_{I(q_{D}^2)}$, we have that 
	$$
	\setlength{\unitlength}{1mm}
	\begin{picture}(101.00,12.00)
		\put(20,7){\vector(1,0){20}}
		\put(46,8){\vector(1,0){27}}
		\put(46,6){\vector(1,0){27}}
		\put(13,7){\makebox(0,0){$I(q_{D}^{k})$}}
		\put(43,7){\makebox(0,0){$D$}}
		\put(80,7){\makebox(0,0){$D\otimes H$}}
		\put(27.5,10){\makebox(0,0){$i_{D}^{k}$}}
		\put(60,11){\makebox(0,0){$(D\otimes y)\co \delta_{D}$}}
		\put(60,3){\makebox(0,0){$D\otimes \eta_{H}$}}
	\end{picture}
	$$
	is an equalizer diagram for $k\in \{1,2\} $ and, as a consequence, we can assume that  $i_{D}^1=i_{D}^2$  and $I(q_{D}^1)=I(q_{D}^2)$.  Then,  $p_{D}^1\circ i_{D}^1=id_{I(q_{D}^1)}=p_{D}^2\circ i_{D}^1$ and, composing with $p_{D}^1$ and $p_{D}^2$ we  obtain the  equalities 

\begin{equation}
\label{1emis21}
p_{D}^2=p_{D}^1\circ q_{D}^2, \;\;\;\; p_{D}^1=p_{D}^2\circ q_{D}^1.
\end{equation}

Therefore, 
\begin{equation}
\label{1emis221}
q_{D}^1=q_{D}^2\circ q_{D}^1, \;\;\;\; q_{D}^2=q_{D}^1\circ q_{D}^2
\end{equation}
also hold.

\end{remark}

\begin{notation}
Taking into account the previous remark, in what follows we will denote the morphism $i_{D}^1=i_{D}^2$
by $i_{D}$ and the objects $I(q_{D}^1)=I(q_{D}^2)$ by $I(q_{D})$.
\end{notation}

\begin{remark}
Note that,  by (\ref{1emis21}) and (\ref{m-i}), we have that 
$$	\mu_{I(q_{D})}^{1}=p_{D}^1\circ \mu_{D}^{1}\circ (i_{D}\otimes i_{D})=p_{D}^2\circ q_{D}^1\circ \mu_{D}^{1}\circ (i_{D}\otimes i_{D})=p_{D}^2\circ \mu_{D}^{1}\circ (i_{D}\otimes i_{D}),$$
and, similarly, 
$$	\mu_{I(q_{D})}^{2}=p_{D}^1\circ \mu_{D}^{2}\circ (i_{D}\otimes i_{D}).$$
\end{remark}

\begin{theorem}
	\label{th-proj-cc1}
	Let $({\mathbb H}, {\mathbb D}, x,y)$ be a projection of Hopf braces where 
	${\mathbb D}$ is cocommutative. Then, the following equality 
	\begin{equation}
		\label{qdqd}
		(q_{D}^1\otimes D)\circ \delta_{D}\circ i_{D}=(q_{D}^2\otimes D)\circ \delta_{D}\circ i_{D}
	\end{equation}
	holds where $q_{D}^1$ and $q_{D}^2$ are the idempotent morphisms introduced in Remark \ref{qd1qd2}.
\end{theorem}

\begin{proof} If $({\mathbb H}, {\mathbb D}, x,y)$ is a projection of Hopf braces with ${\mathbb D}$ cocommutative,  by Lemma \ref{cole}, we have that $q_{D}^1$ is a coalgebra morphism. Then, 
\begin{itemize}
		\item[ ]$\hspace{0.38cm}(q_{D}^2\otimes D )\circ \delta_{D}\circ q_{D}^1 $
		\item [ ]$=((q_{D}^2\circ q_{D}^1)\otimes q_{D}^1)\circ \delta_{D}  $ {\scriptsize (by the condition of coalgebra morphism for $q_{D}^1$)}
		\item [ ]$=(q_{D}^1\otimes q_{D}^1 )\circ \delta_{D}  $ {\scriptsize (by (\ref{1emis221}))}
		\item [ ]$=(q_{D}^1\otimes D)\circ \delta_{D}\circ q_{D}^1 $ {\scriptsize (by the condition of coalgebra morphism for $q_{D}^1$ and $q_{D}^1\circ q_{D}^1=q_{D}^1$)}
\end{itemize}
	holds, and as a consequence, composing with $i_{D}$, we obtain (\ref{qdqd}).
\end{proof}		

\begin{remark}
%\label{claro}
Note that,  if (\ref{qdqd}) holds, using (\ref{id1}) and (\ref{1emis221}), we obtain that 
	\begin{equation}
		\label{qdqdn}
		(q_{D}^1\otimes q_{D}^1)\circ \delta_{D}\circ i_{D}=(q_{D}^2\otimes q_{D}^1)\circ \delta_{D}\circ i_{D}= (q_{D}^2\otimes q_{D}^2)\circ \delta_{D}\circ i_{D}
	\end{equation}
	holds. Then, the idempotent morphisms $q_{D}^1$ and $q_{D}^2$ induce the same coproduct in $I(q_{D})$. By Theorem \ref{th-proj-cc1},  this is the situation that occurs when $({\mathbb H}, {\mathbb D}, x,y)$ is a projection of Hopf braces with ${\mathbb D}$ cocommutative.
\end{remark}

In the following theorem we will prove that, under cocommutative conditions, projections of Hopf braces induces new Hopf braces.

\begin{theorem}
	\label{311}
	Let $({\mathbb H}, {\mathbb D}, x,y)$ be a projection of Hopf braces with
	${\mathbb D}$  cocommutative. Then,
	$$
	{\mathbb I(q_{D})}=(I(q_{D}), \eta_{I(q_{D})}, \mu_{I(q_{D})}^{1}, \eta_{I(q_{D})}^{2}, \mu_{I(q_{D})}^{2}, \varepsilon_{I(q_{D})}, \delta_{I(q_{D})}, \lambda_{I(q_{D})}^{1}, \lambda_{I(q_{D})}^{2})$$
	is a Hopf brace in ${\sf C}$ where:
	\begin{equation}
		\label{3111}
		\eta_{I(q_{D})}=p_{D}^1\circ \eta_{D}=p_{D}^2\circ \eta_{D},
	\end{equation}
	\begin{equation}
		\label{3112}
		\mu_{I(q_{D})}^{1}=p_{D}^1\circ \mu_{D}^{1}\circ (i_{D}\otimes i_{D}),
	\end{equation}
	\begin{equation}
		\label{3113}
		\mu_{I(q_{D})}^{2}=p_{D}^2\circ \mu_{D}^{2}\circ (i_{D}\otimes i_{D}),
	\end{equation}
	\begin{equation}
		\label{31131}
		\varepsilon_{I(q_{D})}=\varepsilon_{D}\circ i_{D},
	\end{equation}
	\begin{equation}
		\label{3114}
		\delta_{I(q_{D})}=(p_{D}^1\otimes p_{D}^1)\circ \delta_{D}\circ 
		i_{D}=(p_{D}^2\otimes p_{D}^2)\circ \delta_{D}\circ 
		i_{D},
	\end{equation}
	\begin{equation}
%\label{3115}
		\lambda_{I(q_{D})}^{1}=p_{D}^1\circ \lambda_{D}^1\circ i_{D}, 
	\end{equation}
	\begin{equation}
%\label{3116}
		\lambda_{I(q_{D})}^{2}=p_{D}^2\circ \lambda_{D}^2\circ i_{D}.
	\end{equation}
\end{theorem}

\begin{proof} First note that $\eta_{I(q_{D})}= p_{D}^1\circ \eta_{D}=p_{D}^2\circ \eta_{D}$  holds because $\eta_{I(q_{D})}^1$ is the unique morphism such that 
$i_{D}\circ \eta_{I(q_{D})}^1=\eta_{D}$ and $\eta_{I(q_{D})}^2$ is the unique morphism such that 
$i_{D}\circ \eta_{I(q_{D})}^2=\eta_{D}$.  Thus, from now on, we will use $\eta_{I(q_{D})}$ to denote the morphism $\eta_{\ID}^{1}=\eta_{\ID}^{2}$. Also, 
by the cocommutativity of ${\mathbb D}$, we can assure that (\ref{3114}) holds. Therefore, by the general theory of Hopf algebra projections  $(I(q_{D}), \eta_{I(q_{D})}, \mu_{I(q_{D})}^{1},  \varepsilon_{I(q_{D})}, \delta_{I(q_{D})}, \lambda_{I(q_{D})}^{1})$ and $(I(q_{D}), \eta_{I(q_{D})}, \mu_{I(q_{D})}^{2},  \varepsilon_{I(q_{D})}, \delta_{I(q_{D})}, \lambda_{I(q_{D})}^{2})$ are  Hopf algebras in ${\sf C}$ because the cocommutativity condition implies that $\rho_{I(q_{D})}=\eta_{H}\ot I(q_{D})$.

	On the other hand,  the equality 
	\begin{equation}
		\label{igii}
		i_{D}\circ \Gamma_{{I(q_{D})}_{1}}= \Gamma_{{D}_{1}}\circ (i_{D}\otimes i_{D})
	\end{equation}
	holds because
	\begin{itemize}
		\item[ ]$\hspace{0.38cm} i_{D}\circ \Gamma_{{I(q_{D})}_{1}} $
		\item [ ]$=\mu_{D}^1\circ (( i_{D}\circ \lambda_{I(q_{D})})\otimes (i_{D}\circ \mu_{I(q_{D})}^2)) \circ (\delta_{I(q_{D})}\otimes I(q_{D})) $ {\scriptsize (by (\ref{m-i}))}
		\item [ ]$=\mu_{D}^1\circ ((\lambda_{D}^{1}\circ i_{D})\otimes ( \mu_{D}^2\circ (i_{D}\otimes i_{D}))) \circ (\delta_{I(q_{D})}\otimes I(q_{D}))   $ {\scriptsize (by (\ref{m-i}) and (\ref{lxy}))}
		\item [ ]$=\mu_{D}^1\circ ((\lambda_{D}^{1}\circ q_{D}^1)\otimes ( \mu_{D}^2\circ (q_{D}^1\otimes i_{D}))) \circ ((\delta_{D}\circ i_{D})\otimes I(q_{D}))  $ {\scriptsize (by (\ref{coalg-id}))}
		\item [ ]$=\Gamma_{{D}_{1}}\circ (i_{D}\otimes i_{D})$ {\scriptsize (by the condition of coalgebra morphism for $q_{D}^1$),}
	\end{itemize}
	and then we have that the identities
	$$
		%\label{igii1}
		\Gamma_{{I(q_{D})}_{1}}=p_{D}^1\circ  \Gamma_{{D}_{1}}\circ (i_{D}\otimes i_{D})=p_{D}^2\circ  \Gamma_{{D}_{1}}\circ (i_{D}\otimes i_{D})
	$$
	hold. As a consequence, 
	\begin{itemize}
		\item[ ]$\hspace{0.38cm}  \mu_{I(q_{D})}^{1}\circ ( \mu_{I(q_{D})}^{2}\otimes \Gamma_{{I(q_{D})}_{1}})\circ ({I(q_{D})}\otimes c_{{I(q_{D})},{I(q_{D})}}\otimes {I(q_{D})})\circ (\delta_{I(q_{D})}\otimes I(q_{D})\otimes I(q_{D}))$
		\item [ ]$=p_{D}^{1}\circ \mu_{D}^{1}\circ ( ( i_{D}\circ \mu_{I(q_{D})}^{2})\otimes ( i_{D}\circ \Gamma_{{I(q_{D})}_{1}}))\circ ({I(q_{D})}\otimes c_{{I(q_{D})},{I(q_{D})}}\otimes {I(q_{D})})\circ (\delta_{I(q_{D})}\otimes I(q_{D})\otimes I(q_{D}))  $ 
		\item[ ]$\hspace{0.38cm}$ {\scriptsize (by (\ref{m-i}))}
		\item [ ]$=p_{D}^{1}\circ\mu_{D}^1\circ (\mu_{D}^2 \otimes \Gamma_{D_{1}} )\circ (D\otimes c_{D,D}\otimes D)\circ (((i_{D}\otimes i_{D})\circ \delta_{I(q_{D})})\otimes i_{D}\otimes i_{D})  $ {\scriptsize (by (\ref{m-i}),  (\ref{igii}) and }
		\item[ ]$\hspace{0.38cm}${\scriptsize  naturality of $c$)}
		
		\item [ ]$=	p_{D}^{1}\circ\mu_{D}^1\circ (\mu_{D}^2 \otimes \Gamma_{D_{1}} )\circ (D\otimes c_{D,D}\otimes D)\circ (((q_{D}^1\otimes q_{D}^1)\circ \delta_{D}\circ i_{D})\otimes i_{D}\otimes i_{D}) $ {\scriptsize (by (\ref{coalg-id}))}
		\item [ ]$=p_{D}^{1}\circ\mu_{D}^1\circ (\mu_{D}^2 \otimes \Gamma_{D_{1}} )\circ (D\otimes c_{D,D}\otimes D)\circ (( \delta_{D}\circ i_{D})\otimes i_{D}\otimes i_{D})  $ {\scriptsize (by the condition of coalgebra  }
		\item[ ]$\hspace{0.38cm}$ {\scriptsize morphism for $q_{D}^1$)}
		\item [ ]$= p_{D}^{1}\circ\mu_{D}^2\circ (D\otimes \mu_{D}^1) \circ (i_{D}\otimes i_{D}\otimes i_{D})$ {\scriptsize (by (iii) of Definition \ref{H-brace})}
		\item [ ]$=  \mu_{I(q_{D})}^{2}\circ (I(q_{D})\otimes  \mu_{I(q_{D})}^{1})$ {\scriptsize (by (\ref{m-i}))}
	\end{itemize}
	hold and ${\mathbb I(q_{D})}$ is a Hopf brace.
\end{proof}

\begin{definition}
	Let $({\mathbb H}, {\mathbb D}, x,y)$ be a projection of Hopf braces in {\sf C}. We will say that it is strong if (\ref{qdqd}) and
%\begin{equation}
%\label{spro1H}
%p_{D}^1\circ \Gamma_{D_{1}}=p_{D}^1\circ \Gamma_{D_{1}}\circ (D\ot q_{D}^1) 
%\end{equation}
%and 
	\begin{equation}
		\label{spro2H-n}
		p_{D}^1\circ \mu_{D}^2\circ (x\otimes D)=p_{D}^2\circ  \mu_{D}^2\circ (x\otimes q_{D}^1)
	\end{equation}
hold. 
	
Note that (\ref{spro2H-n}) implies that 
	\begin{equation}
		\label{spro2H}
		p_{D}^1\circ \mu_{D}^2\circ (x\otimes i_{D})=p_{D}^2\circ  \mu_{D}^2\circ (x\otimes i_{D})
	\end{equation}
holds.
	
Strong projections with morphisms of projections of Hopf braces form a category that we will denote by {\sf SP(HBr)}. In other words, {\sf SP(HBr)} is the full subcategory of {\sf P(HBr)}  whose objects are strong projections.
	
Note that, by (\ref{1emis21}),  the equality (\ref{spro2H}) is equivalent to 
%\begin{equation}
%\label{spro2H1}
$$
p_{D}^1\circ \mu_{D}^2\circ (x\otimes i_{D})=p_{D}^1\circ q_{D}^2\circ  \mu_{D}^2\circ (x\otimes i_{D}).
$$
%\end{equation}
\end{definition}

\begin{theorem}
	\label{th-proj-mod1}
	Let $({\mathbb H}, {\mathbb D}, x,y)$ be a strong projection of Hopf braces. Then,  $$(I(q_{D}),\psi_{I(q_{D})}^1=p_{D}^1\circ \mu_{D}^1\circ (x\ot i_{D}), \psi_{I(q_{D})}^2=p_{D}^2\circ\mu_{D}^2\circ (x\ot i_{D}))$$ 
	is a left ${\mathbb H}$-module and the product $\mu_{I(q_{D})}^1$ defined in (\ref{3112}) is a morphism of left $H_{2}$-modules.
	
	Then, if $$\Psi_{I(q_{D})}^{H_{2}}=(\psi_{I(q_{D})}^2\otimes H)\circ (H\otimes c_{H, I(q_{D})})\circ (\delta_{H}\otimes I(q_{D})),$$ the following equality holds:
	\begin{equation}
		\label{p11}
		(\mu_{I(q_{D})}^1\otimes H)\circ (I(q_{D})\otimes \Psi_{I(q_{D})}^{H_{2}})\circ (\Psi_{I(q_{D})}^{H_{2}}\otimes I(q_{D}))=\Psi_{I(q_{D})}^{H_{2}}\circ (H\otimes \mu_{I(q_{D})}^1).
	\end{equation}
	
	Finally, if  ${\mathbb H}$ is cocommutative and 
	$$\Psi_{I(q_{D})}^{H_{1}}=(\psi_{I(q_{D})}^1\otimes H)\circ (H\otimes c_{H, I(q_{D})})\circ (\delta_{H}\otimes I(q_{D})),$$
	we have that 
	\begin{equation}
		\label{p12}
		(I(q_{D})\otimes \mu_{H}^1)\circ (\Psi_{I(q_{D})}^{H_{1}}\otimes \mu_{H}^2)\circ (H\otimes \Psi_{I(q_{D})}^{H_{2}}\otimes H)\circ (((\Gamma^{\prime}_{H_{1}}\otimes H)\circ (H\otimes c_{H,H})\circ (\delta_{H}\otimes H))\otimes I(q_{D})\otimes H)
	\end{equation}
	$$=(I(q_{D})\otimes \mu_{H}^2)\circ (\Psi_{I(q_{D})}^{H_{2}}\otimes \mu_{H}^1)\circ (H\otimes \Psi_{I(q_{D})}^{H_{1}}\otimes H)
	$$
	holds.
	
\end{theorem}

\begin{proof} If $({\mathbb H}, {\mathbb D}, x,y)$ is a strong projection of Hopf braces $$(I(q_{D}),\psi_{I(q_{D})}^1=p_{D}^1\circ \mu_{D}^1\circ (x\ot i_{D}), \psi_{I(q_{D})}^2=p_{D}^2\circ\mu_{D}^2\circ (x\ot i_{D}))$$ 
is an object in $\;_{\mathbb H}{\sf Mod}$. Indeed, first note that by the general theory of Hopf algebra projections we have that  $(I(q_{D}),\psi_{I(q_{D})}^1=p_{D}^1\circ \mu_{D}^1\circ (x\ot i_{D}))$ is an object in  $\;_{{\sf H}_{1}}{\sf Mod}$ and by similar arguments we can assure that  $(I(q_{D}),\psi_{I(q_{D})}^2=p_{D}^2\circ \mu_{D}^2\circ (x\ot i_{D}))$ is an object in  $\;_{{\sf H}_{2}}{\sf Mod}$. Finally, (\ref{mod-l1}) follows by:
	
\begin{itemize}
	\item[ ]$\hspace{0.38cm} \psi_{I(q_{D})}^1\circ (\mu_{H}^{2}\otimes \Gamma_{I(q_{D})})\circ (H\otimes c_{H,H}\otimes I(q_{D}))\circ (\delta_{H}\otimes H\otimes I(q_{D}))$
	\item [ ]$=p_{D}^{1}\circ \mu_{D}^{1} \circ (\mu_{D}^{2}\otimes \Gamma_{D_{1}})\circ (D\otimes c_{D,D}\otimes D) \circ ((\delta_{D}\circ x)\otimes x\otimes i_{D})$ {\scriptsize (by the condition  of morphism of Hopf}
	\item[ ]$\hspace{0.38cm}${\scriptsize algebras  for $x$, (\ref{id1}) and the the naturality of $c$)}
	\item [ ]$=p_{D}^{1}\circ \mu_{D}^{2} \circ (x\otimes (\mu_{D}^{1}\circ (x\otimes i_{D}) ))$ {\scriptsize (by (iii) of Definition \ref{H-brace})}
	\item [ ]$=\psi_{I(q_{D}^2)}^1\circ (H\otimes \psi_{I(q_{D})}^1)$
	{\scriptsize (by (\ref{spro2H-n})).}
\end{itemize}

The product $\mu_{I(q_{D})}^1$ is a morphism of left $H_{2}$-modules because:
	
	\begin{itemize}
		\item[ ]$\hspace{0.38cm} \mu_{I(q_{D})}^1\circ (\psi_{I(q_{D})}^2\otimes \psi_{I(q_{D})}^2)\circ (H\otimes c_{H, I(q_{D})}\otimes I(q_{D}))\circ (\delta_{H}\otimes I(q_{D})\otimes I(q_{D})) $
		\item [ ]$=p_{D}^{1}\circ \mu_{D}^{1}\circ ((q_{D}^{1}\circ \mu_{D}^{2}\circ (x\otimes i_{D}))\otimes (\mu_{D}^{2}\circ (x\otimes i_{D}^{1}))) \circ (H\otimes c_{H, I(q_{D})}\otimes I(q_{D}))\circ (\delta_{H}\otimes I(q_{D})\otimes I(q_{D}))$ 
		\item[ ]$\hspace{0.38cm}${\scriptsize (by  (\ref{spro2H}) and (\ref{id1}))}
		\item [ ]$=p_{D}^{1}\circ \mu_{D}^{1}\circ ((q_{D}^{1}\circ \mu_{D}^{2}\circ (D\otimes i_{D}))\otimes (\mu_{D}^{2}\circ (D\otimes i_{D}))) \circ (D\otimes c_{D, I(q_{D})}\otimes I(q_{D}))$ 
		\item[ ]$\hspace{0.38cm}\circ ((\delta_{D}\circ x)\otimes I(q_{D})\otimes I(q_{D}))$ {\scriptsize (by the condition of coalgebra morphism for $x$ and the naturality of $c$)}
		\item [ ]$=p_{D}^{1}\circ \mu_{D}^{1}\circ ((\mu_{D}^{1}\circ (\mu_{D}^2\otimes (x\circ \lambda_{H}^{1}\circ \mu_{H}^{2}\circ (y\otimes H))))\circ (D\otimes c_{D,D}\otimes D)\circ (\delta_{D}\otimes ((D\otimes y)\circ \delta_{D}\circ i_{D})))$ 
		\item[ ]$\hspace{0.38cm}\otimes (\mu_{D}^{2}\circ (D\otimes i_{D}))) \circ (D\otimes c_{D, I(q_{D})}\otimes I(q_{D}))\circ ((\delta_{D}\circ x)\otimes I(q_{D})\otimes I(q_{D}))$ {\scriptsize (by the condition of}
		\item[ ]$\hspace{0.38cm}${\scriptsize algebra morphism for $y$ and the condition of coalgebra morphism for $\mu_{D}^2$)}
		\item [ ]$=p_{D}^{1}\circ \mu_{D}^{1}\circ ((\mu_{D}^{1}\circ (\mu_{D}^2\otimes (x\circ \lambda_{H}^{1}\circ y)\circ  (D\otimes c_{D,D})\circ (\delta_{D}\otimes  i_{D}))) \otimes (\mu_{D}^{2}\circ (D\otimes i_{D})))$
		\item[ ]$\hspace{0.38cm}\circ (D\otimes c_{D, I(q_{D})}\otimes I(q_{D}))\circ ((\delta_{D}\circ x)\otimes I(q_{D})\otimes I(q_{D})) $ {\scriptsize (in this equality we used that $i_{D}$ is the equalizer }
		\item[ ]$\hspace{0.38cm}${\scriptsize  morphism of $(D\otimes y)\circ \delta_{D}$ and $D\otimes \eta_{H}$)}
		\item [ ]$= p_{D}^{1}\circ \mu_{D}^{1}\circ (\mu_{D}^{2}\otimes \Gamma_{D_1})\circ (D\otimes c_{D,D}\otimes D)\circ ((\delta_{D}\circ x)\otimes i_{D}\otimes i_{D})$ {\scriptsize (by the condition of Hopf algebra}
		\item[ ]$\hspace{0.38cm}${\scriptsize  morphism for $x$, $y\circ x=id_{H}$ and the naturality of $c$)}
		\item [ ]$=p_{D}^{1}\circ \mu_{D}^{2}\circ (x\otimes (\mu_{D}^{1}\circ (i_{D}\otimes i_{D}))) $ {\scriptsize (by (iii) od Definition \ref{H-brace})}
		\item [ ]$=\psi_{I(q_{D})}^2\circ (H\otimes \mu_{I(q_{D})}^1)$ {\scriptsize (by  (\ref{spro2H}) and (\ref{m-i}))}.
	\end{itemize}
	
On the other hand, (\ref{p11}) follows by

	\begin{itemize}
		\item[ ]$\hspace{0.38cm} (\mu_{I(q_{D})}^1\otimes H)\circ (I(q_{D})\otimes \Psi_{I(q_{D})}^{H_{2}})\circ (\Psi_{I(q_{D})}^{H_{2}}\otimes I(q_{D}))$
		\item [ ]$= ((\mu_{I(q_{D})}^1\circ (\psi_{I(q_{D})}^2\otimes \psi_{I(q_{D})}^2)\circ (H\otimes c_{H, I(q_{D})}\otimes I(q_{D}))\circ (\delta_{H}\otimes I(q_{D})\otimes I(q_{D})))\otimes H)$
		\item[ ]$\hspace{0.38cm}\circ (H\otimes I(q_{D})\otimes c_{H,I(q_{D})}) \circ (H\otimes c_{H,I(q_{D})}\otimes I(q_{D}))\circ (\delta_{H}\otimes I(q_{D})\otimes I(q_{D}))$ {\scriptsize (by  the naturality of }
		\item[ ]$\hspace{0.38cm}$ {\scriptsize  $c$  and the coassociativity of $\delta_{H}$)}
		
		\item [ ]$= ((\psi_{I(q_{D})}^2\circ (H\otimes \mu_{I(q_{D})}^1))\otimes H)\circ (H\otimes I(q_{D})\otimes c_{H,I(q_{D})}) \circ (H\otimes c_{H,I(q_{D})}\otimes I(q_{D}))$
		\item[ ]$\hspace{0.38cm}\circ (\delta_{H}\otimes I(q_{D})\otimes I(q_{D}))$  {\scriptsize (by the condition  of morphism of left $H_{2}$-modules for $\mu_{I(q_{D})}^1$)}
		\item [ ]$=\Psi_{I(q_{D})}^{H_{2}}\circ (H\otimes \mu_{I(q_{D})}^1)  $ {\scriptsize (by the naturality of $c$).}
	\end{itemize}
	
	Finally, the proof of (\ref{p12}) is the following:
	
	\begin{itemize}
		\item[ ]$\hspace{0.38cm} (I(q_{D})\otimes \mu_{H}^1)\circ (\Psi_{I(q_{D})}^{H_{1}}\otimes \mu_{H}^2)\circ (H\otimes \Psi_{I(q_{D})}^{H_{2}}\otimes H)$
		\item[ ]$\hspace{0.38cm}\circ (((\Gamma^{\prime}_{H_{1}}\otimes H)\circ (H\otimes c_{H,H})\circ (\delta_{H}\otimes H))\otimes I(q_{D})\otimes H)$
		\item [ ]$= (\psi_{I(q_{D})}^1\otimes \mu_{H}^1)\circ (H\otimes c_{H, I(q_{D})}\otimes H)$
		\item[ ]$\hspace{0.38cm}\circ (((\Gamma^{\prime}_{H_{1}}\otimes \Gamma^{\prime}_{H_{1}})\circ \delta_{H\otimes H})\otimes ((\psi_{I(q_{D})}^2\otimes \mu_{H}^2)\circ (H\otimes c_{H, I(q_{D})}\otimes H)\circ (\delta_{H}\otimes I(q_{D})\otimes H)))$
		\item[ ]$\hspace{0.38cm}\circ (H\otimes c_{H,H}\otimes I(q_{D})\otimes H)\circ (\delta_{H}\otimes H\otimes \otimes I(q_{D})\otimes H)$ {\scriptsize (by  Lemma \ref{dphi-p})}
		\item [ ]$= ((\psi_{I(q_{D})}^1\circ (\Gamma^{\prime}_{H_{1}}\otimes \psi_{I(q_{D})}^2 )\circ (H\otimes c_{H,H}\otimes I(q_{D})))\otimes (\mu_{H}^1\circ (\Gamma^{\prime}_{H_{1}}\otimes \mu_{H}^2)\circ (H\otimes c_{H,H}\otimes H)))$
		\item[ ]$\hspace{0.38cm}\circ (H\otimes H\otimes ((H\otimes c_{H, I(q_{D})}\otimes H)\circ (c_{H,H}\otimes c_{H, I(q_{D})})\circ (H\otimes c_{H,H}\otimes I(q_{D}))) \otimes H\otimes H)$ 
		\item[ ]$\hspace{0.38cm}\circ (((\delta_{H}\otimes \delta_{H})\circ \delta_{H})\otimes ((H\otimes c_{H, I(q_{D})})\circ (\delta_{H}\otimes I(q_{D})))\otimes H)
		$ {\scriptsize (by  the naturality of   $c$, the cocommutativity }
		\item[ ]$\hspace{0.38cm}$ {\scriptsize  of $\delta_{H}$ and $c_{H,H}\circ c_{H,H}=id_{H}$)}
		\item [ ]$= ((\psi_{I(q_{D})}^1\circ (\Gamma^{\prime}_{H_{1}}\otimes \psi_{I(q_{D})}^2 )\circ (H\otimes c_{H,H}\otimes I(q_{D}))\circ (\delta_{H}\otimes H\otimes I(q_{D})))\otimes (\mu_{H}^1\circ (\Gamma^{\prime}_{H_{1}}\otimes \mu_{H}^2)$
		\item[ ]$\hspace{0.38cm}\circ (H\otimes c_{H,H}\otimes H)\circ (\delta_{H}\otimes H\otimes H)))\circ (H\otimes  ((H\otimes c_{H, I(q_{D})}\otimes H)\circ (c_{H,H}\otimes c_{H, I(q_{D})}))\otimes H)$ 
		\item[ ]$\hspace{0.38cm}\circ (\delta_{H}\otimes \delta_{H}\otimes I(q_{D})\otimes H)$
		{\scriptsize (by the naturality of $c$)}
		\item [ ]$= ((\psi_{I(q_{D})}^2\circ (H\otimes \psi_{I(q_{D})}^1))\otimes (\mu_{H}^2\circ (H\otimes \mu_{H}^1)))\circ (H\otimes  ((H\otimes c_{H, I(q_{D})}\otimes H)\circ (c_{H,H}\otimes c_{H, I(q_{D})}))\otimes H)$ 
		\item[ ]$\hspace{0.38cm}\circ (\delta_{H}\otimes \delta_{H}\otimes I(q_{D})\otimes H)$ {\scriptsize (by  (\ref{mod-l1p}) for $I(q_{D})$ and (\ref{iii-eq}))}
		\item [ ]$=(I(q_{D})\otimes \mu_{H}^2)\circ (\Psi_{I(q_{D})}^{H_{2}}\otimes \mu_{H}^1)\circ (H\otimes \Psi_{I(q_{D})}^{H_{1}}\otimes H) $ {\scriptsize (by the naturality of $c$).}
	\end{itemize}
\end{proof}
 
\begin{remark}
\label{coacI}	
Let $({\mathbb H}, {\mathbb D}, x,y)$ be a projection of Hopf braces. The idempotent morphisms $q_{D}^{1}$ and $q_{D}^{2}$ induce the same coaction on $I(q_{D})$ because, by (\ref{id1}) and (\ref{1emis21}), we have that 
$$\rho_{I(q_{D})}^1=(y\otimes p_{D}^{1})\circ \delta_{D}\circ i_{D}=(y\otimes (p_{D}^{1}\circ q_{D}^2))\circ \delta_{D}\circ i_{D}=(y\otimes p_{D}^{2})\circ \delta_{D}\circ i_{D}=\rho_{I(q_{D})}^2.$$

Then, in the following, we will denote this coaction by $\rho_{I(q_{D})}$.
\end{remark}

\begin{definition}
Let $({\mathbb H}, {\mathbb D}, x,y)$ be a strong projection of Hopf braces in {\sf C}. We will say that it is ${\rm v}_{1}$-strong if 
\begin{equation}
\label{VS2} (\mu_{H}^1\otimes I(q_{D}))\circ (H\otimes c_{I(q_{D}),H})\circ (\rho_{I(q_{D})}\otimes H)=
(\mu_{H}^2\otimes I(q_{D}))\circ (H\otimes c_{I(q_{D}),H})\circ (\rho_{I(q_{D})}\otimes H), 
\end{equation}
holds and the morphism 
\begin{equation}
\label{morphd1}
(\psi_{N}^2\otimes I(q_{D}))\circ (H\otimes c_{I(q_{D}),N})\circ (\rho_{I(q_{D})}\otimes N) 
:I(q_{D})\otimes N\rightarrow N\otimes I(q_{D}) 
\end{equation}
is a morphism of left $H_{1}$-modules for all $(N, \psi_{N}^1, \psi_{N}^2, \rho_{N}) \in \;_{\mathbb H}^{\mathbb H}{\sf WYD}$.
	
These  projections with morphisms of projections of Hopf braces form a category that we will denote by ${\sf V}_{1}${\sf SP(HBr)}, i.e., ${\sf V}_{1}${\sf SP(HBr)} is the full subcategory of {\sf SP(HBr)}  whose objects are ${\rm v}_{1}$-strong projections.

\end{definition}

\begin{theorem}
\label{YDbrace3}  Let $({\mathbb H}, {\mathbb D}, x,y)$ be a ${\rm v}_{1}$-strong projection of Hopf braces. Then, the triple 
$$ (I(q_{D}),\psi_{I(q_{D})}^1=p_{D}^1\circ \mu_{D}^1\circ (x\ot i_{D}), \psi_{I(q_{D})}^2=p_{D}^2\circ\mu_{D}^2\circ (x\ot i_{D}), \rho_{I(q_{D})}=(y\otimes p_{D}^{1})\circ \delta_{D}\circ i_{D})$$
is an object in $\:^{\mathbb H}_{\mathbb H}{\sf YD}.$
\end{theorem}

\begin{proof}
By Theorem \ref{th-proj-mod1}, we know that  $$(I(q_{D}),\psi_{I(q_{D})}^1=p_{D}^1\circ \mu_{D}^1\circ (x\ot i_{D}), \psi_{I(q_{D})}^2=p_{D}^2\circ\mu_{D}^2\circ (x\ot i_{D}))$$ is a left ${\mathbb H}$-module. Also, by Remark \ref{coacI}, the coaction $\rho_{I(q_{D})}$ does not depend on $q_{D}^{1}$ and $q_{D}^{2}$. On the other hand, by the general theory of Hopf algebra projections,  $(I(q_{D}),\psi_{I(q_{D})}^1, \rho_{I(q_{D})}) $ is a left Yetter-Drinfeld module over $H_{1}$ and, similarly, $(I(q_{D}),\psi_{I(q_{D})}^2, \rho_{I(q_{D})}) $ is a left Yetter-Drinfeld module over $H_{2}$. Finally,  (iv) of Definition \ref{YDbrace} and the $H_{1}$-linearity of the morphism  defined in (\ref{morphd1}) follows directly from the condition of  ${\rm v}_{1}$-strong projection.
\end{proof}

\begin{example}
%\label{MP}
In \cite{AGORE} we can find constructions for Hopf braces by means of using matched pairs of Hopf algebras in a category of vector spaces over a field ${\mathbb F}$ or, in a more general setting, in a symmetric monoidal category ${\sf C}$ that we will assume strict without loss of generality. Recall that  a matched pair  of Hopf algebras in  ${\sf C}$ is a system $(A,H, \varphi_{A}, \psi_{H})$, where $A$ and $H$ are Hopf algebras, $A$ is a left $H$-module coalgebra with action $\varphi_{A}:H\otimes A\rightarrow A$, $H$ is a right $A$-module coalgebra with action $\psi_{H}:H\otimes A\rightarrow H$ and the following  conditions hold: 
$$
%\label{MP1}
\varphi_{A}\circ (H\otimes \eta_{A})=\varepsilon_{H}\otimes \eta_{A},
$$
$$
%\label{MP2}
\psi_{H}\circ ( \eta_{H}\otimes A)=\eta_{H}\otimes \varepsilon_{A},
$$
$$
%\label{MP3}
\varphi_{A}\circ (H\otimes \mu_{A})=\mu_{A}\circ (A\otimes \varphi_{A})\circ (\Psi_{A}^{H}\otimes A),
$$
$$
%\label{MP4}
\psi_{H}\circ (\mu_H\otimes A)=\mu_{H}\circ (\psi_{H}\otimes H)\circ (H\otimes \Psi_{A}^{H}),
$$
$$
%\label{MP5}
(\psi_{H}\otimes \varphi_A)\circ \delta_{H\otimes A}=c_{A,H}\circ \Psi_{A}^{H},
$$
	
where $\Psi_{A}^{H}=(\varphi_A\otimes \psi_{H})\circ \delta_{H\otimes A}$.
	
If $(A,H, \varphi_{A}, \psi_{H})$ is a matched pair  of Hopf algebras, the double cross product $A\bowtie H$ of $A$ with $H$ is the Hopf algebra built on the object $A\otimes H$ with product 
$$\mu_{A\bowtie H}=(\mu_{A}\otimes \mu_{H})\circ (A\otimes \Psi_{A}^{H}\otimes H)$$
and tensor product unit, counit, coproduct and antipode 
$$\lambda_{A\bowtie H}=\Psi_{A}^{H}\circ (\lambda_{H}\otimes \lambda_{A})\circ c_{A,H}$$
where $\lambda_{H}$ is the antipode of $H$ and $\lambda_{A}$ is the antipode of $A$.
	
Let $A$ be a Hopf algebra and ${\mathbb H}$ be a cocommutative Hopf brace. If  $(A,H_1, \varphi_{A}, \psi_{H_1})$ is a matched pair of Hopf algebras, $(A, \varphi^{2}_{A})$ is a left $H_2$-module algebra-coalgebra, $$\Gamma_{A}^{H_{2}}=(\varphi^{2}_{A}\otimes H)\circ (H\otimes c_{H,A})\circ (\delta_{H}\otimes A),$$ and the equalities 
%\begin{equation}
%\label{MP6}
$$
\varphi^{2}_{A}\circ (H\otimes \varphi_{A})
$$
%\end{equation}
$$=\varphi_{A}\circ ((\mu_{H}^1\circ (H\otimes \lambda_{H}^1))\otimes \varphi^{2}_{A})\circ (H\otimes \delta_{H}\otimes A)\circ (((\mu_{H}^2\otimes H)\circ (H\otimes c_{H,H})\circ (\delta_{H}\otimes H))\otimes A),$$
%\begin{equation}
%\label{MP7}
$$
\mu^{2}_{H}\circ (H\otimes \psi_{H_1})
$$
%\end{equation}
$$=\mu_{H}^1\circ (\psi_{H_1}\otimes H)\circ ((\mu_{H}^1\circ (H\otimes \lambda_{H}^1))\otimes \Gamma_{A}^{H_{2}})\circ (H\otimes \delta_{H}\otimes A)\circ (((\mu_{H}^2\otimes H)\circ (H\otimes c_{H,H})\circ (\delta_{H}\otimes H))\otimes A),$$
hold, by \cite[Theorem 2.5]{AGORE}, we have that the tensor product $A\otimes H$ with the products 
$$\mu_{A\bowtie H}^1=(\mu_{A}\otimes \mu_{H}^1)\circ (A\otimes \Psi_{A}^{H_1}\otimes H),$$
$$\mu_{A\bowtie H}^2=\mu_{A\sharp H_{2}}=(\mu_{A}\otimes \mu_{H}^2)\circ (A\otimes \Gamma_{A}^{H_{2}}\otimes H),$$
tensor product unit, counit, coproduct and antipodes 
$$\lambda_{A\bowtie H}^1=\Psi_{A}^{H_1}\circ (\lambda_{H}^1\otimes \lambda_{A})\circ c_{A,H},\;\;
\lambda_{A\bowtie H}^2= \Gamma_{A}^{H_{2}}\circ  (\lambda_{H}^2\otimes \lambda_{A})\circ c_{A,H},$$
is a Hopf brace that we will denote by $ A\bowtie {\mathbb H}$.

If in the previous construction we  consider the particular case where $\psi_{H_1}=H\otimes \varepsilon_{A}$ we obtain that 
$({\mathbb H}, A\bowtie {\mathbb H}, x=\eta_{A}\otimes H, y=\varepsilon_{A}\otimes H )$ is a projection of Hopf braces where 
$$q_{A\bowtie H}=q_{A\bowtie H}^{1}=q_{A\bowtie H}^{2}=A\otimes (\eta_{H}\circ \varepsilon_{H}).$$
	
Therefore, $I(q_{A\bowtie H}^{1})=I(q_{A\bowtie H}^{2})=A$,
$$p_{A\bowtie H}^{1}=p_{A\bowtie H}^{2}=A\otimes \varepsilon_{H},$$
and 
$$i_{A\bowtie H}^{1}=i_{A\bowtie H}^{2}=A\otimes \eta_{H}.$$
	
As a consequence of the previous facts, it is easy to show that (\ref{qdqd}),  (\ref{spro2H-n}) and (\ref{VS2}) hold because in this setting 
$$\rho_{I(q_{A\bowtie H})}=(y\otimes p_{A\bowtie H}^{1})\circ \delta_{A\bowtie H}\circ i_{A\bowtie H}^{1}=\eta_{H}\otimes A.$$

As a consequence, we have that 
$$c_{A,N}=(\psi_{N}^2\otimes I(q_{A\bowtie H})\circ (H\otimes c_{I(q_{A\bowtie H}),N})\circ (\rho_{I(q_{A\bowtie H})}\otimes N)$$
and, using the cocommutativity condition, we obtain that it is a morphism of left $H_{1}$-modules for all $(N, \psi_{N}^1, \psi_{N}^2, \rho_{N}) \in \;_{\mathbb H}^{\mathbb H}{\sf WYD}$.
Then  $({\mathbb H}, A\bowtie {\mathbb H}, x=\eta_{A}\otimes H, y=\varepsilon_{A}\otimes H )$ is a ${\rm v}_{1}$-strong projection of Hopf braces such that the object  ${\mathbb I(q_{A\bowtie H})}$ is the Hopf brace ${\mathbb A}_{triv}$ introduced in Theorem \ref{2-th1} because 
$\eta_{I(q_{A\bowtie H})}=\eta_{A}$, $\mu_{I(q_{A\bowtie H})}^1=\mu_{I(q_{A\bowtie H})}^2=\mu_{A}$, $\varepsilon_{I(q_{A\bowtie H})}=\varepsilon_{A}$, $\delta_{I(q_{A\bowtie H})}=\delta_{A}$ and $\lambda_{I(q_{A\bowtie H})}^1=\lambda_{I(q_{A\bowtie H})}^2=\lambda_{A}$.  

On the other hand, by Theorem \ref{YDbrace3}, we know that $ I(q_{A\bowtie H})$ with the two actions 
$$\psi_{I(q_{A\bowtie H})}^1=p_{A\bowtie H}^1\circ \mu_{A\bowtie H}^1\circ (x\ot i_{A\bowtie H}^1)=\varphi_{A},\;\;\psi_{I(q_{A\bowtie H})}^2=p_{A\bowtie H}^2\circ\mu_{A\bowtie H}^2\circ (x\ot i_{A\bowtie H}^1)=\varphi_{A}^2,$$
and trivial coaction $\rho_{I(q_{A\bowtie H})}=\eta_{H}\otimes A$
is an object in $\:^{\mathbb H}_{\mathbb H}{\sf YD}.$ Moreover, by the general theory of Hopf algebra projections, $(A, \eta_{A},  \mu_{A}, \varepsilon_{A}, \delta_{A}, \lambda_{A})$ is a Hopf algebra in $\:^{ H_1}_{H_1}{\sf YD}$ and in $\:^{H_2}_{H_2}{\sf YD}$. Therefore, $\eta_{A}$,  $\mu_{A}$, $\varepsilon_{A}$, $\delta_{A}$, $\lambda_{A}$ are morphisms of left $H_{1}$-modules, left $H_{2}$-modules and left $H$-comodules. As a consequence of these facts we obtain that ${\mathbb A}_{triv}$ is a Hopf brace in $\:^{\mathbb H}_{\mathbb H}{\sf YD}$ because in this case the braiding $t_{A,A}$ in $\:^{\mathbb H}_{\mathbb H}{\sf YD}$ is the symmetry isomorphism $c_{A,A}$. Finally, note that the previous assertions imply that $(A, \varphi_{A})$ is not only a left $H_{1}$-module coalgebra but also a left $H_{1}$-module algebra. 
\end{example}

\begin{remark}
%\label{YDbrace4} 
Let's assume that ${\sf C}$ is symmetric. Let $({\mathbb H}, {\mathbb D}, x,y)$ be a strong projection of Hopf braces with
${\mathbb D}$ cocommutative. Then, the Hopf brace ${\mathbb I(q_{D})}$, introduced in Theorem \ref{311}, with the actions of the previous theorem 
is an object in $\:^{\mathbb H}_{\mathbb H}{\sf YD}$ where $\rho_{I(q_{D})}=\eta_{H}\ot I(q_{D})$. Note that in this case, for all $(N, \psi_{N}^1, \psi_{N}^2, \rho_{N}) \in \;_{\mathbb H}^{\mathbb H}{\sf WYD}$,  $t_{I(q_{D}), N}= c_{I(q_{D}), N}$ is a morphism of left $H_{1}$-modules because if ${\mathbb D}$ is cocommutative, the Hopf brace ${\mathbb H}$ is cocommutative. 
\end{remark}

In the following theorem we present the conditions that permit to obtain, using the bossonization process, Hopf braces in {\sf C} working with Hopf braces in the category of Yetter-Drinfeld modules associated to a cocommutative Hopf brace ${\mathbb H}$ in ${\sf C}$.

\begin{theorem}
\label{mainequiv}
Let's assume that ${\sf C}$ is symmetric. Let ${\mathbb H}$ be a cocommutative Hopf brace in ${\sf C}$ and let  ${\mathbb A}$ be a  Hopf brace in $\:^{\mathbb H}_{\mathbb H}{\sf YD}$.  Let $\Psi_{A}^{H_{1}}:H\otimes A\rightarrow A\otimes H$, $\Psi_{A}^{H_{2}}:H\otimes A\rightarrow A\otimes H$ and 
$\Omega_{H}^{A}:A\otimes H\rightarrow H\otimes A$ be the morphisms defined by
$$\Psi_{A}^{H_{1}}=(\psi_{A}^{1}\otimes H)\circ (H\otimes c_{H,A})\circ (\delta_{H}\otimes A), \;\;\;\Psi_{A}^{H_{2}}=(\psi_{A}^{2}\otimes H)\circ (H\otimes c_{H,A})\circ (\delta_{H}\otimes A), $$
$$\Omega_{H}^{A}=(\mu_{H}^{1}\otimes A)\circ (H\otimes c_{A,H})\circ (\rho_{A}\otimes H).$$

Then, ${\mathbb A}\blacktriangleright\hspace{-0.1cm}\blacktriangleleft {\mathbb H}=(( A\blacktriangleright\hspace{-0.1cm}\blacktriangleleft H)_{1}, ( A\blacktriangleright\hspace{-0.1cm}\blacktriangleleft H)_{2})$, where 
$$\eta_{A\blacktriangleright\hspace{-0.1cm}\blacktriangleleft H}=\eta_{A}\otimes \eta_{H}, \;\;\; \varepsilon_{A\blacktriangleright\hspace{-0.1cm}\blacktriangleleft H}=\varepsilon_{A}\otimes \varepsilon_{H}, $$
$$\mu_{A\blacktriangleright\hspace{-0.1cm}\blacktriangleleft H}^1=(\mu_{A}^{1}\otimes \mu_{H}^1)\circ (A\otimes \Psi_{A}^{H_{1}}\otimes H), \;\;\; \mu_{A\blacktriangleright\hspace{-0.1cm}\blacktriangleleft H}^2=(\mu_{A}^{2}\otimes \mu_{H}^2)\circ (A\otimes \Psi_{A}^{H_{2}}\otimes H),$$
$$\delta_{A\blacktriangleright\hspace{-0.1cm}\blacktriangleleft H}=(A\otimes \Omega_{H}^{A}\otimes H)\circ (\delta_{A}\otimes \delta_{H}), $$
and 
$$\lambda_{A\blacktriangleright\hspace{-0.1cm}\blacktriangleleft H}^1=\Psi_{A}^{H_{1}}\circ (\lambda_{H}^1\otimes \lambda_{A}^1)\circ \Omega_{H}^{A}, \;\;\; \lambda_{A\blacktriangleright\hspace{-0.1cm}\blacktriangleleft H}^2=\Psi_{A}^{H_{2}}\circ (\lambda_{H}^2\otimes \lambda_{A}^2)\circ \Omega_{H}^{A},$$
is a Hopf brace in ${\sf C}$ if, and only if, the following equalities hold:
\begin{equation}
\label{4161} 
(A\otimes \mu_{H}^{1})\circ (\Psi_{A}^{H_{1}} \otimes \mu_{H}^{2})\circ (H\otimes \Psi_{A}^{H_{2}}\otimes H)\circ (((\Gamma_{H_{1}}^{\prime}\otimes H)\circ (H\otimes c_{H,H})\circ (\delta_{H}\otimes H))\otimes A\otimes H)
\end{equation}
$$=(A\otimes \mu_{H}^{2})\circ (\Psi_{A}^{H_{2}} \otimes \mu_{H}^{1})\circ (H\otimes \Psi_{A}^{H_{1}}\otimes H),$$
\begin{equation}
\label{4122} 
\Psi_{A}^{H_{2}} \circ (H\otimes \mu_{A}^{1})=(\mu_{A}^{1}\otimes H)\circ (A\otimes \Psi_{A}^{H_{2}})\circ (\Psi_{A}^{H_{2}}\otimes A), 
\end{equation}
\begin{equation}
\label{4123} 
(\mu_{A}^{1}\otimes H)\circ (A\otimes \Psi_{A}^{H_{1}})\circ (A\otimes ((\Gamma_{H_{1}}^{\prime}\otimes \Gamma_{A_{1}})\circ (H\otimes c_{A,H}\otimes A)\circ (\rho_{A}\otimes H\otimes A)))\circ (\delta_{A}\otimes H\otimes A)
\end{equation}
$$=(\mu_{A}^{2}\otimes H)\circ (A\otimes \Psi_{A}^{H_{1}}).$$
\end{theorem} 

\begin{proof}
First of all prove some equalities that we will need in the proof. More concretely, we will prove that the following equalities hold:
\begin{equation}
\label{1f}
\Omega_{H}^{A}\circ (\eta_{A}\otimes H)=H\otimes \eta_{A} 
\end{equation}
\begin{equation}
\label{2f}
\Omega_{H}^{A}\circ (A\otimes \eta_{H})=\rho_{A} 
\end{equation}
\begin{equation}
\label{1f1}
\delta_{A\blacktriangleright\hspace{-0.1cm}\blacktriangleleft H}\circ (\eta_{A}\otimes H)=(A\otimes H\otimes \eta_{A}\otimes H)\circ (\eta_{A}\otimes \delta_{H}), 
\end{equation}
\begin{equation}
\label{2f1}
\delta_{A\blacktriangleright\hspace{-0.1cm}\blacktriangleleft H}\circ (A\otimes \eta_{H})=((A\otimes \rho_{A})\circ \delta_{A})\otimes \eta_{H}, 
\end{equation}
\begin{equation}
\label{3f}
(\Omega_{H}^{A}\otimes A)\circ (A\otimes \Omega_{H}^{A})\circ (\delta_{A}\otimes H)=(H\otimes \delta_{A})\circ \Omega_{H}^{A}, 
\end{equation}
\begin{equation}
\label{4f}
(H\otimes \Omega_{H}^{A})\otimes (\Omega_{H}^{A}\otimes H)\circ (A\otimes \delta_{H})=( \delta_{H}\otimes A)\circ \Omega_{H}^{A}, 
\end{equation}
\begin{equation}
\label{5f}
\Psi_{A}^{H_{i}}\circ (\eta_{H}\otimes A)=A\otimes \eta_{H}, \;\; i=1,\; 2, 
\end{equation}
\begin{equation}
\label{6f}
\Psi_{A}^{H_{i}}\circ (H\otimes \eta_{A})=\eta_{A}\otimes H, \;\; i=1,\; 2, 
\end{equation}
\begin{equation}
\label{5f1}
\mu_{A\blacktriangleright\hspace{-0.1cm}\blacktriangleleft H}^{i}\circ (A\otimes \eta_{H}\otimes A\otimes H)=\mu_{A}^{i}\otimes H, \;\; i=1,\; 2, 
\end{equation}
\begin{equation}
\label{6f1}
\mu_{A\blacktriangleright\hspace{-0.1cm}\blacktriangleleft H}^{i}\circ (A\otimes H\otimes \eta_{A}\otimes H)=A\otimes \mu_{H}^{i}, \;\; i=1,\; 2, 
\end{equation}
\begin{equation}
\label{7f}
(\mu_{A}^{i}\otimes H)\circ (A\otimes \Psi_{A}^{H_{i}})\circ (\Psi_{A}^{H_{i}}\otimes A)=\Psi_{A}^{H_{i}}\circ (H\otimes \mu_{A}^{i}),
, \;\; i=1,\; 2, 
\end{equation}
\begin{equation}
\label{8f}
(A\otimes \mu_{H}^{i})\circ (\Psi_{A}^{H_{i}}\otimes H)\circ (H\otimes \Psi_{A}^{H_{i}})=\Psi_{A}^{H_{i}}\circ (\mu_{H}^{i}\otimes A),
 \;\; i=1,\; 2, 
\end{equation}
\begin{equation}
	\label{9f}
	(\Psi_{A}^{H_{i}}\otimes H)\circ (H\otimes c_{H,A})\circ (\delta_{H}\otimes A)=(A\otimes \delta_{H})\circ \Psi_{A}^{H_{i}}
	\;\; i=1,\; 2, 
\end{equation}
\begin{equation}
\label{10f}
(\Psi_{A}^{H_{2}}\otimes A)\circ  (H\otimes c_{A,A})\circ (\Omega_{H}^{A}\otimes A)=(A\otimes \Omega_{H}^{A})\circ ( t_{A,A}\otimes H)\circ (A\otimes \Psi_{A}^{H_{2}}), 
\end{equation}
\begin{equation}
\label{11f}
\Gamma_{(A\blacktriangleright\hspace{-0.1cm}\blacktriangleleft H)_1}=(A\otimes \mu_{H}^{1})\circ ((\Psi_{A}^{H_{1}}\circ (\lambda_{H}^1\otimes \Gamma_{A_1})\circ (\Omega_{H}^{A}\otimes A))\otimes \mu_{H}^{2})\circ (A\otimes ((H\otimes  \Psi_{A}^{H_{2}})\circ (\delta_{H}\otimes A))\otimes H), 
\end{equation}
\begin{equation}
\label{11f1}
\Gamma_{(A\blacktriangleright\hspace{-0.1cm}\blacktriangleleft H)_1}\circ (\eta_{A}\otimes H \otimes A\otimes H)=(A\otimes \mu_{H}^{1})\circ (\Psi_{A}^{H_{1}}\otimes \mu_{H}^{2})\circ (\lambda_{H}^1\otimes  \Psi_{A}^{H_{2}}\otimes H)\circ (\delta_{H}\otimes A\otimes H).
\end{equation}

The  equality (\ref{1f}) follows from (\ref{comod-alg}), the naturality of the braiding and the unit properties. The proof of (\ref{2f}) follows from the naturality of the braiding and the unit properties. The equality (\ref{1f1}) is a consequence of (\ref{1f}) and (\ref{2f1}) follows from (\ref{2f}). On the other hand,  we have that 
\begin{itemize}
	\item[ ]$\hspace{0.38cm} (\Omega_{H}^{A}\otimes A)\circ (A\otimes \Omega_{H}^{A})\circ (\delta_{A}\otimes H)$
	\item [ ]$=(\mu_{H}^1\otimes A\otimes A)\circ (H\otimes c_{A,H}\otimes A) \circ (H\otimes A\otimes c_{A,H})\circ (((\mu_{H}^1\otimes A\otimes A)\circ (H\otimes c_{A,H}\otimes A)$
	\item[ ]$\hspace{0.38cm} \circ (\rho_{A}\otimes \rho_{A})\circ \delta_{A})\otimes H)$ {\scriptsize (by the naturality of $c$  and the associativity of $\mu_{H}^1$)}
	\item [ ]$= (\mu_{H}^1\otimes A\otimes A)\circ (H\otimes c_{A,H}\otimes A) \circ (H\otimes A\otimes c_{A,H})\circ (((H\otimes \delta_{A})\circ \rho_{A})\otimes H)$ {\scriptsize (by (\ref{comod-coalg}))}
	\item [ ]$=(H\otimes \delta_{A})\circ \Omega_{H}^{A} $ {\scriptsize (by the naturality  of $c$)}
\end{itemize}
and then, (\ref{3f}) holds. Also, 
\begin{itemize}
	\item[ ]$\hspace{0.38cm}(H\otimes \Omega_{H}^{A})\otimes (\Omega_{H}^{A}\otimes H)\circ (A\otimes \delta_{H})$
	\item [ ]$=(\mu_{H}^1\otimes \mu_{H}^1\otimes A) \circ (H\otimes c_{H,H}\otimes c_{A,H})\circ (\delta_{H}\otimes c_{A,H}\otimes H)\circ (\rho_{A}\otimes \delta_{H})$ {\scriptsize (by the naturality of $c$}
	\item[ ]$\hspace{0.38cm}$ {\scriptsize  and the comodule condition for $A$)}
	\item [ ]$=(((\mu_{H}^1\otimes \mu_{H}^1)\circ \delta_{H\otimes H})\otimes A) \circ  (H\otimes c_{A,H})\circ (\rho_{A}\otimes H) $ {\scriptsize (by the naturality of $c$)}
	\item [ ]$=( \delta_{H}\otimes A)\circ \Omega_{H}^{A} $ {\scriptsize (by the condition of coalgebra morphism for $\mu_{H}^1$)}
\end{itemize}
hold, and we obtain (\ref{4f}). The proof of (\ref{5f}) follows by the condition of coalgebra morphism for $\eta_{H}$, the naturality of $c$ and the condition of left module for $A$. The equality (\ref{6f}) is a consequence of  the naturality of $c$ , the condition of left module algebra for $A$ and the counit properties.  The proof of (\ref{5f1}) is a consequence of (\ref{5f}) and (\ref{6f1}) follows from (\ref{6f}). The identity (\ref{7f}) holds because 
\begin{itemize}
	\item[ ]$\hspace{0.38cm}(\mu_{A}^{i}\otimes H)\circ (A\otimes \Psi_{A}^{H_{i}})\circ (\Psi_{A}^{H_{i}}\otimes A) $
	\item [ ]$=((\mu_{A}^i\circ (\psi_{A}^i\otimes \psi_{A}^i)\circ (H\otimes  c_{H,A}\otimes A)\circ (\delta_{H}\otimes A\otimes A))\otimes H)\circ (H\otimes A\otimes c_{H,A})$
	\item[ ]$\hspace{0.38cm}\circ (H\otimes c_{H,A}\otimes A)\circ (\delta_{H}\otimes A\otimes A) $ {\scriptsize (by the naturality of $c$ and the coassociativity of $\delta_{H}$)}
	\item [ ]$=((\psi_{A}^i\circ (H\otimes \mu_{A}^i))\otimes H)\circ (H\otimes A\otimes c_{H,A})\circ (H\otimes c_{H,A}\otimes A)\circ (\delta_{H}\otimes A\otimes A) $ {\scriptsize (by (\ref{mod-alg}))}
	\item [ ]$=\Psi_{A}^{H_{i}}\circ (H\otimes \mu_{A}^{i}) $ {\scriptsize (by the naturality of $c$)}
\end{itemize}
and (\ref{8f}) follows by 
\begin{itemize}
	\item[ ]$\hspace{0.38cm}(A\otimes \mu_{H}^{i})\circ (\Psi_{A}^{H_{i}}\otimes H)\circ (H\otimes \Psi_{A}^{H_{i}}) $
	\item [ ]$=((\psi_{A}^i\circ (\mu_{H}^i\otimes A))\otimes \mu_{H}^i)\circ (H\otimes H\otimes c_{H,A}\otimes H)\circ (H\otimes c_{H,H}\otimes c_{H,A})\circ (\delta_{H}\otimes \delta_{H}\otimes A) $ {\scriptsize (by the}
	\item[ ]$\hspace{0.38cm}$ {\scriptsize  naturality of $c$ and the condition of left module for $A$)}
	\item [ ]$=(\psi_{A}^i\otimes H)\circ (H\otimes c_{H,A})\circ  (((\mu_{H}^i\otimes \mu_{H}^i)\circ \delta_{H\otimes H})\otimes A) $ {\scriptsize (by the naturality of $c$)}
	\item [ ]$=\Psi_{A}^{H_{i}}\circ (\mu_{H}^{i}\otimes A) $ {\scriptsize (by the condition of coalgebra morphism for $\mu_{H}^i$)}.
\end{itemize}

By the coassociativity of $\delta_{H}$ and the naturality of $c$ we obtain (\ref{9f}). The proof for the equality (\ref{10f}) is the following:
\begin{itemize}
	\item[ ]$\hspace{0.38cm} (\Psi_{A}^{H_{2}}\otimes A)\circ  (H\otimes c_{A,A})\circ (\Omega_{H}^{A}\otimes A) $
	\item [ ]$=(((\psi_{A}^2\otimes H)\circ (H\otimes c_{H,A}))\otimes A)\circ (((\mu_{H}^2\otimes \mu_{H}^i)\circ \delta_{H\otimes H})\otimes c_{A,A})\circ (H\otimes c_{A,H}\otimes A)\circ (\rho_{A}\otimes H\otimes A) $ 
	\item[ ]$\hspace{0.38cm}$ {\scriptsize (by the condition of algebra morphism for $\delta_{H}$)}
	\item [ ]$= ((\psi_{A}^2\circ (\mu_{H}^2\otimes A))\otimes \mu_{H}^2\otimes A)\otimes (H\otimes H\otimes ((c_{H,A}\otimes H)\circ (H\otimes c_{H,A}))\otimes A)\circ (H\otimes H\otimes H\otimes H\otimes c_{A,A})$
	\item[ ]$\hspace{0.38cm}\circ (H\otimes ((c_{H,H}\otimes c_{A,H})\circ (H\otimes c_{A,H}\otimes H)\circ (\rho_{A}\otimes \delta_{H}))\otimes A)\circ (\rho_{A}\otimes H\otimes A) $ {\scriptsize (by the naturality}
	\item[ ]$\hspace{0.38cm}$ {\scriptsize  of $c$ and the condition of comodule for $A$)}
	\item [ ]$=((\psi_{A}^2\circ (H\otimes \psi_{A}^2))\otimes \Omega_{H}^{A})\circ (H\otimes H\otimes c_{A,A}\otimes H)\circ (H\otimes c_{A,H}\otimes c_{H,A})\circ (\rho_{A}\otimes \delta_{H}\otimes A) $ {\scriptsize (by the}
	\item[ ]$\hspace{0.38cm}$ {\scriptsize  naturality of $c$ and the condition of left module for $A$)}
	\item [ ]$=(A\otimes \Omega_{H}^{A})\circ ( t_{A,A}\otimes H)\circ (A\otimes \Psi_{A}^{H_{i}}) $ {\scriptsize (by the naturality of $c$)}
\end{itemize}
and (\ref{11f}) follows by (\ref{3f}) and (\ref{7f}). Finally, (\ref{11f1}) follows by (\ref{11f}), (\ref{1f}) and the unit properties. 

Taking into account the previous equalities, we will prove the theorem. Firstly, 
let's assume that ${\mathbb A}\blacktriangleright\hspace{-0.1cm}\blacktriangleleft {\mathbb H}$ is a Hopf brace in {\sf C}. Then, (iii) of Definition \ref{H-brace} holds for  ${\mathbb A}\blacktriangleright\hspace{-0.1cm}\blacktriangleleft {\mathbb H}$, i.e.,  we have that the following equality:
\begin{equation}
\label{hbsmash}
\mu_{A\blacktriangleright\hspace{-0.1cm}\blacktriangleleft H}^{1}\co (\mu_{A\blacktriangleright\hspace{-0.1cm}\blacktriangleleft H}^{2}\ot \Gamma_{({A\blacktriangleright\hspace{-0.1cm}\blacktriangleleft H})_{1}} )\co (A\otimes H\ot c_{A\otimes H,A\otimes H}\ot A\otimes H)\co (\delta_{A\blacktriangleright\hspace{-0.1cm}\blacktriangleleft H}\ot A\otimes H\ot A\otimes H)
\end{equation}
$$=\mu_{A\blacktriangleright\hspace{-0.1cm}\blacktriangleleft H}^{2}\co (A\otimes H\ot \mu_{A\blacktriangleright\hspace{-0.1cm}\blacktriangleleft H}^{1}).$$

Then composing in (\ref{hbsmash}) with $\eta_{A}\otimes H\otimes \eta_{A}\otimes H\otimes A\otimes H$ we have
\begin{itemize}
	\item[ ]$\hspace{0.38cm}\mu_{A\blacktriangleright\hspace{-0.1cm}\blacktriangleleft H}^{1}\co (\mu_{A\blacktriangleright\hspace{-0.1cm}\blacktriangleleft H}^{2}\ot \Gamma_{({A\blacktriangleright\hspace{-0.1cm}\blacktriangleleft H})_{1}} )\co (A\otimes H\ot c_{A\otimes H,A\otimes H}\ot A\otimes H)\co (\delta_{A\blacktriangleright\hspace{-0.1cm}\blacktriangleleft H}\ot A\otimes H\ot A\otimes H)$
	\item[ ]$\hspace{0.38cm}\circ (\eta_{A}\otimes H\otimes \eta_{A}\otimes H\otimes A\otimes H)$
	\item [ ]$=(A\otimes \mu_{H}^1)\circ (\Psi_{A}^{H_{1}}\otimes H)\circ (\mu_{H}^{2}\otimes \Gamma_{({A\blacktriangleright\hspace{-0.1cm}\blacktriangleleft H})_{1}})\circ (((H\otimes c_{A,H}\otimes H)\circ (\Omega_{H}^{A}\otimes c_{H,H})\circ (\eta_{A}\otimes \delta_{H}\otimes H))$ 
	\item[ ]$\hspace{0.38cm}\otimes A\otimes H ) ${\scriptsize (by the naturality of $c$, the unit properties, (\ref{6f}) and the coalgebra morphism condition for $\eta_{A}$)}
	\item [ ]$=(A\otimes \mu_{H}^1)\circ (\Psi_{A}^{H_{1}}\otimes H)\circ (\mu_{H}^{2}\otimes (\Gamma_{({A\blacktriangleright\hspace{-0.1cm}\blacktriangleleft H})_{1}}\circ (\eta_{A}\otimes H\otimes A\otimes H)))\circ (H\otimes c_{H,H}\otimes A\otimes H)$
	\item[ ]$\hspace{0.38cm}\circ (\delta_{H}\otimes H\otimes A\otimes H) $ {\scriptsize (by the naturality of $c$ and (\ref{1f}))}
	\item [ ]$= (A\otimes \mu_{H}^1)\circ (\Psi_{A}^{H_{1}}\otimes H)\circ (\mu_{H}^{2}\otimes ((A\otimes \mu_{H}^{1})\circ (\Psi_{A}^{H_{1}}\otimes \mu_{H}^{2})\circ (\lambda_{H}^1\otimes  \Psi_{A}^{H_{2}}\otimes H)\circ (\delta_{H}\otimes A\otimes H)))$
	\item[ ]$\hspace{0.38cm}\circ (H\otimes c_{H,H}\otimes A\otimes H) \circ (\delta_{H}\otimes H\otimes A\otimes H) $ {\scriptsize (by (\ref{11f1}))}
	\item [ ]$= (A\otimes \mu_{H}^1)\circ (((A\otimes \mu_{H}^1)\circ (\Psi_{A}^{H_{1}}\otimes H)\circ (H\otimes \Psi_{A}^{H_{1}}))\otimes \mu_{H}^2)\circ (((\mu_{H}^2\otimes \lambda_{H}^{1})\circ (H\otimes c_{H,H})\circ (\delta_{H}\otimes H))$
	\item[ ]$\hspace{0.38cm}\otimes \Psi_{A}^{H_{2}}\otimes H)\circ (H\otimes c_{H,H}\otimes A\otimes H)\circ (\delta_{H}\otimes H\otimes A\otimes H)$ {\scriptsize (by the naturality of $c$, the coassociativity}
	\item[ ]$\hspace{0.38cm}$ {\scriptsize  of $\delta_{H}$ and the associativity of $\mu_{H}^1$)}
	\item [ ]$=(A\otimes \mu_{H}^{1})\circ (\Psi_{A}^{H_{1}} \otimes \mu_{H}^{2})\circ (H\otimes \Psi_{A}^{H_{2}}\otimes H)\circ (((\Gamma_{H_{1}}^{\prime}\otimes H)\circ (H\otimes c_{H,H})\circ (\delta_{H}\otimes H))\otimes A\otimes H) $ 
	\item[ ]$\hspace{0.38cm}${\scriptsize (by (\ref{8f}))}
\end{itemize}
and, on the other hand, using the unit properties
$$\mu_{A\blacktriangleright\hspace{-0.1cm}\blacktriangleleft H}^{2}\co (A\otimes H\ot \mu_{A\blacktriangleright\hspace{-0.1cm}\blacktriangleleft H}^{1})\circ (\eta_{A}\otimes H\otimes \eta_{A}\otimes H\otimes A\otimes H)$$
$$=(A\otimes \mu_{H}^{2})\circ (\Psi_{A}^{H_{2}} \otimes \mu_{H}^{1})\circ (H\otimes \Psi_{A}^{H_{1}}\otimes H)$$

Therefore (\ref{4161}) holds.

If we compose with $\eta_{A}\otimes H\otimes A\otimes \eta_{H}\otimes A\otimes \eta_{H}$ in (\ref{hbsmash}), by the unit properties and (\ref{5f}) we obtain 
$$\mu_{A\blacktriangleright\hspace{-0.1cm}\blacktriangleleft H}^{2}\co (A\otimes H\ot \mu_{A\blacktriangleright\hspace{-0.1cm}\blacktriangleleft H}^{1})\circ (\eta_{A}\otimes H\otimes A\otimes \eta_{H}\otimes A\otimes \eta_{H})$$
$$=\Psi_{A}^{H_{2}} \circ (H\otimes \mu_{A}^{1})$$
and, on the other hand, 
\begin{itemize}
	\item[ ]$\hspace{0.38cm}\mu_{A\blacktriangleright\hspace{-0.1cm}\blacktriangleleft H}^{1}\co (\mu_{A\blacktriangleright\hspace{-0.1cm}\blacktriangleleft H}^{2}\ot \Gamma_{({A\blacktriangleright\hspace{-0.1cm}\blacktriangleleft H})_{1}} )\co (A\otimes H\ot c_{A\otimes H,A\otimes H}\ot A\otimes H)\co (\delta_{A\blacktriangleright\hspace{-0.1cm}\blacktriangleleft H}\ot A\otimes H\ot A\otimes H)$
	\item[ ]$\hspace{0.38cm}\circ (\eta_{A}\otimes H\otimes A\otimes \eta_{H}\otimes A\otimes \eta_{H})$
	\item [ ]$= (\mu_{A}^{1}\otimes \mu_{H}^{1})\circ (A\otimes ((A\otimes \mu_{H}^1)\circ (\Psi_{A}^{H_{1}}\otimes H)\circ (H\otimes \Psi_{A}^{H_{1}}))\otimes H)\circ (\Psi_{A}^{H_{2}}\otimes ((\lambda_{H}^{1}\otimes \Psi_{A}^{H_{2}})\circ (\delta_{H}\otimes A)))$
	\item[ ]$\hspace{0.38cm}\circ (H\otimes c_{H,A}\otimes A)\circ (\delta_{H}\otimes A\otimes A)${\scriptsize (by the naturality of $c$, the unit properties, the coalgebra morphism}
	\item[ ]$\hspace{0.38cm}${\scriptsize  condition for $\eta_{A}$, the associativity of $\mu_{H}^1$, (\ref{1f}), (\ref{11f1}) and the unit properties)}
	\item [ ]$=(\mu_{A}^{1}\otimes \mu_{H}^{1})\circ (A\otimes (\Psi_{A}^{H_{1}}\circ ((id_{H}\ast \lambda_{H}^{1})\otimes A))\otimes H)\circ (A\otimes H\otimes \Psi_{A}^{H_{2}}) \circ (((A\otimes \delta_{H})\circ \Psi_{A}^{H_{2}})\otimes A)$ 
	\item[ ]$\hspace{0.38cm}${\scriptsize (by (\ref{8f}) and (\ref{9f}))}
	\item [ ]$=  (\mu_{A}^{1}\otimes H)\circ (A\otimes \Psi_{A}^{H_{2}})\circ (\Psi_{A}^{H_{2}}\otimes A)${\scriptsize (by the unit properties, (\ref{antipode}) and (\ref{5f}))}.
\end{itemize}

Therefore the equality (\ref{4122}) holds.

Finally, the proof for (\ref{4123}) is the following: Composing with $A\otimes \eta_{H}\otimes \eta_{A}\otimes H\otimes A\otimes \eta_{H}$ in (\ref{hbsmash}), by the unit properties and (\ref{5f1}), we obtain 
$$\mu_{A\blacktriangleright\hspace{-0.1cm}\blacktriangleleft H}^{2}\co (A\otimes H\ot \mu_{A\blacktriangleright\hspace{-0.1cm}\blacktriangleleft H}^{1})\circ (A\otimes \eta_{H}\otimes \eta_{A}\otimes H\otimes A\otimes \eta_{H})$$
$$=(\mu_{A}^{2}\otimes H)\circ (A\otimes \Psi_{A}^{H_{1}})$$
and, on the other hand, 
\begin{itemize}
	\item[ ]$\hspace{0.38cm}\mu_{A\blacktriangleright\hspace{-0.1cm}\blacktriangleleft H}^{1}\co (\mu_{A\blacktriangleright\hspace{-0.1cm}\blacktriangleleft H}^{2}\ot \Gamma_{({A\blacktriangleright\hspace{-0.1cm}\blacktriangleleft H})_{1}} )\co (A\otimes H\ot c_{A\otimes H,A\otimes H}\ot A\otimes H)\co (\delta_{A\blacktriangleright\hspace{-0.1cm}\blacktriangleleft H}\ot A\otimes H\ot A\otimes H)$
	\item[ ]$\hspace{0.38cm}\circ (A\otimes \eta_{H}\otimes \eta_{A}\otimes H\otimes A\otimes \eta_{H})$
	\item [ ]$=(\mu_{A}^{1}\otimes \mu_{H}^{1})\circ (A\otimes \Psi_{A}^{H_{1}}\otimes H)\circ (A\otimes H\otimes ((A\otimes \mu_{H}^{1})\circ ((\Psi_{A}^{H_{1}}\circ (\lambda_{H}^{1}\otimes \Gamma_{A_{1}})\circ (\Omega_{H}^{A}\otimes A))\otimes H)$ 
	\item[ ]$\hspace{0.38cm}\circ (A\otimes H\otimes \Psi_{A}^{H_{2}})\circ (A\otimes \delta_{H}\otimes A)))\circ (A\otimes ((\mu_{H}^2\otimes A\otimes H)\circ (H\otimes c_{A,H}\otimes H)\circ (\rho_{A}\otimes c_{H,H}))\otimes A)$ 
	\item[ ]$\hspace{0.38cm}\circ (\delta_{A}\otimes \eta_{H}\otimes H\otimes A) ${\scriptsize (by the naturality of $c$, (\ref{6f1}), the condition of coalgebra morphism of $\eta_{H}$, (\ref{2f}),(\ref{11f})}
	\item[ ]$\hspace{0.38cm}${\scriptsize  and the unit properties)}
	\item [ ]$=(\mu_{A}^{1}\otimes H)\circ (A\otimes ((A\otimes \mu_{H}^{1})\circ (\Psi_{A}^{H_{1}}\otimes H)\circ (H\otimes \Psi_{A}^{H_{1}})))\circ (A\otimes \mu_{H}^{2}\otimes ((\lambda_{H}^{1}\otimes \Gamma_{A_{1}})\circ (\rho_{A}\otimes A)))$
	\item[ ]$\hspace{0.38cm}\circ (A\otimes H\otimes c_{A,H}\otimes A)\circ (((A\otimes \rho_{A})\circ \delta_{A})\otimes H\otimes A) ${\scriptsize (by the naturality of $c$, the condition of coalgebra}
	\item[ ]$\hspace{0.38cm}$ {\scriptsize  morphism for $\eta_{H}$, the unit properties, (\ref{2f}) and (\ref{5f}))}
	\item [ ]$=  (\mu_{A}^{1}\otimes H)\circ (A\otimes \Psi_{A}^{H_{1}})\circ (A\otimes ((\Gamma_{H_{1}}^{\prime}\otimes \Gamma_{A_{1}})\circ (H\otimes c_{A,H}\otimes A)\circ (\rho_{A}\otimes H\otimes A)))\circ (\delta_{A}\otimes H\otimes A)${\scriptsize (by }
	\item[ ]$\hspace{0.38cm}${\scriptsize  the naturality of $c$, the comodule condition for $A$ and (\ref{8f}))}
\end{itemize}

Conversely, let's assume that (\ref{4161}), (\ref{4122}) and (\ref{4123}) hold. By the bosonization process we know that 
$$(A\blacktriangleright\hspace{-0.1cm}\blacktriangleleft H)_{1}=(A\otimes H, \eta_{A\blacktriangleright\hspace{-0.1cm}\blacktriangleleft H}, \mu_{A\blacktriangleright\hspace{-0.1cm}\blacktriangleleft H}^1, \varepsilon_{A\blacktriangleright\hspace{-0.1cm}\blacktriangleleft H}, \delta_{A\blacktriangleright\hspace{-0.1cm}\blacktriangleleft H}, \lambda_{A\blacktriangleright\hspace{-0.1cm}\blacktriangleleft H}^1)$$
and 
$$(A\blacktriangleright\hspace{-0.1cm}\blacktriangleleft H)_{2}=(A\otimes H, \eta_{A\blacktriangleright\hspace{-0.1cm}\blacktriangleleft H}, \mu_{A\blacktriangleright\hspace{-0.1cm}\blacktriangleleft H}^2, \varepsilon_{A\blacktriangleright\hspace{-0.1cm}\blacktriangleleft H}, \delta_{A\blacktriangleright\hspace{-0.1cm}\blacktriangleleft H}, \lambda_{A\blacktriangleright\hspace{-0.1cm}\blacktriangleleft H}^2)$$
are Hopf algebras in {\sf C}. Then, to finish the proof we only need to show that (iii) of Definition \ref{H-brace} holds for ${\mathbb A}\blacktriangleright\hspace{-0.1cm}\blacktriangleleft {\mathbb H}$. 
Indeed, first note that if (\ref{4123}) holds, we have that
\begin{equation}
\label{12f}
(\Gamma_{A_{1}}\otimes H)\circ (A\otimes \Psi_{A}^{H_{1}}) =
\Psi_{A}^{H_{1}}\circ (\Gamma_{H_{1}}^{\prime}\otimes \Gamma_{A_{1}})\circ (H\otimes c_{A,H}\otimes A)\circ (\rho_{A}\otimes H\otimes A)
\end{equation}
also holds because 
\begin{itemize}
	\item[ ]$\hspace{0.38cm}(\Gamma_{A_{1}}\otimes H)\circ (A\otimes \Psi_{A}^{H_{1}}) $
	\item [ ]$=(\mu_{A}^{1}\otimes H)\circ (\lambda_{A}^{1}\otimes ((\mu_{A}^{1}\otimes H)\circ (A\otimes \Psi_{A}^{H_{1}})\circ (A\otimes ((\Gamma_{H_{1}}^{\prime}\otimes \Gamma_{A_{1}})\circ (H\otimes c_{A,H}\otimes A)\circ (\rho_{A}\otimes H\otimes A)))$
	\item[ ]$\hspace{0.38cm}\circ (\delta_{A}\otimes H\otimes A)))\circ (\delta_{A}\otimes H\otimes A) $ {\scriptsize (by (\ref{4123}))}
	\item [ ]$=\Psi_{A}^{H_{1}}\circ (\Gamma_{H_{1}}^{\prime}\otimes \Gamma_{A_{1}})\circ (H\otimes c_{A,H}\otimes A)\circ (\rho_{A}\otimes H\otimes A) ${\scriptsize (by the coassociativity of $\delta_{A}$, the associativity)}
	\item[ ]$\hspace{0.38cm}$ {\scriptsize  of $\mu_{A}^{1}$, (\ref{antipode}) and the unit and counit properties)}.
\end{itemize}

Then, as a consequence of (\ref{12f}), we can prove the identity
\begin{equation}
\label{12f1}
(\Gamma_{A_{1}}\otimes H)\circ (A\otimes (\Psi_{A}^{H_{1}}\circ (\Gamma_{H_{1}}^{\prime}\otimes A))) =
\Psi_{A}^{H_{1}}\circ (\Gamma_{H_{1}}^{\prime}\otimes \Gamma_{A_{1}})\circ (H\otimes c_{A,H}\otimes A)\circ (\Omega_{H}^{A}\otimes H\otimes A)
\end{equation}
because
\begin{itemize}
	\item[ ]$\hspace{0.38cm}\Psi_{A}^{H_{1}}\circ (\Gamma_{H_{1}}^{\prime}\otimes \Gamma_{A_{1}})\circ (H\otimes c_{A,H}\otimes A)\circ (\Omega_{H}^{A}\otimes H\otimes A) $
	\item [ ]$=\Psi_{A}^{H_{1}}\circ (\Gamma_{H_{1}}^{\prime}\otimes \Gamma_{A_{1}})\circ (H\otimes ((\Gamma_{H_{1}}^{\prime}\otimes A)\circ (H\otimes c_{A,H})\circ (c_{A,H}\otimes H))\otimes A)\circ (\rho_{A}\otimes H\otimes H\otimes A) $ 
	\item[ ]$\hspace{0.38cm}${\scriptsize (by the condition of $H_{2}$-module with action $\Gamma_{H_{1}}^{\prime}$ for $H_{1}$ )}
	\item [ ]$=\Psi_{A}^{H_{1}}\circ (\Gamma_{H_{1}}^{\prime}\otimes \Gamma_{A_{1}})\circ (H\otimes c_{A,H}\otimes A)\circ (\rho_{A}\otimes \Gamma_{H_{1}}^{\prime}\otimes A) $ {\scriptsize (by the naturality of $c$)}
	\item [ ]$=(\Gamma_{A_{1}}\otimes H)\circ (A\otimes (\Psi_{A}^{H_{1}}\circ (\Gamma_{H_{1}}^{\prime}\otimes A))) $ {\scriptsize (by (\ref{12f})).}
\end{itemize}

Therefore, 
\begin{itemize}
	\item[ ]$\hspace{0.38cm} \mu_{A\blacktriangleright\hspace{-0.1cm}\blacktriangleleft H}^{1}\co (\mu_{A\blacktriangleright\hspace{-0.1cm}\blacktriangleleft H}^{2}\ot \Gamma_{({A\blacktriangleright\hspace{-0.1cm}\blacktriangleleft H})_{1}} )\co (A\otimes H\ot c_{A\otimes H,A\otimes H}\ot A\otimes H)\co (\delta_{A\blacktriangleright\hspace{-0.1cm}\blacktriangleleft H}\ot A\otimes H\ot A\otimes H)  $
	\item [ ]$=(\mu_{A}^1\otimes \mu_{H}^{1}) \circ (A\otimes ((A\otimes \mu_{H}^{1})\circ (\Psi_{A}^{H_{1}}\otimes H)\circ (H\otimes \Psi_{A}^{H_{1}}))\otimes H)\circ (\mu_{A\blacktriangleright\hspace{-0.1cm}\blacktriangleleft H}^{2}\otimes \lambda_{H}^{1}\otimes \Gamma_{A_{1}}\otimes \mu_{H}^{2})$
	\item[ ]$\hspace{0.38cm}\circ (A\otimes H\otimes c_{H\otimes A,A\otimes H}\otimes \Psi_{A}^{H_{2}}\otimes H)\circ (A\otimes ((H\otimes \Omega_{H}^{A})\circ (\Omega_{H}^{A}\otimes H)\circ (A\otimes \delta_{H}))\otimes ((A\otimes c_{H,H})$
	\item[ ]$\hspace{0.38cm}\circ (c_{H,A}\otimes H))\otimes A\otimes H)\circ (\delta_{A}\otimes \delta_{H}\otimes A\otimes H\otimes A\otimes H)
	${\scriptsize (by the naturality of $c$, coassociativity of $\delta_{H}$, }
	\item[ ]$\hspace{0.38cm}${\scriptsize  associativity of $\mu_{H}^1$ and (\ref{11f}))}
	
	\item [ ]$=(\mu_{A}^1\otimes \mu_{H}^{1}) \circ (\mu_{A}^2\otimes (\Psi_{A}^{H_{1}}\circ ((\mu_{H}^{1}\circ (\mu_{H}^{2}\otimes \lambda_{H}^{1})\circ (H\otimes c_{H,H}))\otimes \Gamma_{A_{1}}))\otimes \mu_{H}^2) \circ (A\otimes ((\Psi_{A}^{H_{2}}\otimes H)$
	\item[ ]$\hspace{0.38cm}\circ (H\otimes c_{H,A})\circ (\delta_{H}\otimes A))\otimes c_{A,H}\otimes \Psi_{A}^{H_{2}}\otimes H)\circ (A\otimes H\otimes c_{A,A}\otimes c_{H,H}\otimes A\otimes H)$
	\item[ ]$\hspace{0.38cm}\circ (A\otimes \Omega_{H}^{A}\otimes c_{H,A}\otimes H\otimes A\otimes H )\circ (\delta_{A}\otimes \delta_{H}\otimes A\otimes H\otimes A\otimes H)$ {\scriptsize (by (\ref{4f}) and (\ref{8f}))}
	
	\item [ ]$= (\mu_{A}^1\otimes \mu_{H}^{1}) \circ (\mu_{A}^2\otimes (\Psi_{A}^{H_{1}}\circ (\Gamma_{H_{1}}^{\prime}\otimes \Gamma_{A_{1}})\circ (H\otimes c_{A,H}\otimes A))\otimes \mu_{H}^2) \circ (A\otimes ((\Psi_{A}^{H_{2}}\otimes A)\circ (H\otimes c_{A,A})$
	\item[ ]$\hspace{0.38cm}\circ (\Omega_{H}^{A}\otimes A))\otimes ((H\otimes \Psi_{A}^{H_{2}})\circ (c_{H,H}\otimes A))\otimes H)\circ (\delta_{A}\otimes ((H\otimes c_{H,A})\circ (\delta_{H}\otimes A))\otimes H\otimes A\otimes H)$ 
	\item[ ]$\hspace{0.38cm}${\scriptsize (by (\ref{9f}))}
	
	\item [ ]$=(\mu_{A}^1\otimes \mu_{H}^{1}) \circ (\mu_{A}^2\otimes (\Psi_{A}^{H_{1}}\circ (\Gamma_{H_{1}}^{\prime}\otimes \Gamma_{A_{1}})\circ (H\otimes c_{A,H}\otimes A)\circ (\Omega_{H}^{A}\otimes H\otimes A))\otimes \mu_{H}^2)$
	\item[ ]$\hspace{0.38cm}\circ (A\otimes t_{A,A}\otimes H\otimes ((H\otimes \Psi_{A}^{H_{2}})\circ (c_{H,H}\otimes A))\otimes H)$
	\item[ ]$\hspace{0.38cm}\circ (\delta_{A}\otimes ((\Psi_{A}^{H_{2}}\otimes H)\circ (H\otimes c_{H,A}))\circ (\delta_{H}\otimes A))\otimes H\otimes A\otimes H)$ {\scriptsize (by (\ref{10f}))}
	
    \item [ ]$=(\mu_{A}^1\otimes \mu_{H}^{1}) \circ (\mu_{A}^2\otimes ((\Gamma_{A_{1}}\otimes H)\circ (A\otimes (\Psi_{A}^{H_{1}}\circ (\Gamma_{H_{1}}^{\prime}\otimes A))))\otimes \mu_{H}^2)$
    \item[ ]$\hspace{0.38cm}\circ (A\otimes t_{A,A}\otimes H\otimes ((H\otimes \Psi_{A}^{H_{2}})\circ (c_{H,H}\otimes A))\otimes H)\circ (\delta_{A}\otimes ((A\otimes \delta_{H})\circ \Psi_{A}^{H_{2}})\otimes H\otimes A\otimes H)$ 
    \item[ ]$\hspace{0.38cm}${\scriptsize (by (\ref{9f}) and (\ref{12f1}))}
    \item [ ]$=(\mu_{A}^2\otimes H) \circ (A\otimes ((\mu_{A}^1\otimes \mu_{H}^2)\circ (A\otimes \Psi_{A}^{H_{2}}\otimes \mu_{H}^{1})\circ (\Psi_{A}^{H_{2}}\otimes\Psi_{A}^{H_{1}}\otimes H))) $ {\scriptsize (by the (iii) of Definition}
    \item[ ]$\hspace{0.38cm}$ {\scriptsize  \ref{H-brace} for ${\mathbb A}$ and (\ref{4161}))}
	\item [ ]$=\mu_{ A\blacktriangleright\hspace{-0.1cm}\blacktriangleleft  H}^{2}\co (A\otimes H\ot \mu_{ A\blacktriangleright\hspace{-0.1cm}\blacktriangleleft  H}^{1}) ${\scriptsize (by (\ref{4122}))}

\end{itemize}
\end{proof}

\begin{theorem}
\label{ex111}
Let's assume that ${\sf C}$ is symmetric. Let ${\mathbb H}$ be a cocommutative Hopf brace in ${\sf C}$ and let  ${\mathbb A}$ be a  Hopf brace in $\:^{\mathbb H}_{\mathbb H}{\sf YD}$ satisfying (\ref{4161}), (\ref{4122}) and (\ref{4123}). Then, $$({\mathbb H}, {\mathbb A}\blacktriangleright\hspace{-0.1cm}\blacktriangleleft {\mathbb H}, x=\eta_{A}\otimes H, y=\varepsilon_{A}\otimes H)$$ is a ${\rm v}_{1}$-strong projection of Hopf braces.
\end{theorem}
\begin{proof} First note that if ${\mathbb A}$ be a  Hopf brace in $\:^{\mathbb H}_{\mathbb H}{\sf YD}$ we have that $x$ and $y$ are Hopf brace morphisms and $y\circ x=id_{H}$. On the other hand, by (\ref{2f}), (\ref{u-antip}), (\ref{6f}), (\ref{comod-coalg}), the unit properties, the naturality of $c$ and (\ref{antipode})  we have that 
$$q^{1}_{ A\blacktriangleright\hspace{-0.1cm}\blacktriangleleft  H}=A\otimes (\eta_{H}\circ \varepsilon_{H})=q^{2}_{ A\blacktriangleright\hspace{-0.1cm}\blacktriangleleft H}.$$

Then, 
$$p^{1}_{ A\blacktriangleright\hspace{-0.1cm}\blacktriangleleft  H}=p^{2}_{ A\blacktriangleright\hspace{-0.1cm}\blacktriangleleft  H}=A\otimes \varepsilon_{H}, \;\;\; i^{1}_{ A\blacktriangleright\hspace{-0.1cm}\blacktriangleleft  H}=i^{2}_{ A\blacktriangleright\hspace{-0.1cm}\blacktriangleleft  H}=A\otimes \eta_{H}$$
and
$$I(q^{1}_{ A\blacktriangleright\hspace{-0.1cm}\blacktriangleleft  H})=I(q^{2}_{ A\blacktriangleright\hspace{-0.1cm}\blacktriangleleft  H})=A.$$

Then, we have an unique idempotent that we can denote by  $q_{ A\blacktriangleright\hspace{-0.1cm}\blacktriangleleft  H}$ and, also, with $p_{ A\blacktriangleright\hspace{-0.1cm}\blacktriangleleft  H}$ and $i_{ A\blacktriangleright\hspace{-0.1cm}\blacktriangleleft  H}$ we will denote the associated  projection and injection respectively. 

Therefore, by a routine calculus we have that 
$$\eta_{I(q_{ A\blacktriangleright\hspace{-0.1cm}\blacktriangleleft  H})}=\eta_{A}, \;\;\mu_{I(q_{A\blacktriangleright\hspace{-0.1cm}\blacktriangleleft  H})}^1=\mu_{A}^1, \;\; \mu_{I(q_{A\blacktriangleright\hspace{-0.1cm}\blacktriangleleft  H})}^2=\mu_{A}^2 $$
$$\varepsilon_{I(q_{ A\blacktriangleright\hspace{-0.1cm}\blacktriangleleft  H})}=\varepsilon_{A}, \;\;\delta_{I(q_{ A\blacktriangleright\hspace{-0.1cm}\blacktriangleleft  H})}=\delta_{A}, $$
$$\lambda_{I(q_{ A\blacktriangleright\hspace{-0.1cm}\blacktriangleleft  H})}^1=\lambda_{A}^1, \;\;\lambda_{I(q_{ A\blacktriangleright\hspace{-0.1cm}\blacktriangleleft  H})}^2=\lambda_{A}^2$$
and 
$$\psi_{I(q_{A\blacktriangleright\hspace{-0.1cm}\blacktriangleleft  H})}^1=\psi_{A}^1, \;\;\; \psi_{I(q_{A\blacktriangleright\hspace{-0.1cm}\blacktriangleleft  H})}^2=\psi_{A}^2, \;\;\;\rho_{I(q_{A\blacktriangleright\hspace{-0.1cm}\blacktriangleleft  H})}=\rho_{A}.$$

On the other hand, the condition (\ref{qdqd})  holds because $q^{1}_{ A\blacktriangleright\hspace{-0.1cm}\blacktriangleleft  H}=q^{2}_{ A\blacktriangleright\hspace{-0.1cm}\blacktriangleleft  H}$. 
Finally, (\ref{spro2H-n}) holds trivially and (\ref{VS2}) and the left $H_{1}$-linearity condition of the morphism defined in (\ref{morphd1}) follow from the following facts: $I(q^{2}_{ A\blacktriangleright\hspace{-0.1cm}\blacktriangleleft  H})=A$, $\rho_{I(q_{ A\blacktriangleright\hspace{-0.1cm}\blacktriangleleft  H})}=\rho_{A}$ and $(A, \psi_{A}^1, \psi_{A}^2, \rho_{A})$ is an object in $\:^{\mathbb H}_{\mathbb H}{\sf YD}$.

\end{proof}

\begin{remark}
\label{rem12f}
Note that in the conditions of Theorem \ref{mainequiv}, if (\ref{4123}) holds, we proved that  (\ref{12f}) holds. Moreover, using (\ref{eb2}), it is easy to show that if (\ref{12f}) holds we can obtain (\ref{4123}). Therefore (\ref{12f}) and (\ref{4123}) are equivalent.
\end{remark}

\begin{definition}
Let's assume that ${\sf C}$ is  symmetric. Let ${\mathbb H}$ be a cocommutative Hopf brace in ${\sf C}$ and let  ${\mathbb A}$ be a  Hopf brace in $\:^{\mathbb H}_{\mathbb H}{\sf YD}$.	We will say that ${\mathbb A}$ is bosonizable if satisfies (\ref{4161}), (\ref{4122}) and (\ref{4123}).

These  Hopf braces  with morphisms of  Hopf braces in $\:^{\mathbb H}_{\mathbb H}{\sf YD}$ form a category that we will denote by {\sf B-HBr}$(^{\mathbb H}_{\mathbb H}{\sf YD})$.
\end{definition}

\begin{definition}
\label{v2strong}
Let $({\mathbb H}, {\mathbb D}, x,y)$ be a ${\rm v}_{1}$-strong projection of Hopf braces in {\sf C}. We will say that it is ${\rm v}_{2}$-strong if $\lambda_{I(q_{D})}^1$ is a morphism of left  $H_{2}$-modules, $\mu_{I(q_{D})}^2$ and $\lambda_{I(q_{D})}^2$ are morphisms of left  $H_{1}$-modules and the following equalities
\begin{equation}
\label{VS21-q} \mu_{I(q_{D})}^1\circ \delta_{I(q_{D})}=
(\mu_{I(q_{D})}^1\otimes \mu_{I(q_{D})}^1)\circ (I(q_{D})\otimes t_{I(q_{D}), I(q_{D})}\otimes I(q_{D}))\circ (\delta_{I(q_{D})}\otimes \delta_{I(q_{D})}),
\end{equation}
 \begin{equation}
\label{VS21} p_{D}^1\circ \mu_{D}^{1}\circ (\alpha_{D}\otimes \mu_{D}^2)\circ (D\otimes c_{D,D}\otimes D)\circ ((\delta_{D}\circ i_{D})\otimes i_{D}\otimes i_{D})=\mu_{I(q_{D})}^2\circ (I(q_{D})\otimes \mu_{I(q_{D})}^1)
\end{equation}
hold, where 
$$\alpha_{D}=\mu_{D}^{1}\circ ((q_{D}^2\circ\mu_{D}^{2})\otimes r_{D})\circ (D\otimes c_{D,D})\otimes (\delta_{D}\otimes D)$$
and 
$$r_{D}=q^1_{D}\circ \mu_{D}^1\circ ((x\circ y)\otimes \lambda_{D}^{1})\circ \delta_{D}\circ q^1_{D}.$$
	
These  projections with morphisms of projections of Hopf braces form a category that we will denote by ${\sf V}_{2}${\sf SP(HBr)}, i.e., ${\sf V}_{2}${\sf SP(HBr)} is the full subcategory of ${\sf V}_{1}${\sf SP(HBr)}  whose objects are ${\rm v}_{2}$-strong projections.
\end{definition}

\begin{remark}
If $r_{D}$ is the morphism introduced in the previous definition, it's easy to show that 
$$r_{D}=\mu_{D}^1\circ ((x\circ y)\otimes \lambda_{D}^{1})\circ \delta_{D}.$$
\end{remark}

\begin{remark}
In the conditions of the previous definition we have that 
$$(\mu_{I(q_{D})}^1\otimes \mu_{I(q_{D})}^1)\circ (I(q_{D})\otimes t_{I(q_{D}), I(q_{D})}^1\otimes I(q_{D}))\circ (\delta_{I(q_{D})}\otimes \delta_{I(q_{D})})$$
$$=(\mu_{I(q_{D})}^1\otimes \mu_{I(q_{D})}^1)\circ (I(q_{D})\otimes t_{I(q_{D}), I(q_{D})}^2\otimes I(q_{D}))\circ (\delta_{I(q_{D})}\otimes \delta_{I(q_{D})}),$$
because $I(q_{D})$ is a Hopf algebra in the category of left Yetter-Drinfeld modules over $H_1$ and by $t_{I(q_{D}), I(q_{D})}^2=t_{I(q_{D}), I(q_{D})}$.
\end{remark}

\begin{theorem}
\label{ex113}
Let's assume that ${\sf C}$ is symmetric. Let ${\mathbb H}$ be a cocommutative Hopf brace in ${\sf C}$ and let  ${\mathbb A}$ be a  bosonizable Hopf brace in $\:^{\mathbb H}_{\mathbb H}{\sf YD}$. Then, $$({\mathbb H}, {\mathbb A}\blacktriangleright\hspace{-0.1cm}\blacktriangleleft {\mathbb H}, x=\eta_{A}\otimes H, y=\varepsilon_{A}\otimes H)$$ is a ${\rm v}_{2}$-strong projection of Hopf braces.
\end{theorem}

\begin{proof}Note that, by Theorem \ref{ex111} we have that $({\mathbb H}, {\mathbb A}\blacktriangleright\hspace{-0.1cm}\blacktriangleleft {\mathbb H}, x=\eta_{A}\otimes H, y=\varepsilon_{A}\otimes H)$ is a ${\rm v}_{1}$-strong projection of Hopf braces.
On the other hand, (\ref{VS21-q}) holds because $(A, \eta_{A}, \mu_{A}^{1}, \varepsilon_{A}, \delta_{A})$ is a bialgebra in $\:^{\mathbb H}_{\mathbb H}{\sf YD}$.
Also, using the properties of  $\varepsilon_{H}$ and the naturality of $c$, we have the identity
\begin{equation}
\label{1ex113}
(A\otimes \varepsilon_{H})\circ \mu_{A\blacktriangleright\hspace{-0.1cm}\blacktriangleleft  H}^i=(\mu_{A}^i\circ (A\otimes \psi_{A}^i))\otimes \varepsilon_{H}, \;\; i=1,2.
\end{equation}
On the other hand, by the coassociativity of $\delta_{A}$ and the condition of left $H$-comodule coalgebra of $A$, we obtain that 
\begin{equation}
\label{2ex113}
(A\otimes \Omega_{H}^{A}\otimes A)\circ (\delta_{A}\otimes \rho_{A})\circ \delta_{A}=(A\otimes ((H\otimes \delta_{A})\circ \rho_{A}))\circ \delta_{A}.
\end{equation}

As a consequence of these facts we can obtain the following formulations for  the morphisms $r_{A\blacktriangleright\hspace{-0.1cm}\blacktriangleleft  H}$ and $\alpha_{A\blacktriangleright\hspace{-0.1cm}\blacktriangleleft  H}$: 
\begin{equation}
\label{3ex113}
r_{A\blacktriangleright\hspace{-0.1cm}\blacktriangleleft  H}=\lambda_{A}^1\otimes (\eta_{H}\circ\varepsilon_{H}), 
\end{equation}
\begin{equation}
\label{4ex113}
\alpha_{A\blacktriangleright\hspace{-0.1cm}\blacktriangleleft  H}=(\mu_{A}^1\circ (\mu_{A}^2\otimes \lambda_{A}^1)\circ (A\otimes ((\psi_{A}^2\otimes A)\circ (H\otimes c_{A,A})\circ (\Omega_{H}^{A}\otimes A)))\circ (\delta_{A}\otimes H\otimes A)\otimes (\eta_{H}\circ\varepsilon_{H}).	
\end{equation}

Indeed, (\ref{3ex113}) follows by 
\begin{itemize}
	\item[ ]$\hspace{0.38cm}r_{A\blacktriangleright\hspace{-0.1cm}\blacktriangleleft  H} $
	\item [ ]$=(A\otimes (\eta_{H}\circ\varepsilon_{H}))\circ \Psi_{A}^{H_{1}}\circ ((id_{H}\ast \lambda_{H}^{1})\otimes \lambda_{A}^{1})\circ \Omega_{H}^{A}\circ  (A\otimes (\eta_{H}\circ\varepsilon_{H}))$ {\scriptsize (by the unit and counit }
	\item[ ]$\hspace{0.38cm}$ {\scriptsize  properties, (\ref{4f}) and (\ref{8f}))}
	\item [ ]$=(A\otimes (\eta_{H}\circ\varepsilon_{H}))\circ \Psi_{A}^{H_{1}}\circ ((\eta_{H}\circ\varepsilon_{H})\otimes \lambda_{A}^{1})\circ \Omega_{H}^{A}\circ  (A\otimes (\eta_{H}\circ\varepsilon_{H})) ${\scriptsize (by  (\ref{antipode}))}
	\item [ ]$=((\varepsilon_{H}\otimes \lambda_{A}^{1})\circ \rho_{A})\otimes (\eta_{H}\circ\varepsilon_{H})$ {\scriptsize (by the condition of algebra morphism for $\varepsilon_{H}$, (\ref{2f}) and (\ref{5f}))}
	\item [ ]$=\lambda_{A}^1\otimes (\eta_{H}\circ\varepsilon_{H})$ {\scriptsize (by the comodule condition for $A$)}
\end{itemize}
and, by (\ref{1ex113}), the unit and counit properties and  (\ref{5f1}), we obtain (\ref{4ex113}).

Then, 
\begin{itemize}
	\item[ ]$\hspace{0.38cm}p_{A\blacktriangleright\hspace{-0.1cm}\blacktriangleleft  H}^1\circ \mu_{A\blacktriangleright\hspace{-0.1cm}\blacktriangleleft  H}^{1}\circ (\alpha_{A\blacktriangleright\hspace{-0.1cm}\blacktriangleleft  H}\otimes \mu_{A\blacktriangleright\hspace{-0.1cm}\blacktriangleleft  H}^2)\circ (A\otimes H\otimes c_{A\otimes H,A\otimes H}\otimes A\otimes H)$
	\item[ ]$\hspace{0.38cm}\circ ((\delta_{A\blacktriangleright\hspace{-0.1cm}\blacktriangleleft  H}\circ i_{A\blacktriangleright\hspace{-0.1cm}\blacktriangleleft  H}^1)\otimes i_{A\blacktriangleright\hspace{-0.1cm}\blacktriangleleft  H}^1\otimes i_{A\blacktriangleright\hspace{-0.1cm}\blacktriangleleft  H}^1) $
	\item [ ]$=\mu_{A}^1 \circ ((\mu_{A}^1\circ ((\mu_{A}^2\circ (A\otimes \psi_{A}^2))\otimes \lambda_{A}^1))\otimes \mu_{A}^2)\circ (A\otimes H\otimes ((c_{A,A}\otimes A)\circ (A\otimes c_{A,A}))\otimes A)$
	\item[ ]$\hspace{0.38cm}\circ (((A\otimes \Omega_{H}^{A}\otimes A)\circ (\delta_{A}\otimes \rho_{A})\circ \delta_{A})\otimes A\otimes A) $ {\scriptsize (by the condition of algebra morphism for $\varepsilon_{H}$, (\ref{2f}) and}
	\item[ ]$\hspace{0.38cm}${\scriptsize  (\ref{5f}))} 
	\item [ ]$=\mu_{A}^1 \circ (\mu_{A}^2 \otimes \Gamma_{A_{1}})\circ (A\otimes t_{A,A}\otimes A)\circ (\delta_{A}\otimes A\otimes A) $ {\scriptsize (by  (\ref{2ex113}))}
	\item [ ]$=\mu_{A}^2\circ (A\otimes \mu_{A}^{1}) $ {\scriptsize (by (iii) of Definition \ref{H-brace} for $A$)}
	\item [ ]$=\mu_{I(q^{1}_{ A\blacktriangleright\hspace{-0.1cm}\blacktriangleleft  H})}^{2}\circ (I(q^{1}_{ A\blacktriangleright\hspace{-0.1cm}\blacktriangleleft  H})\otimes \mu_{I(q^{1}_{ A\blacktriangleright\hspace{-0.1cm}\blacktriangleleft  H})}^{1}) $ {\scriptsize (by the identities of the proof of Theorem \ref{ex111})}
\end{itemize}
and, as a consequence, $({\mathbb H}, {\mathbb A}\blacktriangleright\hspace{-0.1cm}\blacktriangleleft {\mathbb H}, x=\eta_{A}\otimes H, y=\varepsilon_{A}\otimes H)$ is a ${\rm v}_{2}$-strong projection of Hopf braces. Note that, $\lambda_{A}^1$ is a morphism of left  $H_{2}$-modules, $\mu_{A}^2$ and $\lambda_{A}^2$ are morphisms of left  $H_{1}$-modules because ${\mathbb A}$ is a Hopf brace in $\:^{\mathbb H}_{\mathbb H}{\sf YD}$.
\end{proof}

\begin{theorem}
\label{311YD}
Let's assume that ${\sf C}$ is  symmetric. Let ${\mathbb H}$ be a cocommutative Hopf brace in ${\sf C}$. Let $({\mathbb H}, {\mathbb D}, x,y)$ be a ${\rm v}_{2}$-strong projection of Hopf braces. Then,
$$
{\mathbb I(q_{D})}=(I(q_{D}), \eta_{I(q_{D})}, \mu_{I(q_{D})}^{1},  \mu_{I(q_{D})}^{2}, \varepsilon_{I(q_{D})}, \delta_{I(q_{D})}, \lambda_{I(q_{D})}^{1}, \lambda_{I(q_{D})}^{2})$$
is a Hopf brace in $\:^{\mathbb H}_{\mathbb H}{\sf YD}$, where $\eta_{I(q_{D})}$ is defined as in (\ref{3111}), $\mu_{I(q_{D})}^{1}$ as in (\ref{3112}), $\mu_{I(q_{D})}^{2}$ as in (\ref{3113}), $\varepsilon_{I(q_{D})}$ as in (\ref{31131}), $\delta_{I(q_{D})}$ as in (\ref{3114}), 
$$
%\label{31151}
\lambda_{I(q_{D})}^{1}= \psi_{I(q_{D})}^1\circ (H\otimes (p_{D}^{1}\circ \lambda_{D}^1\circ i_{D}))\circ \rho_{I(q_{D})},
$$
and
$$
%\label{31161}
\lambda_{I(q_{D})}^{2}= \psi_{I(q_{D})}^2\circ (H\otimes (p_{D}^{2}\circ \lambda_{D}^2\circ i_{D}))\circ \rho_{I(q_{D})}, 
$$
being $\psi_{I(q_{D})}^1$, $\psi_{I(q_{D})}^2$ and $\rho_{I(q_{D})}$ the actions and the coaction introduced in Theorem \ref{YDbrace3}.

\end{theorem}

\begin{proof} By Theorem \ref{YDbrace3} we know that the triple $$ (I(q_{D}),\psi_{I(q_{D})}^1=p_{D}^1\circ \mu_{D}^1\circ (x\ot i_{D}), \psi_{I(q_{D})}^2=p_{D}^2\circ\mu_{D}^2\circ (x\ot i_{D}), \rho_{I(q_{D})}=(y\otimes p_{D}^{1})\circ \delta_{D}\circ i_{D})$$
is an object in $\;^{\mathbb H}_{\mathbb H}{\sf YD}.$ Also, by the theory of Hopf algebra projections, 
$$(I(q_{D}), \eta_{I(q_{D})}, \mu_{I(q_{D})}^{1},  \varepsilon_{I(q_{D})}, \delta_{I(q_{D})}, \lambda_{I(q_{D})}^{1})$$ is a Hopf algebra in $\;^{H_{1}}_{H_{1}}{\sf YD}$  and $$ (I(q_{D}), \eta_{I(q_{D})}, \mu_{I(q_{D})}^{2},  \varepsilon_{I(q_{D})}, \delta_{I(q_{D})}, \lambda_{I(q_{D})}^{2})$$ is a Hopf algebra in $\;^{H_{2}}_{H_{2}}{\sf YD}$. 
Moreover, by Theorem \ref{th-proj-mod1} and the conditions of the theorem, we know that  $\eta_{I(q_{D})},\;$ $\mu_{I(q_{D})}^{1},\;$  $\mu_{I(q_{D})}^{2},\;$ $ \varepsilon_{I(q_{D})},$ $ \delta_{I(q_{D})},\;$  $\lambda_{I(q_{D})}^{1}$ and  $\lambda_{I(q_{D})}^{2}$ are morphisms in $\:^{\mathbb H}_{\mathbb H}{\sf YD}.$ Therefore, by (\ref{VS21-q}), $$(I(q_{D}), \eta_{I(q_{D})}, \mu_{I(q_{D})}^{1},  \varepsilon_{I(q_{D})}, \delta_{I(q_{D})}, \lambda_{I(q_{D})}^{1}),\;\;
(I(q_{D}), \eta_{I(q_{D})}, \mu_{I(q_{D})}^{2},  \varepsilon_{I(q_{D})}, \delta_{I(q_{D})}, \lambda_{I(q_{D})}^{2})$$ are Hopf algebras in $\:^{\mathbb H}_{\mathbb H}{\sf YD}.$

Then, to finish the proof we only need to check that (iii) of Definition \ref{H-brace} holds for 
${\mathbb I(q_{D})}$ in $\:^{\mathbb H}_{\mathbb H}{\sf YD}.$ Indeed, first note that using  
the coalgebra morphism condition for $y$, the algebra morphism condition for $x$, the associativity of $\mu_{D}^{i}$, the coassociativity of $\delta_{D}$, (\ref{antipode}) and the unit and counit properties, we obtain that 
\begin{equation}
\label{auxf}
q_{D}^i\ast (x\circ y)=id_{D}, \;\; i=1,2.
\end{equation}

Then, 
\begin{itemize}
	\item[ ]$\hspace{0.38cm}\mu_{I(q_{D})}^{1}\circ ( \mu_{I(q_{D})}^{2}\otimes \Gamma_{{I(q_{D})}_{1}})\circ ({I(q_{D})}\otimes t_{{I(q_{D})},{I(q_{D})}}\otimes {I(q_{D})})\circ (\delta_{I(q_{D})}\otimes I(q_{D})\otimes I(q_{D}))$ 
	\item [ ]$= p_{D}^{1}\circ\mu_{D}^{1}\circ (D\otimes (\mu_{D}^{1}\circ (r_{D}\otimes \mu_{D}^{2})\circ (\delta_{D}\otimes D)))\circ (((q_{D}^{2}\circ \mu_{D}^2)\otimes D)\circ (q_{D}^{2}\otimes ((\mu_{D}^2\otimes D)\circ ((x\circ y)\otimes c_{D,D}))
	$
	\item[ ]$\hspace{0.38cm}\circ (D\otimes\delta_{D}\otimes D))\otimes D)\circ ((\delta_{D}\circ i_{D})\otimes i_{D}\otimes i_{D})$ {\scriptsize (by (\ref{alg-id}), (\ref{m-i}), (\ref{coalg-id}), (\ref{id1}) and (\ref{qdqdn}))}
	\item [ ]$=p_{D}^{1}\circ\mu_{D}^{1}\circ (D\otimes (\mu_{D}^{1}\circ (r_{D}\otimes \mu_{D}^{2})\circ (\delta_{D}\otimes D)))\circ ((((q_{D}^{2}\circ \mu_{D}^2)\otimes D)\circ ((q_{D}^2\ast (x\circ y))\otimes c_{D,D})$
	\item[ ]$\hspace{0.38cm}\circ (\delta_{D}\otimes D))\otimes D)\circ (i_{D}\otimes i_{D}\otimes i_{D}) $ {\scriptsize (by the associativity of $\mu_{D}^2$ and the coassociativity of $\delta_{D}$)}
	\item [ ]$= p_{D}^1\circ \mu_{D}^{1}\circ (\alpha_{D}\otimes \mu_{D}^2)\circ (D\otimes c_{D,D}\otimes D)\circ ((\delta_{D}\circ i_{D})\otimes i_{D}\otimes i_{D})$ {\scriptsize (by (\ref{auxf}),  the naturality of $c$, the }
	\item[ ]$\hspace{0.38cm}$ {\scriptsize  associativity of $\mu_{D}^1$ and the coassociativity of $\delta_{D}$)} 
	\item [ ]$=\mu_{I(q_{D})}^2\circ (I(q_{D})\otimes \mu_{I(q_{D})}^1) $ {\scriptsize (by (\ref{VS21}))}
\end{itemize}
and, as a consequence, ${\mathbb I(q_{D})}$ is a Hopf brace in $\:^{\mathbb H}_{\mathbb H}{\sf YD}$.
\end{proof}

\begin{definition}
%\label{v3strong}
Let $({\mathbb H}, {\mathbb D}, x,y)$ be a ${\rm v}_{2}$-strong projection of Hopf braces in {\sf C}. We will say that it is ${\rm v}_{3}$-strong if the following equality 
\begin{equation}
\label{VS22} ((p_{D}^1\circ \mu_{D}^{1})\otimes H)\circ (D\otimes c_{H,D})\circ (\gamma_{D\otimes H}\otimes \mu_{D}^2)\circ (D\otimes c_{D,H}\otimes D)\circ ((\delta_{D}\circ i_{D})\otimes H\otimes i_{D})
\end{equation}
$$=((p_{D}^1\circ \beta_{D})\otimes H))\circ (i_{D}\otimes (((q_{D}^1\circ \mu_{D}^1\circ (x\otimes D))\otimes H )\circ (H\otimes c_{H,D})\circ (\delta_{H}\otimes i_{D})))$$
holds, where 
$$\gamma_{D\otimes H}=(\mu_{D}^1\otimes H)\circ (x\otimes c_{H,D})\circ ((\delta_{H}\circ \Gamma_{H_{1}}^{\prime})\otimes r_{D})\circ (y\otimes c_{D,H})\circ (\delta_{D}\otimes H),$$
$$\beta_{D}=\mu_{D}^1\circ (r_{D}\otimes \mu_{D}^2)\circ (\delta_{D}\otimes D)$$
and $r_{D}$ is the morphism introduced in Definition \ref{v2strong}.
	
These  projections with morphisms of projections of Hopf braces form a category that we will denote by ${\sf V}_{3}${\sf SP(HBr)}, i.e., ${\sf V}_{3}${\sf SP(HBr)} is the full subcategory of ${\sf V}_{2}${\sf SP(HBr)}  whose objects are ${\rm v}_{3}$-strong projections.
\end{definition}

\begin{theorem}
\label{ex114}
Let's assume that ${\sf C}$ is symmetric. Let ${\mathbb H}$ be a cocommutative Hopf brace in ${\sf C}$ and let  ${\mathbb A}$ be a  bosonizable Hopf brace in $\:^{\mathbb H}_{\mathbb H}{\sf YD}$. Then, $$({\mathbb H}, {\mathbb A}\blacktriangleright\hspace{-0.1cm}\blacktriangleleft {\mathbb H}, x=\eta_{A}\otimes H, y=\varepsilon_{A}\otimes H)$$ is a ${\rm v}_{3}$-strong projection of Hopf braces.
\end{theorem}

\begin{proof} By Theorem \ref{ex113} we only need to show that (\ref{VS22}) holds. First note that, by the unit and counit properties, the naturality of $c$ and (\ref{9f}) we have that 
\begin{equation} 
\label{AHH}
\gamma_{(A\blacktriangleright\hspace{-0.1cm}\blacktriangleleft H)\otimes H}=(A\otimes \delta_{H})\circ \Psi_{A}^{H_{1}}\circ (\Gamma_{H_{1}}^{\prime}\otimes \lambda_{A}^{1})\circ (H\otimes c_{A,H})\circ (\Omega_{H}^{A}\otimes H),
\end{equation}

\begin{equation} 
\label{AHH2}
\beta_{A\blacktriangleright\hspace{-0.1cm}\blacktriangleleft H}=(\Gamma_{A_{1}}\otimes \mu_{H}^2)\circ (A\otimes \Psi_{A}^{H_{2}}\otimes H ).
\end{equation}

Then, 
\begin{itemize}
	\item[ ]$\hspace{0.38cm} ((p_{A\blacktriangleright\hspace{-0.1cm}\blacktriangleleft H}^1\circ \mu_{A\blacktriangleright\hspace{-0.1cm}\blacktriangleleft H}^{1})\otimes H)\circ (A\otimes H\otimes c_{H,A\otimes H})\circ (\gamma_{(A\blacktriangleright\hspace{-0.1cm}\blacktriangleleft H)\otimes H}\otimes \mu_{A\blacktriangleright\hspace{-0.1cm}\blacktriangleleft H}^2)\circ (A\otimes H\otimes c_{A\otimes H,H}\otimes A\otimes H)$
	\item[ ]$\hspace{0.38cm}\circ ((\delta_{A\blacktriangleright\hspace{-0.1cm}\blacktriangleleft H}\circ i_{A\blacktriangleright\hspace{-0.1cm}\blacktriangleleft H}^1)\otimes H\otimes i_{A\blacktriangleright\hspace{-0.1cm}\blacktriangleleft H}^{1})$
	\item [ ]$= (\mu_{A}^1\otimes H)\circ (A\otimes \Psi_{A}^{H_{1}})\circ (\Psi_{A}^{H_{1}}\otimes A)\circ (\Gamma_{H_{1}}^{\prime}\otimes \lambda_{A}^{1}\otimes \mu_{A}^{2})\circ (H\otimes c_{A,H}\otimes A\otimes A)\circ (H\otimes A\otimes c_{A,H}\otimes A)$
	\item[ ]$\hspace{0.38cm}\circ (((\Omega_{H}^{A}\otimes A)\circ (A\otimes \rho_{A})\circ \delta_{A})\otimes H\otimes A) $ {\scriptsize (by the condition of coalgebra morphism for $\eta_{H}$, the unit and}
	\item[ ]$\hspace{0.38cm}${\scriptsize  counit properties, the naturality of $c$, (\ref{2f}) and (\ref{AHH}))}  
	\item [ ]$=\Psi_{A}^{H_{1}}\circ (\Gamma_{H_{1}}^{\prime}\otimes \Gamma_{A_{1}})\circ (H\otimes c_{A,H}\otimes A)\circ (\rho_{A}\otimes H\otimes A) $ {\scriptsize (by the condition of comodule coalgebra, the}
	\item[ ]$\hspace{0.38cm}$ {\scriptsize  naturality of $c$ and (\ref{7f}))}
	\item [ ]$=(\Gamma_{A_{1}}\otimes H)\circ (A\otimes \Psi_{A}^{H_{1}}) $ {\scriptsize (by (\ref{12f})))}
	\item [ ]$= ((p_{A\blacktriangleright\hspace{-0.1cm}\blacktriangleleft H}^1\circ \beta_{A\blacktriangleright\hspace{-0.1cm}\blacktriangleleft H})\otimes H))\circ (i^{1}_{A\blacktriangleright\hspace{-0.1cm}\blacktriangleleft H}\otimes (((q_{A\blacktriangleright\hspace{-0.1cm}\blacktriangleleft H}^1\circ \mu_{A\blacktriangleright\hspace{-0.1cm}\blacktriangleleft H}^1\circ (x\otimes A\otimes H))\otimes H )\circ (H\otimes c_{H,A\otimes H})$
	\item[ ]$\hspace{0.38cm}\circ (\delta_{H}\otimes i^{1}_{A\blacktriangleright\hspace{-0.1cm}\blacktriangleleft H})))$ {\scriptsize (by (\ref{AHH2}), (\ref{3ex113}), the condition of algebra morphism for $\varepsilon_{H}$, the condition of coalgebra morphism}
	\item[ ]$\hspace{0.38cm}$ {\scriptsize for $\eta_{H}$, the unit and counit properties, the naturality of $c$ and (\ref{2f})).}
\end{itemize}

Therefore, $({\mathbb H}, {\mathbb A}\blacktriangleright\hspace{-0.1cm}\blacktriangleleft {\mathbb H}, x=\eta_{A}\otimes H, y=\varepsilon_{A}\otimes H)$ is a ${\rm v}_{3}$-strong projection of Hopf braces.
\end{proof}

\begin{theorem}
\label{312YD}
Let's assume that ${\sf C}$ is  symmetric. Let ${\mathbb H}$ be a cocommutative Hopf brace in ${\sf C}$.  Let $({\mathbb H}, {\mathbb D}, x,y)$ be a ${\rm v}_{3}$-strong projection of Hopf braces. Then, the Hopf brace  ${\mathbb I(q_{D})}$, introduced in Theorem \ref{311YD}, is bosonizable.
\end{theorem}

\begin{proof} By Theorem \ref{311YD} we know that ${\mathbb I(q_{D})}$ is a Hopf brace in $\:^{\mathbb H}_{\mathbb H}{\sf YD}$. Also, by Theorem \ref{th-proj-mod1} we have that (\ref{4161}) and (\ref{4122}) hold. Then to finish the proof, by Remark \ref{rem12f}, we only need to show that (\ref{12f}) holds. Indeed, 
\begin{itemize}
	\item[ ]$\hspace{0.38cm} \Psi_{I(q_{D})}^{H_{1}}\circ (\Gamma_{H_{1}}^{\prime}\otimes \Gamma_{I(q_{D})_{1}})\circ (H\otimes c_{I(q_{D}),H}\otimes I(q_{D}))\circ (\rho_{I(q_{D})}\otimes H\otimes I(q_{D}))$
	\item [ ]$=((p_{D}^1\circ \mu_{D}^1)\otimes H)\circ (D\otimes c_{H,D}) \circ (((x\otimes H)\circ \delta_{H}\circ \Gamma_{H_{1}}^{\prime}\circ (y\otimes H))\otimes \beta_{D}) \circ (D\otimes c_{D,H}\otimes D)$
	\item[ ]$\hspace{0.38cm}\circ ((\delta_{D}\circ i_{D})\otimes H\otimes i_{D})$ {\scriptsize (by  (\ref{id1}) and (\ref{m-i}))}  
	\item [ ]$=((p_{D}^1\circ \mu_{D}^1)\otimes H)\circ (D\otimes c_{H,D}) \circ (\gamma_{D\otimes H}\otimes\mu_{D}^2) \circ (D\otimes c_{D,H}\otimes D) \circ ((\delta_{D}\circ i_{D})\otimes H\otimes i_{D})$ {\scriptsize (by  }
	\item[ ]$\hspace{0.38cm}$ {\scriptsize   the  naturality of $c$, the associativity of $\mu_{D}^1$  and coassociativity of $\delta_{D}$)}
	\item [ ]$=((p_{D}^1\circ \beta_{D})\otimes H))\circ (i^{1}_{D}\otimes (((q_{D}^1\circ \mu_{D}^1\circ (x\otimes D))\otimes H )\circ (H\otimes c_{H,D})\circ (\delta_{H}\otimes i^{1}_{D}))) $ {\scriptsize (by (\ref{VS22})))}
	\item [ ]$= (\Gamma_{I(q_{D})_{1}}\otimes H)\circ (I(q_{D})\otimes \Psi_{I(q_{D})}^{H_{1}}) $ {\scriptsize (by  (\ref{id1}) and (\ref{m-i})).}
\end{itemize}

Thus, ${\mathbb I(q_{D})}$ is bosonizable.

\end{proof}

\begin{definition}
%\label{v4strong} 
Let $({\mathbb H}, {\mathbb D}, x,y)$ be a ${\rm v}_{3}$-strong projection of Hopf braces in {\sf C}. We will say that it is ${\rm v}_{4}$-strong if $q_{D}^{1}=q_{D}^2$.

These  projections with morphisms of projections of Hopf braces form a category that we will denote by ${\sf V}_{4}${\sf SP(HBr)}, i.e., ${\sf V}_{4}${\sf SP(HBr)} is the full subcategory of ${\sf V}_{3}${\sf SP(HBr)}  whose objects are ${\rm v}_{4}$-strong projections. With ${\mathbb H}\rule[2.5pt]{0.1cm}{1pt}{\sf V}_{4}${\sf SP(HBr)}  we will denote the  subcategory of ${\sf V}_{4}${\sf SP(HBr)} whose objects are ${\rm v}_{4}$-strong projections with ${\mathbb H}$ fixed and whose morphisms are the ones with the first component equal to the identity of $H$.
\end{definition}

\begin{theorem}
\label{ex115}
Let's assume that ${\sf C}$ is symmetric. Let ${\mathbb H}$ be a cocommutative Hopf brace in ${\sf C}$ and let  ${\mathbb A}$ be a  bosonizable Hopf brace in $\:^{\mathbb H}_{\mathbb H}{\sf YD}$. Then, $$({\mathbb H}, {\mathbb A}\blacktriangleright\hspace{-0.1cm}\blacktriangleleft {\mathbb H}, x=\eta_{A}\otimes H, y=\varepsilon_{A}\otimes H)$$ is a ${\rm v}_{4}$-strong projection of Hopf braces.
\end{theorem}

\begin{proof}
The proof follows by the identities of the proof of Theorem \ref{ex111}. 
\end{proof}

\begin{theorem}
\label{313YD}
Let's assume that ${\sf C}$ is  symmetric. Let ${\mathbb H}$ be a cocommutative Hopf brace in ${\sf C}$. Let $({\mathbb H}, {\mathbb D}, x,y)$ be a ${\rm v}_{4}$-strong projection of Hopf braces. Then, the Hopf brace  ${\mathbb I(q_{D})}\blacktriangleright\hspace{-0.1cm}\blacktriangleleft {\mathbb H}$  is isomorphic to ${\mathbb D}$.
\end{theorem}

\begin{proof} If  $({\mathbb H}, {\mathbb D}, x,y)$ is a ${\rm v}_{4}$-strong projection of Hopf braces, we have that 
$p_{D}^{1}=p_{D}^2$. Then, by the general theory of Hopf algebra projections (see (\ref{nnu})), $\nu_{D}=\nu_{D}^{1}=(p_{D}^1\otimes y)\co \delta_{D}=( p_{D}^2\otimes y)\co \delta_{D}=\nu_{D}^{2}$ and, as a consequence, it is a Hopf algebra isomorphism between $( I(q_{D})\blacktriangleright\hspace{-0.1cm}\blacktriangleleft  H)_{1}$ and $D_{1}$ and between $( I(q_{D})\blacktriangleright\hspace{-0.1cm}\blacktriangleleft  H)_{2}$ and $D_{2}$. Therefore, $\nu_{D}$ is a Hopf brace isomorphism between  ${\mathbb I(q_{D})}\blacktriangleright\hspace{-0.1cm}\blacktriangleleft {\mathbb H}$  and ${\mathbb D}$.
\end{proof}

\begin{remark}
In the conditions of the previous theorem we have  the equality 
\begin{equation}
\label{chfin}
\mu_{D}^{1}\circ ( i_{D}\otimes x)=\mu_{D}^{2}\circ ( i_{D}\otimes x)
\end{equation} 
because $\mu_{D}^{1}\circ ( i_{D}\otimes x)$ is the inverse of $\nu_{D}^{1}$
and $\mu_{D}^{2}\circ ( i_{D}\otimes x)$ is the inverse of $\nu_{D}^{2}$.
\end{remark}

\begin{corollary}
\label{fin}
Let's assume that ${\sf C}$ is  symmetric. Let ${\mathbb H}$ be a cocommutative Hopf brace in ${\sf C}$. The categories ${\mathbb H}\rule[2.5pt]{0.1cm}{1pt}{\sf V}_{4}${\sf SP(HBr)} and {\sf B-HBr}$(^{\mathbb H}_{\mathbb H}{\sf YD})$ are equivalent.
\end{corollary}

\begin{proof} By theorems \ref{312YD} and \ref{ex115}  it is easy to show that there exists two functors 
	$${\sf F}^{coinv}: {\mathbb H}\rule[2.5pt]{0.1cm}{1pt}{\sf V}_{4}{\sf SP(HBr)}\rightarrow {\sf B \rule[2.5pt]{0.1cm}{1pt}HBr}(^{\mathbb H}_{\mathbb H}{\sf YD}), \;\;\; {\sf G}^b:  {\sf B\rule[2.5pt]{0.1cm}{1pt}HBr}(^{\mathbb H}_{\mathbb H}{\sf YD})\rightarrow  {\mathbb H}\rule[2.5pt]{0.1cm}{1pt}{\sf V}_{4}{\sf SP(HBr)},$$ 
defined on objects by 
$$ {\sf F}^{coinv}(({\mathbb H}, {\mathbb D}, x,y))={\mathbb I(q_{D})}, \;\; {\sf G}^b({\mathbb A})=({\mathbb H}, {\mathbb A}\blacktriangleright\hspace{-0.1cm}\blacktriangleleft {\mathbb H}, x=\eta_{A}\otimes H, y=\varepsilon_{A}\otimes H)$$
and on morphisms by  the following: Let  $(id_{H} ,t): ({\mathbb H}, {\mathbb D}, x,y)\rightarrow ({\mathbb H}, {\mathbb D}^{\prime}, x^{\prime},y^{\prime})$ be a morphism in ${\mathbb H}\rule[2.5pt]{0.1cm}{1pt}{\sf V}_{4}{\sf SP(HBr)}$.  Taking into account that $q_{D}^{1}=q_{D}^2$, we will denote the idempotent morphism by $q_{D}$, the injection by $i_{D}$, the projection by $p_{D}$ and the image by $I(q_{D})$.  Define $$t_{D}:=p_{D^{\prime}}\circ t\circ  i_{D}:I(q_{D})\rightarrow I(q_{D^{\prime}})$$

Then, using that $t$ is Hopf algebra morphisms and (\ref{eqpHBr}), we have that 
\begin{center}
	\unitlength=0.75mm \special{em:linewidth 0.4pt}
	\linethickness{0.4pt}
\begin{picture}(100.00,54.00)
		\put(30.00,50.00){\vector(1,0){50.00}}
		\put(30.00,20.00){\vector(1,0){50.00}}
		\put(29.00,46.00){\vector(3,-1){19.00}}
		\put(62.00,39.00){\vector(3,1){20.00}}
		\put(29.00,16.00){\vector(3,-1){19.00}}
		\put(63.00,9.00){\vector(3,1){19.00}}
		\put(55.00,33.00){\vector(0,-1){22.00}}
		\put(86.00,47.00){\vector(0,-1){24.00}}
		\put(25.00,47.00){\vector(0,-1){25.00}}
		\put(25.00,50.00){\makebox(0,0)[cc]{$D$}}
		\put(86.00,50.00){\makebox(0,0)[cc]{$D$}}
		\put(55.00,36.00){\makebox(0,0)[cc]{$I(q_{D})$}}
		\put(25.00,19.00){\makebox(0,0)[cc]{$D^{\prime}$}}
		\put(86.00,20.00){\makebox(0,0)[cc]{$D^{\prime}$}}
		\put(55.00,54.00){\makebox(0,0)[cc]{$q_{D}$}}
		\put(43.00,45.00){\makebox(0,0)[cc]{$p_{D}$}}
		\put(67.00,45.00){\makebox(0,0)[cc]{$i_{D}$}}
		\put(22.00,34.00){\makebox(0,0)[cc]{$t$}}
		\put(60.00,26.00){\makebox(0,0)[cc]{$t_{D}$}}
		\put(90.00,35.00){\makebox(0,0)[cc]{$t$}}
		\put(43.00,24.00){\makebox(0,0)[cc]{$q_{D^{\prime}}$}}
		\put(37.00,9.00){\makebox(0,0)[cc]{$p_{D^{\prime}}$}}
		\put(73.00,8.50){\makebox(0,0)[cc]{$i_{D^{\prime}}$}}
		\put(55.00,7.00){\makebox(0,0)[cc]{$I(q_{D^{\prime}})$}}
\end{picture}
\end{center}
is a commutative diagram and, as a consequence, by a similar proof that the one used in \cite[Theorem 3.4]{NikaRamon4}, we can obtain that  $t_{D}$ is Hopf algebra morphism in $^{H_{i}}_{H_{i}}{\sf YD}$ for $i=1,2.$ Therefore, $t_{D}$ is a morphism of Hopf braces in $^{\mathbb H}_{\mathbb H}{\sf YD}$ and we define $${\sf F}^{coinv}((id_{H},t))=t_{D}.$$

On the other hand, if $s:{\mathbb A}\rightarrow {\mathbb A}^{\prime}$ is a morphism in {\sf B-HBr}$(^{\mathbb H}_{\mathbb H}{\sf YD})$, the pair $(id_{H}, s\otimes H)$ is a morphism between  $({\mathbb H}, {\mathbb A}\blacktriangleright\hspace{-0.1cm}\blacktriangleleft {\mathbb H}, x=\eta_{A}\otimes H, y=\varepsilon_{A}\otimes H)$ and $({\mathbb H}, {\mathbb A}^{\prime}\blacktriangleright\hspace{-0.1cm}\blacktriangleleft {\mathbb H}, x=\eta_{A^{\prime}}\otimes H, y=\varepsilon_{A^{\prime}}\otimes H)$. Then, we define $${\sf G}^b(s):=(id_{H}, s\otimes H).$$

Finally, following the same techniques used in the proof of \cite[Theorem 3.4]{NikaRamon4} and the isomorphism of Theorem \ref{313YD}, we can assure that  ${\sf F}^{coinv}$ and ${\sf G}^b$ induce an equivalence of categories because for all $({\mathbb H}, {\mathbb D}, x,y)$ we have that 
$$({\mathbb H}, {\mathbb D}, x,y)\simeq ({\mathbb H}, {\mathbb I(q_{D})}\blacktriangleright\hspace{-0.1cm}\blacktriangleleft {\mathbb H}, x=\eta_{I(q_{D})}\otimes H, y=\varepsilon_{I(q_{D})}\otimes H)= ({\sf G}^b\circ {\sf F}^{coiv})(({\mathbb H}, {\mathbb D}, x,y))$$
and, for all ${\mathbb A}$, $({\sf F}^{coinv}\circ {\sf G}^b)({\mathbb A})={\mathbb A}$.
\end{proof}

\begin{remark}
\label{timo}
In \cite[Theorem 5.4]{Zhu} the author works with a projection of Hopf braces $({\mathbb H}, {\mathbb D}, x,y)$ such that $I(q_{D})=I(q_{D}^2)$ and proves that there exists a Hopf brace structure on the tensor product $I(q_{D})\otimes H$ isomorphic to ${\mathbb D}$. If we study the proof in detail, we see that the author uses the identity $\nu_{D}^1\circ (\nu_{D}^{2})^{-1}=id_{I(q_{D})\otimes H}$ where $\nu_{D}^1=( p_{D}^{1}\otimes y)\circ \delta_{D}$ and $\nu_{D}^2=( p_{D}^{2}\otimes y)\circ \delta_{D}$ are the isomorphisms defined in (\ref{nnu}). Then, $\nu_{D}^1=\nu_{D}^2 $ and this implies that $q_{D}^1=q_{D}^2$ (note that this condition is not assumed in the statement of the theorem). Also, $\nu_{D}^1=\nu_{D}^2$ implies that their inverses are the same and then (\ref{chfin}) also holds. Moreover, in the statement  it is not assumed either that $\delta_{I(q_{D})\blacktriangleright\hspace{-0.1cm}\blacktriangleleft H}^1=\delta_{I(q_{D})\blacktriangleright\hspace{-0.1cm}\blacktriangleleft H}^2$ and, however, this is also used. Therefore, assuming a correct formulation of the conditions for \cite[Theorem 5.4]{Zhu},  in this theorem all that is done is to transfer the Hopf brace structure from ${\mathbb D}$ to $I(q_{D})\otimes H$ using the isomorphism $\nu_{D}$.
\end{remark}

\section*{Data Availability Statement}
Data sharing is not applicable to this article as no new data were created or analyzed in this study.

\section*{Disclosure Statement}
There are not elevant financial or non-financial competing interests.

\section*{Funding}
The  authors were supported by  Ministerio de Ciencia e Innovaci\'on of Spain. Agencia Estatal de Investigaci\'on. Uni\'on Europea - Fondo Europeo de Desarrollo Regional (FEDER). Grant PID2020-115155GB-I00: Homolog\'{\i}a, homotop\'{\i}a e invariantes categ\'oricos en grupos y \'algebras no asociativas.

\bibliographystyle{amsalpha}

\end{document}